\documentclass[a4paper,12pt,leqno]{amsart} 
\usepackage[a4paper, top=2cm, bottom=2cm, left=2cm, right=2cm]{geometry}
\usepackage[utf8]{inputenc}

\usepackage{setspace}
\usepackage{amsmath,amssymb,amsthm,mathtools}
\usepackage{caption, subcaption,cleveref}
\usepackage{amsrefs}
\usepackage{tikz, tikz-cd}
\usepackage{docmute}

\newtheorem{theorem}{Theorem}[section] 
\newtheorem{lemma}[theorem]{Lemma} 
\newtheorem{prop}[theorem]{Proposition} 
\newtheorem{cor}[theorem]{Corollary} 
\theoremstyle{definition} 
\newtheorem{rem}[theorem]{Remark} 
\newtheorem{defin}[theorem]{Definition} 
 
\newtheorem{claim}[theorem]{Claim}

\title{Systems of arcs on a torus with two punctures}
\author{Denali Relles}

\begin{document}

\maketitle

\begin{abstract}
	For a compact surface \( S \) with a finite set of marked points \( P \), we define a 1-system to be a collection of arcs which are pairwise non-homotopic and intersect pairwise at most once. We prove that, up to equivalence, there are exactly 23 maximal 1-systems on \( (S, P) \) when \( S \) is a torus and \( |P| = 2 \). 

	Along the way, we generalize some of the results of \cite{przytycki_2015} to the context of surfaces with boundary. In particular, we prove that the maximal cardinality of a 1-system on \( (S, P) \) is \( 2 |\chi| (|\chi| +  1) - \frac{v}{2} \), where \( \chi \) is the Euler characteristic of \( (S, P) \) and \( v \) is the number of marked points of \( P \) in the boundary of \( S \).
\end{abstract}

\section{Introduction}
\label{sec:intro}

Let $S$ be a compact surface, possibly with boundary, and let $P \subset S$ be a finite set of marked points. Note that marked points may be in the boundary of the surface. We define an \emph{arc on \( (S, P) \)} to be a map $u : [0,1] \to S$ such that $u(\{0,1\}) \subset P$ and $u((0,1)) \subset S \setminus P$. 

Let $u$ be an arc on \( (S, P) \).  An arc is called \textit{simple} if it is injective on $(0,1)$. An arc $u$ is called \textit{non-essential} if it can homotoped to a constant map or into \( \partial S \), relative to its endpoints, and avoiding marked points. 
Unless otherwise noted, all arcs in this paper will be simple and essential. Similarly, we call two arcs \( u, v \) \textit{homotopic as arcs} if there exists a homotopy between them, relative to their endpoints, and avoiding marked points. 

We define a \emph{$k$-system} on $(S, P)$ to be a collection of arcs \( \mathcal{A} \) such that:
\begin{itemize}
	\item 
		any pair of arcs \( u, v \in \mathcal{A} \) are not homotopic as arcs, and

	\item
		any pair of arcs \( u, v \in \mathcal{A} \) intersect at most \( k \) times.

\end{itemize}


We call two $k$-systems $\mathcal{A}, \mathcal{B}$ on \( (S, P) \) \textit{equivalent} if there exists a bijection $f: \mathcal{A} \to \mathcal{B}$ and a homeomorphism $h: S \to S$ such that $h(P) = P$, and for each $a \in \mathcal{A}$, the arc $h \circ a$ is homotopic as an arc to $f(a)$. 

We say a \( k \)-system \( \mathcal{A} \) on \( (S, P) \) is \emph{maximal} if it contains at least as many arcs as any other \( k \)-system on \( (S, P) \). 

We now state our main result.

\begin{theorem}
	Let \( S \) be a torus and let \( P = \left\{ x, y \right\} \) be any set of two points. Then there are exactly 23 equivalence classes of maximal \( 1 \)-systems on \( (S, P) \). 
	\label{thm:main thm}
\end{theorem}

\begin{figure}
	\centering
	\begin{subfigure}[b]{0.3\textwidth}
		\centering
		\begin{tikzpicture}[scale=1.6]
			\draw[dashed] (0,0) rectangle (2,2);

			\draw (1.1,1.4) -- (0,1.3) (2,1.3) -- (1.7,1.4); 
			\draw (1.7,1.4) to[out=-135,in=0] (0,1.15) (2,1.15) -- (1.7,1.4) 
			(1.7,1.4) to[out=135,in=10] (0,1.4) (2,1.4) -- (1.7,1.4); 

			\draw (1.1,1.4) -- (0.7,2) (0.7,0) -- (0,0.7) (2,0.7) -- (1.7,1.4);

			\draw (1.7,1.4) to[out=-150,in=0] (1.1,1.2) to[out=180,in=-135]  (1.3,2) (1.3,0) -- (1.7,1.4);
			\draw (1.1,1.4) -- (1.5,2) (1.5,0) -- (1.7,1.4)
			(1.7,0) -- (1.7,2);

			\draw (1.1,1.4) to[out=35,in=145] (1.7,1.4) to[out=-145,in=-35] (1.1,1.4);
			\draw[fill] (1.1,1.4) circle (0.06)
			(1.7,1.4) circle (0.06);

			\draw (1.1,1.4) -- (0.2,0) (0.2,2) -- (0,1.8) (2,1.8) -- (1.8,1.8) to[out=180,in=50] (1.1,1.4); 
			\draw (1.1,1.4) -- (0.3,0) (0.3,2) -- (0,1.6) (2,1.6) -- (1.7,1.4); 
			\draw (1.7,1.4) -- (0.4,0) (0.4,2) -- (0,1.7) (2,1.7) -- (1.8,1.7) to[out=180,in=40] (1.1,1.4); 
			\draw (1.7,1.4) -- (0.5,0) (0.5,2) -- (0,1.5) (2,1.5) -- (1.7,1.4); 

		\end{tikzpicture}
		\caption{}
	\end{subfigure}
	\begin{subfigure}[b]{0.3\textwidth}
		\centering
		\begin{tikzpicture}[scale=1.6]
			\draw[dashed] (0,0) rectangle (2,2);

			\draw (1.1,1.4) -- (1.1,2) (1.1,0) -- (1.1,1.4)
			(1.1,1.4) -- (1.3,2) (1.3,0) -- (1.7,1.4)
			(1.7,1.4) -- (1.5,2) (1.5,0) -- (1.1,1.4)
			(1.7,1.4) -- (1.7,2) (1.7,0) -- (1.7,1.4);

			\draw (1.1,1.4) -- (0,1.3) (2,1.3) -- (1.7,1.4) 
			(1.1,1.4) -- (0,1.2) (2,1.2) to[out=-170,in=-40] (1.1,1.4) 
			(1.7,1.4) to[out=135,in=10] (0,1.4) (2,1.4) -- (1.7,1.4); 

			\draw (1.1,1.4) -- (0.4,0) (0.4,2) -- (0,1.6) (2,1.6) to[out=-120,in=-45] (1.1,1.4)
			(1.1,1.4) -- (0.35,0) (0.35,2) -- (0,1.65) (2,1.65) -- (1.7,1.4)
			(1.1,1.4) -- (0.3,0) (0.3,2) -- (0,1.7) (2,1.7) to[out=-180,in=45] (1.1,1.4);

			\draw (1.1,1.4) -- (2,0.7) (0,0.7) -- (0.7,0) (0.7,2) -- (1.1,1.4);

			\draw (1.1,1.4) to[out=35,in=145] (1.7,1.4) to[out=-145,in=-35] (1.1,1.4);
			\draw[fill] (1.1,1.4) circle (0.06)
			(1.7,1.4) circle (0.06);

		\end{tikzpicture}
		\caption{}
	\end{subfigure}
	\begin{subfigure}[b]{0.3\textwidth}
		\centering
		\begin{tikzpicture}[scale=1.6]
			\draw[dashed] (0,0) rectangle (2,2);

			\draw (1.7,1.4) to[out=140,in=20] (1,1.55) to[out=-160,in=160] (1,1.25) to[out=-20,in=-140] (1.7,1.4);

			\draw (1.1,1.4) -- (0,1.3) (2,1.3) -- (1.7,1.4) 
			(1.7,1.4) to[out=-135,in=0] (0,1.15) (2,1.15) -- (1.7,1.4); 
			\draw (1.7,1.4) to[out=135,in=10] (0,1.4) (2,1.4) -- (1.7,1.4); 

			\draw (1.7,0) -- (1.7,2)
			(1.6,0) -- (1.7,1.4) (1.1,1.4) -- (1.6,2)
			(1.5,0) -- (1.7,1.4) to[out=-140,in=-20] (1,1.2) to[out=160,in=-140] (1,1.7) to[out=40,in=-150] (1.5,2);

			\draw (1.7,1.4) -- (2,1) (0,1) -- (1,0) (1,2) -- (1.1,1.4);

			\draw (1.1,1.4) to[out=45,in=180] (2,1.7) (0,1.7) -- (0.4,2) (0.4,0) -- (1.7,1.4)
			(0,1.6) -- (0.3,2) (0.3,0) -- (1.1,1.4) (1.7,1.4) -- (2,1.6);
			\draw (0,1.5) -- (0.5,2) (0.5,0) -- (1.7,1.4) -- (2,1.5);

			\draw (1.1,1.4) to[out=35,in=145] (1.7,1.4) to[out=-145,in=-35] (1.1,1.4);
			\draw[fill] (1.1,1.4) circle (0.06)
			(1.7,1.4) circle (0.06);

		\end{tikzpicture}
		\caption{}
	\end{subfigure}
	\caption{The 3 1-systems obtained when \( J \) consists of a single non-loop arc.}
	\label{fig:J12 finals}
\end{figure}
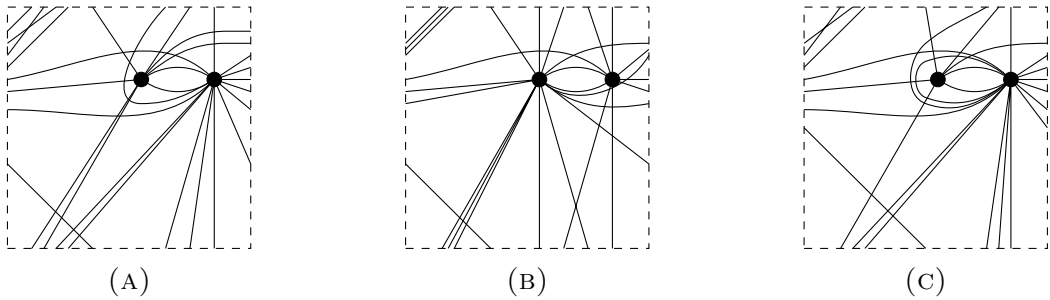

Three of these systems are described in \Cref{fig:J3 enum}. One is described in \Cref{fig:J21 final}. Three are described in \Cref{fig:J22 X1 syss} and two are described in \Cref{fig:J22 X2 syss}. Six systems are described in \Cref{fig:J23 third}. Two are described in \Cref{fig:J11 finals}. Three are described in \Cref{fig:J12 finals}. Finally, three are described in \Cref{fig:J0 finals}.

This investigation of maximal systems of arcs on surfaces originated with the study of systems of closed curves on surfaces. In~\cite{juvan_1996}, Juvan, Malni\v{c}, and Mohar introduce the term ``$k$-system'' to mean a collection of closed curves which intersect pairwise at most $k$ times. They then proved that for a fixed surface, the size of such a collection is uniformly bounded.

Przytycki~\cite{przytycki_2015} used $k$-systems of arcs to prove results about $k$-systems of curves. He proved that, for a closed surface \( S \) with marked points \( P \), the maximal cardinality of a 1-system of arcs on \( (S, P) \) is \( 2 | \chi | (|\chi| + 1) \), where \( \chi \) is the Euler characteristic of \( S \setminus P \). Using this fact, he was able to prove a non-sharp bound on the cardinality of 1-systems of curves, cubic in \( |\chi| \). 

In this paper, we also prove a generalization of a result proved in~\cite{przytycki_2015}. For a compact surface, possibly with boundary, and a finite set of marked points \( P \subset S \), let \( g \) be the genus of \( S \), let \( b \) be the number of boundary components of \( S \), and set \( p = |P \cap \text{int}(S)| \) and \( v = |P \cap \partial S| \). Define \( \chi\left( S, P \right) \vcentcolon = 2 - 2g - b - p - \frac{v}{2} \). We prove the following theorem, which may be of independent interest, to aid in our analysis. The proof of the theorem is presented in \Cref{sec:bounds}.

\begin{theorem}
	A maximal 1-system on \( (S, P) \) has cardinality
	\[ 2 |\chi| (|\chi| + 1) - \frac{v}{2} .\]
	\label{lem:main lem}
\end{theorem}

In~\cite{aougab_2019}, Aougab, Biringer, and Gaster used Przytycki's bound on the cardinality of 1-systems of arcs to prove that the cardinality of a 1-system of curves is bounded by \( O\left(\frac{|\chi| ^3}{(\log|\chi|)^2}\right) \). In ~\cite{greene_2019}, Greene improved this further to \( O \left( |\chi|^2 \log|\chi| \right) \).

We note that classifying maximal 0-systems is equivalent to classifying triangulations, and classifying 1-systems reveals structure of the surface in a similar way. An analysis of this kind was done by Tee in~\cite{tee_2021}, who considered the case where \( S \) is a sphere and \( P \) is a set of four marked points. He found 9 equivalence classes of maximal 1-systems. We remark that the maximal 1-systems on a torus with two punctures is a more diverse collection, with 23 members. 

Tee defined the \textit{non-intersecting subset} \( J \) for a 1-system \( \mathcal{A} \) to be the set of arcs disjoint from every other arc. Like in~\cite{tee_2021}, we use this property to break up our analysis. 


The rest of the paper is organized as follows. In \Cref{sec:prelim}, we develop some tools. We devote \Cref{sec:bounds} to proving our result about 1-systems on surfaces with boundary. In \Cref{sec:Jge2} we analyze the case where \( |J| \ge 2 \). In \Cref{sec:J1} we analyze the case where \( |J| = 1 \), and in \Cref{sec:J0} we analyze the case where \( J = \emptyset \).

\section{Preliminaries}
\label{sec:prelim}

Let \( S \) be a compact surface, possibly with boundary, and let \( P \subset S \) be a finite set of marked points. Let \( \mathcal{A} \) be a maximal 1-system on \( (S, P) \). 

\begin{rem}
	If a 1-system \( \mathcal{A} \) on \( (S, P) \) is maximal, then it is \emph{saturated}. That is, for any arc \( u \) on \( (S, P) \) not in \( \mathcal{A} \), \( \mathcal{A} \cup \left\{ u \right\} \) is not a 1-system. In particular, \( u \) is homotopic to some arc in \( \mathcal{A} \) or \( u \) intersects some arc of \( \mathcal{A} \) at least twice.
\end{rem}

Let \( u \) be an arc on \( (S, P) \).
If $u(0) = u(1)$, we say that $u$ is a \textit{loop} and that $u$ is \textit{based at} the marked point $u(0) = u(1)$. Otherwise, we say that $u$ is a \textit{non-loop arc} and that $u$ is \textit{between} the marked points $u(0)$ and $u(1)$.

\begin{defin}
	Define a (possibly disconnected) surface $S'$ with marked points \( P' \) to be \emph{almost embedded in $(S, P)$} if it admits a map (an \emph{almost embedding}) $p: S' \to S$ which is an embedding when restricted to $S' \setminus \partial S'$ such that \( p^{-1}(P) = P' \). 
\end{defin}

The following remark follows from the fact that almost embeddings are embeddings on their interior.
\begin{rem}
	Let \( p: S' \to S \) be an almost embedding. If $u$ is an arc on $(S, P)$, then the set $p^{-1}(u)$ is a disjoint collection of paths. Furthermore, if \( u(0,1) \subset p(S' \setminus \partial S') \), then the closure of \( p^{-1} (u(0,1)) \) is an arc on \( (S', P') \). 
\end{rem}

Let $ B$ be a non-empty 0-system on $(S, P)$. Define $S \setminus B$, or \textit{$S$ cut along $B$} to be the completion of $S \setminus \bigcup_{u \in B} u$. Define the \textit{gluing map} $p: S \setminus B \to S$ to be the identity away from the boundary, and the unique extension on the boundary. Note that the gluing map $p$ is an almost embedding.

If \( p: S \setminus B \to S \) is a gluing map, then \( p^{-1} (P) \) is a finite set. Then, each arc $u$ on \( (S \setminus B, p^{-1} (P)) \) naturally induces an arc on \( (S, P) \) which we call \( p(u) \). We may refer to this arc also as \( u \), abusing notation. 

\begin{rem}
	Every 1-system on \( (S, P) \) which is disjoint from \( B \) is induced by a 1-system on \( (S \setminus B, p^{-1}(P)) \). 
\end{rem}

We make the following additional note: each 1-system on \( S \setminus B \) may induce multiple non-equivalent 1-systems on \( S \), because \( S \setminus B \) may admit homeomorphisms which do not extend to \( S \).

Let \( J \) be the non-intersecting set of \( \mathcal{A} \). We delay the proof of \Cref{lem:help conn J} and \Cref{lem:help small J} to another section. 

\begin{lemma}
	\( S \setminus J \) is connected.
	\label{lem:help conn J}
\end{lemma}

\begin{cor}
	\label{lem:help small J}
	If \( S \) is a torus with \( |P| = 2 \), we have \( |J| \le 3 \).
\end{cor}

In this paper we may sometimes assume that a pair of arcs \( u, v \) are in minimal position: that is, they intersect the fewest number of times out of any pair \( u', v' \), where \( u', v' \) are homotopic as arcs to \( u, v \), respectively. In particular, if we equip the surface \( S \setminus P \) with a complete hyperbolic metric, then each arc may be realized as a geodesic, and in this situation all arcs will be pairwise in minimal position. See \cite{farb_2012}*{Prop 1.3 and Cor 1.9} and the comments on arcs on page 35 therein. 
So, for any 1-system there exists an equivalent 1-system whose arcs are pairwise in minimal position.

To determine if arcs are in minimal position, we use the bigon condition from~\cite{farb_2012}*{Prop 1.7}, taking into account the comments on arcs on page 35 of the same. This states that a pair of arcs are in minimal position if and only if they form no bigons, as shown in \Cref{fig:help bigons a}, or half-bigons, as shown in \Cref{fig:help bigons b}. 

\begin{figure}
	\centering
	\begin{subfigure}{0.2\textwidth}
		\centering
		\begin{tikzpicture}
			\draw (0,0) circle (1);
			\fill[color=white!70!red] (0,0.3) to[out=-135,in=135] (0,-0.3) to[out=45,in=-45] (0,0.3);

			\draw[green] (0.3,0.6) -- (0,0.3) to[out=-135,in=135] (0,-0.3) -- (0.3,-0.6);
			\draw (-0.3,0.6) -- (0,0.3) to[out=-45,in=45] (0,-0.3) -- (-0.3,-0.6);

			\draw[fill] (0.3,0.6) circle (0.04);
			\draw[fill] (-0.3,0.6) circle (0.04);
			\draw[fill] (-0.3,-0.6) circle (0.04);
			\draw[fill] (0.3,-0.6) circle (0.04);
		\end{tikzpicture}
		\caption{}
		\label{fig:help bigons a}
	\end{subfigure}
	\begin{subfigure}{0.2\textwidth}
		\centering
		\begin{tikzpicture}
			\draw (0,0) circle (1);
			\fill[color=white!70!red] (0,0.3) to[out=-135,in=135] (0,-0.3) to[out=45,in=-45] (0,0.3);

			\draw[green] (0.3,0.6) -- (0,0.3) to[out=-135,in=135] (0,-0.3);
			\draw (-0.3,0.6) -- (0,0.3) to[out=-45,in=45] (0,-0.3);
			\draw[fill] (0.3,0.6) circle (0.04);
			\draw[fill] (-0.3,0.6) circle (0.04);
			\draw[fill] (0,-0.3) circle (0.04);
		\end{tikzpicture}
		\caption{}
		\label{fig:help bigons b}
	\end{subfigure}
	\caption{}
	\label{fig:help bigons}
\end{figure}

Let $S'$ be a surface which is a closed disk with two marked points, at least one of which is on the boundary. That is, one of the two surfaces shown in \Cref{fig:help lem twos}. Let $u$ be the unique arc (up to homotopy) connecting the two marked points of $S'$ inside $S'$. This is shown in red in \Cref{fig:help lem twos}. We call \( u \) the \emph{internal arc} of \( (S', P') \). 
\begin{lemma} 
\label{lem:twos}
Let \( S', P', u \) be as above. Suppose $p: S' \to S$ is an almost embedding and let $w$ be an arc on $(S, P)$. If \( u,w \) are in minimal position, then \( p^{-1} (w) \) must intersect \( \partial S' \setminus P' \) in at least twice as many points as it intersects \( u \). 
\end{lemma}

\begin{figure}
	\centering
	\begin{subfigure}{0.3\textwidth}
		\centering
		\begin{tikzpicture}[scale=0.6]
			\draw (0,0) circle (1.5);

			\draw[red] (0,1.5) -- (0,-1.5);

			\draw[fill] (0,1.5) circle (0.06)
			(0,-1.5) circle (0.06);

		\end{tikzpicture}
		\caption{}
		\label{fig:help lem twos a}
	\end{subfigure}
	\begin{subfigure}{0.3\textwidth}
		\centering
		\begin{tikzpicture}[scale=0.6]
			\draw (0,0) circle (1.5);
			\draw[red] (0,1) -- (0,-1.5);
			\draw[fill] (0,1) circle (0.06)
			(0,-1.5) circle (0.06);
		\end{tikzpicture}
		\caption{}
		\label{fig:help lem twos b}
	\end{subfigure}
	\caption{}
	\label{fig:help lem twos}
\end{figure}

\begin{proof}
	Let \( x \) be a connected component of \( p^{-1} (w) \) which intersects \( u \). Then \( x \) is a path whose endpoints must either be at a marked point or on the boundary of \( S' \). 

	First suppose \( S' \) is the surface shown in \Cref{fig:help lem twos a}. Up to symmetry, homotopy, and moving endpoints along \( \partial S' \setminus P' \), the only paths on \( S' \) are those shown in \Cref{fig:help lem twos second}. The only one of these paths which intersects \( u \) is the one shown in \Cref{fig:help lem twos second a}. It must be that \( x \) is this path, so we see that \( x \) intersects \( u \) exactly once and has both endpoints on \( \partial S' \setminus P' \). 

\begin{figure}
	\centering
	\begin{subfigure}{0.2\textwidth}
		\centering
		\begin{tikzpicture}[scale=0.6]
			\draw (0,0) circle (1.5);

			\draw[blue] (-1.5,0) -- (1.5,0);
			\draw[red] (0,1.5) -- (0,-1.5);
			\draw[fill] (0,1.5) circle (0.06)
			(0,-1.5) circle (0.06);
		\end{tikzpicture}
		\caption{}
		\label{fig:help lem twos second a}
	\end{subfigure}
	\begin{subfigure}{0.2\textwidth}
		\centering
		\begin{tikzpicture}[scale=0.6]
			\draw (0,0) circle (1.5);
			\draw[blue] (0,-1.5) to[out=120,in=-120] (0,1.5);
			\draw[red] (0,1.5) -- (0,-1.5);
			\draw[fill] (0,1.5) circle (0.06)
			(0,-1.5) circle (0.06);
		\end{tikzpicture}
		\caption{}
		\label{fig:help lem twos second b}
	\end{subfigure}
	\begin{subfigure}{0.2\textwidth}
		\centering
		\begin{tikzpicture}[scale=0.6]
			\draw (0,0) circle (1.5);
			\draw[blue] (0,-1.5) to[out=120,in=-30] (-1.5,0);
			\draw[red] (0,1.5) -- (0,-1.5);
			\draw[fill] (0,1.5) circle (0.06)
			(0,-1.5) circle (0.06);
		\end{tikzpicture}
		\caption{}
		\label{fig:help lem twos second c}
	\end{subfigure}
	\begin{subfigure}{0.2\textwidth}
		\centering
		\begin{tikzpicture}[scale=0.6]
			\draw (0,0) circle (1.5);
			\draw[blue] (-60:1.5) to[out=120,in=-120] (60:1.5);
			\draw[red] (0,1.5) -- (0,-1.5);
			\draw[fill] (0,1.5) circle (0.06)
			(0,-1.5) circle (0.06);
		\end{tikzpicture}
		\caption{}
		\label{fig:help lem twos second d}
	\end{subfigure}
	\caption{}
	\label{fig:help lem twos second}
\end{figure}

Next suppose \( S' \) is the surface shown in \Cref{fig:help lem twos b}. The only paths on \( S' \) are those shown in \Cref{fig:help lem twos third}. The only one of these paths which intersects \( u \) is the one shown in \Cref{fig:help lem twos third a}. It must be that \( x \) is this path, so again \( x \) intersects \( u \) exactly once and has both endpoints on \( \partial S' \setminus P' \).

	\begin{figure}
		\centering
		\begin{subfigure}{0.18\textwidth}
			\centering
			\begin{tikzpicture}[scale=0.6]
				\draw (0,0) circle (1.5);
				\draw[blue] (-1.5,0) -- (1.5,0);
				\draw[red] (0,0.5) -- (0,-1.5);
				\draw[fill] (0,0.5) circle (0.06)
				(0,-1.5) circle (0.06);
			\end{tikzpicture}
			\caption{}
			\label{fig:help lem twos third a}
		\end{subfigure}
		\begin{subfigure}{0.18\textwidth}
			\centering
			\begin{tikzpicture}[scale=0.6]
				\draw (0,0) circle (1.5);
				\draw[blue] (0,0.5) -- (45:1.5);
				\draw[red] (0,0.5) -- (0,-1.5);
				\draw[fill] (0,0.5) circle (0.06)
				(0,-1.5) circle (0.06);
			\end{tikzpicture}
			\caption{}
			\label{fig:help lem twos third b}
		\end{subfigure}
		\begin{subfigure}{0.18\textwidth}
			\centering
			\begin{tikzpicture}[scale=0.6]
				\draw (0,0) circle (1.5);
				\draw[blue] (0,-1.5) to[out=60,in=-60] (0,0.5);
				\draw[red] (0,0.5) -- (0,-1.5);
				\draw[fill] (0,0.5) circle (0.06)
				(0,-1.5) circle (0.06);
			\end{tikzpicture}
			\caption{}
			\label{fig:help lem twos third c}
		\end{subfigure}

		\begin{subfigure}{0.18\textwidth}
			\centering
			\begin{tikzpicture}[scale=0.6]
				\draw (0,0) circle (1.5);
				\draw[blue] (0,-1.5) to[out=120,in=-30] (-1.5,0);
				\draw[red] (0,0.5) -- (0,-1.5);
				\draw[fill] (0,0.5) circle (0.06)
				(0,-1.5) circle (0.06);
			\end{tikzpicture}
			\caption{}
			\label{fig:help lem twos third d}
		\end{subfigure}
		\begin{subfigure}{0.18\textwidth}
			\centering
			\begin{tikzpicture}[scale=0.6]
				\draw (0,0) circle (1.5);
				\draw[blue] (0,-1.5) to[out=105,in=180] (0,1) to[out=0,in=75] (0,-1.5);
				\draw[red] (0,0.5) -- (0,-1.5);
				\draw[fill] (0,0.5) circle (0.06)
				(0,-1.5) circle (0.06);
			\end{tikzpicture}
			\caption{}
			\label{fig:help lem twos third e}
		\end{subfigure}
		\begin{subfigure}{0.18\textwidth}
			\centering
			\begin{tikzpicture}[scale=0.6]
				\draw (0,0) circle (1.5);
				\draw[blue] (-60:1.5) to[out=120,in=-120] (60:1.5);
				\draw[red] (0,0.5) -- (0,-1.5);
				\draw[fill] (0,0.5) circle (0.06)
				(0,-1.5) circle (0.06);
			\end{tikzpicture}
			\caption{}
			\label{fig:help lem twos third f}
		\end{subfigure}
		\caption{}
		\label{fig:help lem twos third}
	\end{figure}

	In both cases, each connected component of \( p^{-1} (w) \) corresponds to either exactly one intersection with \( x \) and two intersections with \( \partial S' \setminus P' \) or no intersections with \( x \) and possibly some intersections with \( \partial S' \setminus P' \). So, we are done.
\end{proof}

\begin{cor}
\label{cor:twoss}
Let \( S', P' \) be as above and let \( p: S' \to S \) be an almost embedding which is injective at all but finitely many points. If \( p(\partial S') \) is contained in the union of three arcs of \( \mathcal{A} \), then the internal arc of \( (S', P') \) must be included in \( \mathcal{A} \). 
\end{cor}

\begin{proof}
	Let \( u \) be the internal arc of \( (S', P') \). 
	Suppose for contradiction that \( u \notin \mathcal{A} \). Since \( \mathcal{A} \) is saturated, there must be some arc \( w \in \mathcal{A} \) which intersects \( u \) twice. 

	By \Cref{lem:twos}, the set \( E := p^{-1} (w) \cap (\partial S' \setminus P') \) must contain at least four points. First suppose \( p \) is injective on \( E \). If this is the case, then \( w \) intersects \( p(\partial S') \setminus P \) in at least four points. Since \( p(\partial S) \) in contained in three arcs, it must be that \( w \) intersects one of these arcs twice. This contradicts the assumption that \( \mathcal{A} \) is a 1-system.

	Now suppose \( p \) maps two of the points in \( E \) to the same point of \( S \). Since we assume that \( p \) is injective on all but finitely many points, this can only happen if that point on \( S \) is the intersection of two or more arcs, as shown in \Cref{fig:help cor twoss}. So, when this happens, the two points in \( E \) still correspond to intersections of \( w \) with two of the arcs which contain \( p(\partial S') \). Therefore, we find the same contradiction as above. 
	\begin{figure}
		\centering
		\begin{tikzpicture}
			\fill[color=white!70!red] (30:0.6) to[out=45,in=90] (1,0) to[out=-90,in=-45] (-30:0.6) -- (0,0) -- (-150:0.6) to[out=-45,in=-135] (1,-0.6) to[out=45,in=-45] (1,0.6) to[out=135,in=45] (150:0.6) -- (0,0);
			\draw[green,thick] (0,-0.4) -- (0,0.4);
			\draw (30:1) -- (30:-1)
			(-30:1) -- (-30:-1);
		\end{tikzpicture}
		\caption{The image of \( S' \) is shown shaded in red, a subarc of \( w \) is shown in green, and the black paths are subarcs of two of the three arcs which contain \( p(\partial S') \).}
		\label{fig:help cor twoss}
	\end{figure}
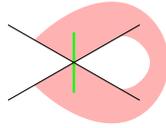
\end{proof}

\begin{cor}
	Suppose \( \mathcal{A} \) contains a loop arc which bounds a disk with one marked point on the interior. Then, the internal arc of that disk must be in \( \mathcal{A} \). Furthermore, that arc is disjoint from every other arc in \( \mathcal{A} \). 
	\label{cor:help ones}
\end{cor}

\begin{proof}
	Since the arc \( u \) bounds a disk with one marked point on the interior, we have satisfied the hypotheses of \Cref{cor:twoss} and we conclude that the internal arc of that disk, \( w \), is in \( \mathcal{A} \). 

	For the second part, suppose for contradiction that some other arc \( x \in \mathcal{A} \) intersects \( w \). Applying \Cref{lem:twos}, this implies that \( x \) intersects \( u \) at least twice. This contradicts the assumption that \( \mathcal{A} \) is a 1-system. 
\end{proof}

\section{Maximal 1-system Formula}
\label{sec:bounds}

In this section, we prove \Cref{lem:main lem}, as well as \Cref{lem:help conn J} and \Cref{lem:help small J}. In order to do so, we must first prove some intermediate results.

Let \( S \) be a compact connected surface, possibly with boundary. Let \( P\subset S \) be a finite set of marked points such that each component of \( \partial S \) contains at least one marked point. We equip \( S \setminus P \) with a complete hyperbolic metric such that the boundary is geodesic. As discussed in \Cref{sec:prelim}, we assume that all arcs on \( (S, P) \) are geodesics with respect to this metric. 

Let \( g \) be the genus of \( S \), let \( b \) be the number of boundary components of \( S \), and set \( p = |P \cap \text{int}(S)| \) and \( v = |P \cap \partial S| \). Define \( \chi\left( S, P \right) \vcentcolon = 2 - 2g - b - p - \frac{v}{2} \). We observe that when \( P \subset \text{int} (S) \), we have \( \chi(S, P) = \chi (S \setminus P) \), where \( \chi(S) \) is the usual Euler characteristic. If we write \( \chi \) this is understood to mean \( \chi(S, P) \). 

\begin{figure}
	\centering
	\begin{tikzpicture}
		\draw (0,0) circle (1);

		\draw[fill] (30:1) circle (0.05)
		(150:1) circle (0.05)
		(-90:1) circle (0.05);

		\node at (1.6,0) {$\to$};
	\end{tikzpicture}
	\begin{tikzpicture}
		\draw (-90:1.4) to[out=90,in=-150] (30:1.4) to[out=-150,in=-30] (150:1.4) to[out=-30,in=90] (-90:1.4);
	\end{tikzpicture}
	\caption{An example of \( S \setminus P \) when \( S \) is a disk and \( P \) consists of three points in the boundary of \( S \).}
	\label{fig:nons cusps}
\end{figure}
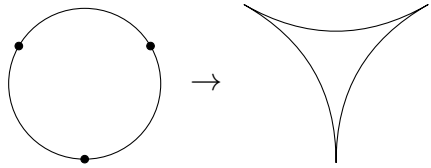

We treat punctures on the boundary of our surface as cusps, with 0 degree angles. For example, see \Cref{fig:nons cusps}. Using the Gauss-Bonnet theorem, we see that the hyperbolic area of the surface \( S \setminus P \) is equal to \( 2 \pi \chi(S, P) \). Therefore, if \( B \) is a 0-system on \( (S, P) \), then we have the natural almost embedding \( q: S \setminus B \to S \), and
\[ \chi(S, P) = \sum_{i=1}^{k} \chi(S_i, P_i), \]
where \( \left\{ S_i \right\}_{i=1}^k \) are the \( k \) connected components of \( S \setminus B \), and \( P_i = q^{-1} (P) \cap S_i \). In addition, we have \( \chi (S_i, P_i) < 0 \).

Fix a 1-system \( \mathcal{A} \) on \( (S, P) \). The following definition generalizes~\cite{przytycki_2015}*{Def 2.1}. 

\begin{defin}
	A \emph{tip} \( \tau \) of \( \mathcal{A} \) at a marked point \( x \in P \) is a pair \( (\alpha, \beta) \), where each of \( \alpha, \beta \) is either an oriented arc of \( \mathcal{A} \) or a component of \( \partial S \setminus P \), such that \( \alpha, \beta \) are \emph{consecutive}. 
\end{defin}

\begin{figure}
	\centering
	\begin{subfigure}[c]{0.3\textwidth}
		\centering
		\begin{tikzpicture}
			\fill[color=white!70!red] (0,-1) -- (0.18,-0.4) -- (-0.18,-0.4);
			\draw (0,-1) arc (-90:90:1);
			\draw[blue,thick] (0,-1) -- (0.3,0);

			\draw[blue,thick] (-0.3,0) -- (0,-1);
			\draw (0,-1) arc (270:90:1);

			\draw[green] (0.3,0) to[out=-160,in=-20] (-0.3,0);

			\draw[fill] (0.3,0) circle (0.04)
			(-0.3,0) circle (0.04)
			(0,-1) circle (0.04)
			(0,1) circle (0.04);
		\end{tikzpicture}
		\caption{A tip where both elements are arcs.}
	\end{subfigure}
	\begin{subfigure}[c]{0.3\textwidth}
		\centering
		\begin{tikzpicture}
			\fill[color=white!70!red] (0,-1) arc (-90:-120:1) -- (-0.18,-0.4);
			\draw (0,-1) arc (-90:90:1);
			\draw (0,-1) -- (0.3,0);

			\draw[blue,thick] (-0.3,0) -- (0,-1);
			\draw[blue,thick] (0,-1) arc (270:90:1);

			\draw[green] (-0.3,0) to[out=95,in=-110] (0,1);

			\draw[fill] (0.3,0) circle (0.04)
			(-0.3,0) circle (0.04)
			(0,-1) circle (0.04)
			(0,1) circle (0.04);
		\end{tikzpicture}
		\caption{A tip where one element is in the boundary.}
	\end{subfigure}
	\begin{subfigure}[c]{0.3\textwidth}
		\centering
		\begin{tikzpicture}
			\fill[color=white!70!red] (0.3,0) ++(-16.7:0.05) -- ++(-106.7:0.25) arc (-106.7:73.3:0.25) -- ++(163.3:0.1) arc(73.3:253.3:0.25) -- ++(73.3:0.25);
			\draw (0,-1) arc (-90:90:1);
			\draw[blue,thick] (0,-1) -- (0.3,0);

			\draw (-0.3,0) -- (0,-1);
			\draw (0,-1) arc (270:90:1);

			\draw[green] (0,-1) to[out=90,in=-107] (0,0.24) to[out=73,in=73] (0.7,0.05) to[out=-107,in=56] (0,-1);
			\draw[fill] (0.3,0) circle (0.04)
			(-0.3,0) circle (0.04)
			(0,-1) circle (0.04)
			(0,1) circle (0.04);
		\end{tikzpicture}
		\caption{A tip where \( \alpha = \beta \).} 
	\end{subfigure}
	\caption{Some tips. The elements of the tip are colored blue, and a disc witnessing the tip is shown in red. The image of the third side of the associated nib is shown in green.}
	\label{fig:nons tips}
\end{figure}
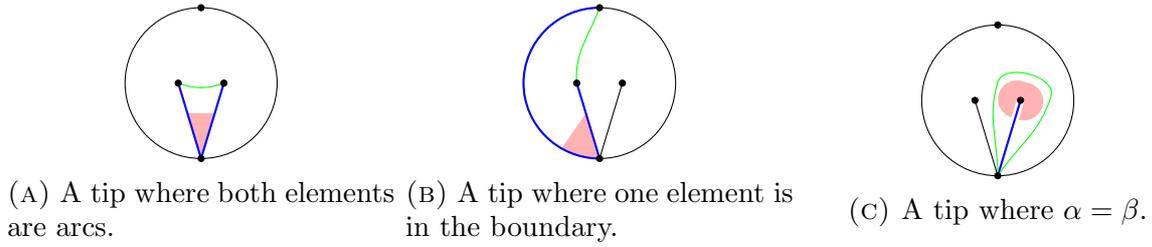

By consecutive we mean that \( S \) admits an orientation preserving embedding of the disk \( [0,1] \times [0,1] \) such that \( \left\{ 0 \right\} \times [0,1] \) maps into \( \alpha \), \( (0,0) \) maps to \( x \), \( [0,1] \times \left\{ 0 \right\} \) maps into \( \beta \), and \( (0,1) \times (0,1) \) intersects neither \( P \) nor any arc of \( \mathcal{A} \). See \Cref{fig:nons tips}.

Let \( \tau = (\alpha, \beta) \) be a tip at \( x \in P \) and let \( N_\tau \) be an abstract open ideal hyperbolic triangle with vertices \( a, t, b \). We may associate to \( \tau \) a unique local isometry \( \nu_\tau : N_\tau \to S \setminus P \) which sends \( t \) to \( x \), \( ta \) to \( \alpha \), and \( tb \) to \( \beta \). Call \( \nu_\tau \) the \emph{nib} of \( \tau \). This is the same as the concept of nibs introduced in~\cite{przytycki_2015}.

The following proposition is the same as~\cite{przytycki_2015}*{Prop 2.2}. 
\begin{prop}
	Let \( \mathcal{A} \) be a 1-system on \( (S, P) \) and let \( \nu : N = \bigsqcup N_\tau \to S \) be the disjoint union of all the nibs. Then for each \( s \in \text{int}(S) \) the preimage \( \nu^{-1} (s) \) consists of at most \( 2 (|\chi| + 1) \) points.
	\label{prop:2.2}
\end{prop}

The proof of this proposition relies on the following three statements. They are similar to statements found in~\cite{przytycki_2015}, adapted for our setting.

\begin{lemma}
	Let \( x \in P \) be a marked point in the interior of \( S \). If a 0-system \( B \) contains only arcs joining \( x \) to a marked point in \( P \setminus \left\{ x \right\} \) then \( |B| \le 2|\chi| \).
	\label{lem:2.3}
\end{lemma}
\begin{proof}
	Without loss of generality, suppose \( B \) is maximal. Then the components of \( S \setminus B \) are either squares, or triangles where one side is in the boundary of \( S \). The area of \( S \setminus P \) is \( 2\pi|\chi| \). The number of triangles is equal to \( v \), the number of marked points in the boundary, and each has area \( \pi \). Therefore the squares have total area \( \pi (2 |\chi| - v) \), and so there are \( |\chi| - \frac{v}{2} \) squares. Each square has four arcs in its boundary, each triangle has two arcs in its boundary (the other side of the triangle was in \( \partial S \)), and each arc contributes to this total twice. So, the number of arcs is 
	\[ \left( 4 \left( |\chi| - \frac{v}{2} \right) + 2v \right) / 2 = 2|\chi|.\qedhere \]
\end{proof}

Let \( n \in N_\tau \) for some tip \( \tau \). We define the \emph{slit} at \( n \) to be the restriction of \( \nu_\tau \) to the geodesic ray starting at \( n \) and going towards the vertex \( t \). This is the same definition as the one given in~\cite{przytycki_2015}.

\begin{lemma}
	A slit is an embedding.
	\label{lem:2.5}
\end{lemma}

\begin{proof}
	Let \( n \in N_\tau \) for some tip \( \tau \). We define the doubled surface \( \tilde{S} \) to be the disjoint union of two copies of \( S \), glued along the boundary, and we define \( \tilde{P} \) similarly to be twice \( P \), identified when the marked point is on \( \partial S \). 

	The inclusion \( S \hookrightarrow \tilde{S} \) is an isometric embedding. In particular, the map takes geodesics to geodesics.

	If \( \alpha \) is an arc on \( (S, P) \), define \( \tilde{\alpha} \) to be the inclusion of that arc into \( \tilde{S} \). If instead we have \( \alpha \subset \partial S \), then define \( \tilde{\alpha} \) to be the arc induced by the inclusion of \( \alpha \) into \( \big( \tilde{S}, \tilde{P} \big) \). 

	We define \( \tilde{\beta} \) similarly. Now, \( \tilde{\tau} := \big( \tilde{\alpha}, \tilde{\beta} \big) \) is a tip, and its nib is exactly the inclusion of \( \nu_\tau \) into \( \tilde{S} \). So, by~\cite{przytycki_2015}*{Lem 2.5}, the slit at \( n \) is an embedding into \( \tilde{S} \), and hence into \( S \). 
\end{proof}

\begin{lemma}
	Suppose the arcs in \( \mathcal{A} \) intersect pairwise at most once. If for distinct \( n, n' \in N \) we have \( \nu(n) = \nu(n') \), then the slits at \( n, n' \) are disjoint except at the endpoints. 
	\label{lem:2.6}
\end{lemma}

\begin{proof}
	This proof proceeds similarly to the previous one. We consider the doubled marked surface \( \big( \tilde{S}, \tilde{P} \big) \). If \( n_1, n_2 \) are in the nibs of \( \tau_1, \tau_2 \), respectively (and it might be that \( \tau_1 = \tau_2 \)), we may associate to them tips \( \tilde{\tau_1}, \tilde{\tau_2} \) whose elements are arcs on \( \big( \tilde{S}, \tilde{P} \big) \). By~\cite{przytycki_2015}*{Lem 2.6}, the inclusion of the slits at \( n, n' \) into \( \tilde{S} \), and hence into \( S \), are disjoint except at the endpoint. 
\end{proof}

\begin{proof}[Proof of \Cref{prop:2.2}]
	Let \( P' = P \cup \left\{ s \right\} \). Then \( |\chi (S, P')| = |\chi(S, P)| + 1 \). Every slit at \( s \) in \( (S, P) \) gives an arc on \( (S, P') \), so let \( \mathcal{S} \) be the collection of these arcs. The arcs in \( \mathcal{S} \) are simple by \Cref{lem:2.5} and pairwise disjoint by \Cref{lem:2.6}. Since \( s \) is in the interior of \( S \) and \( P' = P \cup \left\{ s \right\} \), by \Cref{lem:2.3} we have \( |\mathcal{S}| \le 2 |\chi (S, P')| = 2(|\chi(S, P)| +1)  \). 
\end{proof}

\begin{proof}[Proof of \Cref{lem:main lem}]
	First, we use a construction similar to the one in~\cite{przytycki_2015} to show that this cardinality is realized. We choose a 0-system \( B \) such that \( S \setminus B \) is a disk with no marked points on the interior. Let \( q: S \setminus B \to S \) be the natural almost embedding. We know \( \chi(S \setminus B, q^{-1}(P)) = \chi(S, P) \), and additionally \( S \setminus B \) has genus 0, 1 boundary component, and no marked points on the interior. Therefore, \( S \setminus B \) has \( 2 + 2|\chi| \) marked points on the boundary. Since cutting along an arc always adds two to the number of marked points on the boundary, we have \( 2|B| + v = 2|\chi| + 2 \), and \( S \setminus B \) has \( \frac{(2|\chi| + 2)(2|\chi| - 1)}{2} \) diagonals. If we let \( \mathcal{A} \) be the union of \( B \) and the diagonals of \( S \setminus B \), then:
	\begin{align*}
		|\mathcal{A}| &= |B| + \frac{(2|\chi| + 2)(2|\chi| - 1)}{2} \\
		&= |\chi| + 1 - \frac{v}{2} + (|\chi| + 1)(2|\chi| - 1) \\
		&= 2|\chi|(|\chi| + 1) - \frac{v}{2} \\
	\end{align*}
	So, it remains to give an upper bound. 

	To bound the number of arcs in \( \mathcal{A} \), we count the tips. Each arc in \( \mathcal{A} \) is the first element of exactly two tips, once in each orientation, and each component of \( \partial S \setminus P \) is the first element of exactly one tip. The number of components of \( \partial S \setminus P \) is exactly \( v \), the number of marked points in the boundary. So, the area of \( N \) is \( (2 |\mathcal{A}| + v) \pi \). The area of \( S \setminus P \) is \( 2 |\chi| \pi\). By \Cref{prop:2.2}, the map \( \nu \) is a most \( 2(|\chi| + 1) \) to 1, so we have:
	\begin{align*}
		(2 |\mathcal{A}| + v) \pi &\le 2 \pi |\chi| \cdot 2 (|\chi| + 1)\\
		2 |\mathcal{A}| + v &\le 4 |\chi| (|\chi| + 1) \\
		2 |\mathcal{A}| &\le 4 |\chi| (|\chi| + 1) - v \\
		|\mathcal{A}| &\le 2 |\chi| (|\chi| + 1) - \frac{v}{2} \qedhere
	\end{align*}
\end{proof}

\begin{proof}[Proof of \Cref{lem:help conn J}]
	For contradiction, suppose \( S \setminus J \) has two or more connected components. Let \( q : S \setminus J \to S \) be the natural almost embedding and let \( P' = q ^{-1} (P) \). Up to taking a subset of \( J \), we assume that \( S \setminus J \) has exactly two components, \( S_1, S_2 \). Let \( v = |P \cap \partial S| \). For \( i = 1,2 \), let \( P_i = S_i \cap P' \) and let \( v_i = |P_i \cap \partial S_i| \). 

	Using \Cref{lem:main lem}, we know that the maximum cardinality of a 1-system on \( (S_i,P_i) \) is \( 2 |\chi_i| (|\chi_i| + 1) - \frac{v_i}{2} \). If we choose a maximal 1-system on \( S_1, S_2 \), we obtain a 1-system on \( (S, P) \) by projecting all the arcs in both 1-systems. In fact, given \( J \), every 1-system on \( (S, P) \) arises this way, since every arc in \( \mathcal{A} \setminus J \) lies entirely within a single connected component of \( S \setminus J \). 

	We may also note that \( v_1 + v_2 = |P'| = 2|J| + v  \), since each arc we cut along adds 2 to the total number of marked points on the boundary of some connected component.

	Since \( \mathcal{A} \) is a maximal 1-system on \( (S, P) \), it must have cardinality \( 2 | \chi | (|\chi| +1) - \frac{v}{2} \). Conversely, if \( \mathcal{A} \) is obtained by projecting 1-systems from \( (S_1, P_1) \) and \( (S_2, P_2) \), then:
	\begin{align*}
		|\mathcal{A}| & \le 2 |\chi_1| (|\chi_1| + 1) - \frac{v_1}{2} + 2 |\chi_2| (|\chi_2| + 1) - \frac{v_2}{2} + |J| \\
		& = 2 \left(|\chi_1| (|\chi_1| + 1) + |\chi_2| (|\chi_2| + 1) \right) - \frac{v_1+v_2}{2} + |J| \\
		& = 2 \left(|\chi_1|^2 + |\chi_1| + |\chi_2|^2 +|\chi_2| \right) - \frac{v}{2} \\
		& < 2 \left(|\chi_1|^2 + |\chi_1| + |\chi_2|^2 +|\chi_2| + 2 |\chi_1||\chi_2| \right) - \frac{v}{2} \\
		&= 2 (|\chi_1| + |\chi_2|) \left( |\chi_1| + |\chi_2| + 1 \right) - \frac{v}{2} \\
		&= 2 | \chi | (|\chi| + 1) - \frac{v}{2}
	\end{align*}
	We have a strict inequality because, as noted above, \( \chi(S_i, P_i) < 0 \) and in particular \( \chi(S_i, P_i) \ne 0 \). 
	This contradicts the maximality of \( \mathcal{A} \). So, we are done.
\end{proof}

\begin{proof}[Proof of \Cref{lem:help small J}]
	Suppose for contradiction that \( |J| > 3 \). Let \( g, b \) be the genus and number of boundary components of \( S \setminus J \), respectively. We have \( g \ge 0 \) and \( b \ge 1 \), since genus is a non-negative number, and cutting creates at least one boundary component.

	Let \( q : S \setminus J \to S \) be the natural almost embedding. Define \( P' = q^{-1} (P) \). Let \( p = |P' \cap \text{int} S| \) and let \( v = |P' \cap \partial S| \). We have \( p \ge 0 \), since it is the size of a set, and as we remark above \( 2|J| = v \).

	\begin{align*}
		-2 = \chi(S, P) &= \chi(S \setminus J, P') \\
		&= 2 -2g - b - p - \frac{v}{2}\\
		&= 2 - 2g - b - p - |J|\\
		& < -1 -2g - b - p\\
		& \le -2
	\end{align*}
	This is a contradiction. So, we are done. 
\end{proof}

\section{$|J| \ge 2$}
\label{sec:Jge2}

Let $S$ be a torus and let $|P| = 2$. By \Cref{lem:help small J}, we have \( |J| \le 3 \). So, we analyze the case where \( |J| = 3 \), then we consider the three possible cases when \( |J| = 2 \).

\subsection{$|J| = 3$}

\begin{figure}
	\centering
	\begin{subfigure}{0.6\textwidth}
		\centering
		\begin{tikzpicture}
			\draw[dashed] (0,0) -- (2,0) -- (2,2) -- (0,2) -- (0,0);

			\draw[fill] (1.3,0.7) circle (0.06)
			(0.7,1.3) circle (0.06);

			\draw (0.7,0) -- (0.7,2)
			(0,1.3) -- (2,1.3)
			(0.7,1.3) -- (1.3,0.7);

		\end{tikzpicture}
		\quad
		\begin{tikzpicture}
			\draw[dashed] (0,0) -- (2,0) -- (2,2) -- (0,2) -- (0,0);

			\draw[fill] (1.3,0.7) circle (0.06)
			(0.7,1.3) circle (0.06);

			\draw (0.7,0) -- (0.7,2)
			(0.7,1.3) -- (1.3,0.7)
			(1.3,0.7) -- (2,1) (0,1) -- (0.7,1.3);

		\end{tikzpicture}
		\quad
		\begin{tikzpicture}
			\draw[dashed] (0,0) -- (2,0) -- (2,2) -- (0,2) -- (0,0);

			\draw[fill] (1.3,0.7) circle (0.06)
			(0.7,1.3) circle (0.06);

			\draw (0.7,1.3) -- (1,2) (1,0) -- (1.3,0.7)
			(0.7,1.3) -- (1.3,0.7)
			(1.3,0.7) -- (2,1) (0,1) -- (0.7,1.3);

		\end{tikzpicture}
		\caption{non-separating zero systems}
		\label{fig:J3 enum a}
	\end{subfigure}
	\begin{subfigure}{0.3\textwidth}
		\centering
		\begin{tikzpicture}
			\foreach \x in {0,60,...,300} {
				\draw[blue] (\x:1) -- ({\x+180}:1) -- ({\x+60}:1);
			}
			\foreach \x in {0,60,...,300}{
				\draw[fill] (\x:1) circle (0.06);
				\draw (\x:1) -- ({\x + 60}:1);
			}
		\end{tikzpicture}
		\caption{}
		\label{fig:J3 enum b}
	\end{subfigure}
	\caption{}
	\label{fig:J3 enum}
\end{figure}

By \Cref{lem:help conn J}, $J$ is non-separating. So we consider all possible non-separating collections of three disjoint arcs. These are shown in Figure~\ref{fig:J3 enum a}.

In all three cases, the surface $S \setminus J$ is a disc with six marked points on the boundary: a hexagon. Then, the only arcs which could be in $\mathcal{A}$ are the diagonals of the hexagon. There are nine diagonals, none of which intersect two or more times, since the arcs can be realized as straight lines. This is shown in \Cref{fig:J3 enum b}. So, for all three possible $J$ of this size, $\mathcal{A}$ must be composed of the three arcs in $J$ and the nine diagonals of the hexagon $S \setminus J$. 


\subsection{$J$ is two loops}

Since \( J \) is non-separating, it must consist of two loops based at the same vertex. 
Then, the surface $S - J$ is a disc with four marked points on the boundary and one marked point in the interior, as shown in \Cref{fig:J21 intro b}. Let \( h \) denote rotation of the disk by 90 degrees counterclockwise.

\begin{figure}
	\centering
	\begin{subfigure}[b]{0.3\textwidth}
		\centering
		\begin{tikzpicture}
			\draw[red] (0.7,0) -- (0.7,2)
			(0,1.4) -- (2,1.4);

			\draw[fill] (1.4,0.7) circle (.05)
			(0.7,1.4) circle (.05);

			\draw[dashed] (0,0) -- (2,0) -- (2,2) -- (0,2) -- (0,0);
		\end{tikzpicture}
		\caption{}
		\label{fig:J21 intro a}
	\end{subfigure}
	\begin{subfigure}[b]{0.3\textwidth}
		\centering
		\begin{tikzpicture}[scale=0.5]
			\draw[fill] (0,0) circle (.1)
			(4,0) circle (.1)
			(0,4) circle (.1)
			(4,4) circle (.1)
			(2,2) circle (.1);

			\draw (0,0) -- (4,0) -- (4,4) -- (0,4) -- (0,0);
		\end{tikzpicture}
		\caption{}
		\label{fig:J21 intro b}
	\end{subfigure}
	\caption{}
	\label{fig:J21 intro}
\end{figure}

We see that an arc on \( S - J \) with both endpoints at the marked point in the interior cannot be essential. So, we write \( A \setminus J = X \sqcup Y \), where \( X \) contains the arcs which have exactly one endpoint on the boundary of \( S - J \), and \( Y \) contains the arcs which have both endpoints on the boundary of \( S - J \).


\begin{figure}
	\centering
	\begin{subfigure}[b]{0.3\textwidth}
		\centering
		\begin{tikzpicture}[scale=0.5]
			\draw[red] (0,0) -- (2,2) -- (4,0);

			\draw[fill] (0,0) circle (.1)
			(4,0) circle (.1)
			(0,4) circle (.1)
			(4,4) circle (.1)
			(2,2) circle (.1);

			\draw (0,0) -- (4,0) -- (4,4) -- (0,4) -- (0,0);

			\node[red] at (1,1.6) {a};
			\node[red] at (3.3,1.6) {h(a)};

		\end{tikzpicture}
		\caption{}
		\label{fig:J21 class a}
	\end{subfigure}
	\begin{subfigure}[b]{0.6\textwidth}
		\centering
		\begin{tikzpicture}[scale=0.5]
			\draw[fill] (0,0) circle (.1)
			(4,0) circle (.1)
			(0,4) circle (.1)
			(4,4) circle (.1)
			(2,2) circle (.1);

			\draw (0,0) -- (4,0) -- (4,4) -- (0,4) -- (0,0);

			\draw[red]
			(0,0) to[out=75,in=-165] (4,4);

			\node[red] at (1,3) {b};
		\end{tikzpicture}
		\begin{tikzpicture}[scale=0.5]
			\draw[red] (0,0) to[out=65,in=-130] (1.3, 2.1) 
			to[out=50,in=180] (2,2.4) to[out=0,in=130]
			(2.7,2.1) to[out=-50,in=115] (4,0);
			\draw[fill] (0,0) circle (.1)
			(4,0) circle (.1)
			(0,4) circle (.1)
			(4,4) circle (.1)
			(2,2) circle (.1);

			\draw (0,0) -- (4,0) -- (4,4) -- (0,4) -- (0,0);

			\node[red] at (2,3) {c};

		\end{tikzpicture}
		\begin{tikzpicture}[scale=0.5]
			\draw[fill] (0,0) circle (.1)
			(4,0) circle (.1)
			(0,4) circle (.1)
			(4,4) circle (.1)
			(2,2) circle (.1);

			\draw (0,0) -- (4,0) -- (4,4) -- (0,4) -- (0,0);

			\draw[red] (0,0) to[out=55,in=-135] (1.8,2.2) to[out=45,in=135] (2.2,2.2) to[out=-45,in=45] (2.2,1.8) to[out=-135,in=35] (0,0);

			\node[red] at (2,1) {d};
		\end{tikzpicture}
		\caption{}
		\label{fig:J21 class b}
	\end{subfigure}
	\caption{}
	\label{fig:J21 class}
\end{figure}

We see that any arc in \( X \) must be of the form \( h^i (a) \) where \( a \) is the arc shown in \Cref{fig:J21 class a} and \( i \in \left\{ 0,1,2,3 \right\} \). In particular, this means \( |X| \le 4 \). 

Similarly, any arc in \( Y \) must be one of \( h^i(b), h^i(c), h^i(d) \) where \( b, c, d \) are the arcs shown in \Cref{fig:J21 class b} and \( i \in \left\{ 0,1,2,3 \right\} \).

We note that \( d \) bounds a disk with one marked point on the interior. So $\mathcal{A} \setminus J$ cannot include \( h^i(d) \) for any \( i \), since this would contradict the assumption that \( |J| = 2 \), by \Cref{cor:help ones}. In addition, we see that the pair \( c, h^2(c) \) intersects twice, and so does the pair \( h(c), h^3(c) \). This is shown in \Cref{fig:J21 crosses}. So, we see that \( |Y| \le 6 \). We conclude that in fact \( |X| = 4, |Y| = 6 \). 


\begin{figure}
	\centering
	\begin{tikzpicture}[scale=0.5]
		\draw[dashed] (0,0) to[out=65,in=-130] (1.3, 2.1) 
		to[out=50,in=180] (2,2.4) to[out=0,in=130]
		(2.7,2.1) to[out=-50,in=115] (4,0);
		\draw[dashed] (0,4) to[out=-65,in=130] (1.3, 1.9) 
		to[out=-50,in=180] (2,1.6) to[out=0,in=-130]
		(2.7,1.9) to[out=50,in=-115] (4,4);
		\draw[fill] (0,0) circle (.1)
		(4,0) circle (.1)
		(0,4) circle (.1)
		(4,4) circle (.1)
		(2,2) circle (.1);

		\draw (0,0) -- (4,0) -- (4,4) -- (0,4) -- (0,0);

	\end{tikzpicture}
	\quad
	\begin{tikzpicture}[scale=0.5]
		\draw[dashed] (0,0) to[out=25,in=-140] (2.1,1.3) 
		to[out=40,in=-90] (2.4,2) to[out=90,in=-40]
		(2.1,2.7) to[out=140,in=-25] (0,4);
		\draw[dashed] (4,0) to[out=165,in=-40] (1.9,1.3) 
		to[out=140,in=-90] (1.6,2) to[out=90,in=-140]
		(1.9,2.7) to[out=40,in=-165] (4,4);
		\draw[fill] (0,0) circle (.1)
		(4,0) circle (.1)
		(0,4) circle (.1)
		(4,4) circle (.1)
		(2,2) circle (.1);

		\draw (0,0) -- (4,0) -- (4,4) -- (0,4) -- (0,0);
	\end{tikzpicture}
	\caption{}
	\label{fig:J21 crosses}
\end{figure}

Without loss of generality, we have 
\[ \mathcal{A} \setminus J = \left\{ a, h(a), h^2(a), h^3(a), b, h(b), h^2(b), h^3(b), c, h(c) \right\}. \]
So there is only one maximal 1-system, up to equivalence. This is shown in \Cref{fig:J21 final}.

\begin{figure}
	\centering
	\begin{tikzpicture}[scale=0.8]
		\draw[dashed] (0,0) to[out=65,in=-130] (1.3, 2.1) 
		to[out=50,in=180] (2,2.4) to[out=0,in=130]
		(2.7,2.1) to[out=-50,in=115] (4,0);
		\draw[dotted] (0,0) to[out=25,in=-140] (2.1,1.3) 
		to[out=40,in=-90] (2.4,2) to[out=90,in=-40]
		(2.1,2.7) to[out=140,in=-23] (0,4);

		\draw[fill] (0,0) circle (.1)
		(4,0) circle (.1)
		(0,4) circle (.1)
		(4,4) circle (.1)
		(2,2) circle (.1);

		\draw (0,0) -- (4,0) -- (4,4) -- (0,4) -- (0,0)
		-- (4,4)
		(4,0) -- (0,4);

		\draw
		(0,4) to[out=-75,in=165] (4,0);
		\draw
		(4,0) arc (10:80:4.9)
		(4,4) arc (100:170:4.9)
		(0,0) arc (-80:-10:4.9);
	\end{tikzpicture}
	\caption{}
	\label{fig:J21 final}
\end{figure}

\subsection{$J$ contains a loop and a non-loop arc}


Here, the surface $S \setminus J$ is an annulus with three marked points in one boundary component and one marked point in the other boundary component, as shown in \Cref{fig:J22 intro}.

\begin{figure}
	\centering
	\begin{tikzpicture}[baseline=(current bounding box.center),scale=1.2]
		\draw[dashed] (0,0) -- (2,0) -- (2,2) -- (0,2) -- (0,0);

		\draw[fill] (0.6,1.4) circle (0.06)
		(1.4,0.6) circle (0.06);

		\draw (1.4, 0) -- (1.4,2)
		(1.4,0.6) -- (0.6,1.4);

		\node at (0.3,1.4) {a};
		\node at (1.7,0.6) {b};

		\node at (3,1) {$\to$};
	\end{tikzpicture}
	\quad
	\begin{tikzpicture}[baseline=(current bounding box.center),scale=.35]
		\draw (0,0) circle (4)
		(0,0) circle (1);
		\draw[fill] (0,1) circle (0.1)
		(0,4) circle (0.1)
		(4,0) circle (0.1)
		(-4,0) circle (0.1);

		\node at (0,1.5) {b};
		\node at (0,4.5) {a};
		\node at (4.5,0) {b};
		\node at (-3.5,0) {b};
	\end{tikzpicture}
	\caption{}
	\label{fig:J22 intro}
\end{figure}

Write \( \mathcal{A} \setminus J = X \sqcup Y \), where \( X \) is the set of arcs which have both endpoints on the same boundary component of \( S \setminus J \), and \( Y \) is the set of arcs which have one endpoint on each boundary component.

If an arc has both endpoints on the boundary component with one marked point, then it must be either homotopic to a constant map, or homotopic to the boundary. In either case, such an arc is not essential. So, we may assume that every arc \( x \in X \) has both endpoints on the boundary component with three marked points.

\begin{figure}
	\centering
	\begin{subfigure}[b]{0.6\textwidth}
		\centering
		\begin{tikzpicture}[scale=.4]
			\draw (0,0) circle (4)
			(0,0) circle (1);
			\draw[fill] (0,1) circle (0.1)
			(0,4) circle (0.1)
			(4,0) circle (0.1)
			(-4,0) circle (0.1);

			\draw (0,4) to[out=-50,in=90] (2,0) to[out=-90,in=0] (0,-2) to[out=180,in=-90] (-2,0) to[out=90,in=-130] (0,4);

			\node at (2.8,0) {$x_1$};
		\end{tikzpicture}
		\begin{tikzpicture}[scale=.4]
			\draw (0,0) circle (4)
			(0,0) circle (1);
			\draw[fill] (0,1) circle (0.1)
			(0,4) circle (0.1)
			(4,0) circle (0.1)
			(-4,0) circle (0.1);

			\draw (-4,0) to[out=-60,in=180] (0,-3) to[out=0,in=-90] (3,0) to[out=90,in=-25] (0,4);

			\node at (2.1,0) {$x_2$};
		\end{tikzpicture}
		\caption{}
		\label{fig:J22 one Xs}
	\end{subfigure}
	\begin{subfigure}[b]{0.3\textwidth}
		\centering
		\begin{tikzpicture}[scale=.4]
			\draw (0,0) circle (4)
			(0,0) circle (1);
			\draw[fill] (0,1) circle (0.1)
			(0,4) circle (0.1)
			(4,0) circle (0.1)
			(-4,0) circle (0.1);

			\draw (0,1) -- (0,4);

			\node at (1,2.7) {\( y \)};
		\end{tikzpicture}
		\caption{}
		\label{fig:J22 one Ys}
	\end{subfigure}
	\caption{}
	\label{fig:J22 one}
\end{figure}

Let \( h : S \setminus J \to S \setminus J \) be the homeomorphism that rotates the three marked points on the outer boundary component counterclockwise. Let \( x \in X \). Up to applying \( h \), \( x \) must be one of \( x_1, x_2 \) as shown in \Cref{fig:J22 one Xs}.

Note that \( x_1 \) intersects twice each \( h(x_1), h^2(x_1) \). In addition, \( x_1 \) intersects twice \( h(x_2) \). So, if \( x_1 \in X \) then we have \( X \subset X_1 = \left\{ x_1, x_2, h(x_2) \right\} \). If instead none of \( x_1, h(x_1), h^2 (x_1) \) are in \( X \), then \( X \subset X_2 = \left\{ x_2, h(x_2), h^2(x_2) \right\} \). In either case, we have \( |X| \le 3 \). The two maximal possibilities up to the action of \( h \) are shown in \Cref{fig:J22 max xs}.

Now note the following. Up to application of \( h \), there is a unique arc \( y \) that has one endpoint on each boundary component, as shown in \Cref{fig:J22 one Ys}. In fact, every member of \( Y \) is of the form \( y_k := h^{3k}(y) \) for some \( k \in \frac{1}{3} \mathbb{Z} \). We have a formula for the number of intersections of \( y_{k},  y_{\ell} \):
\begin{equation}
	i(y_k, y_\ell) = \left\lceil |k - \ell| \right\rceil - 1
	\label{eq:hk1 hk2 int}
\end{equation}

Using \cref{eq:hk1 hk2 int}, we find that \( |Y| \le 7 \). From the argument above, we conclude that in fact \( |Y| = 7 \) and \( |X| = 3 \).

\begin{figure}
	\centering
    \begin{tikzpicture}[scale=.4]
        \draw (0,0) circle (4)
        (0,0) circle (1);
        \draw[fill] (0,1) circle (0.1)
        (0,4) circle (0.1)
        (4,0) circle (0.1)
        (-4,0) circle (0.1);
        
        \draw[red] (4,0) to[out=-120,in=0] (0,-3) to[out=180,in=-90] (-3,0) to[out=90,in=-155] (0,4);
        \draw[blue] (-4,0) to[out=-60,in=180] (0,-3) to[out=0,in=-90] (3,0) to[out=90,in=-25] (0,4);
        \draw (0,4) to[out=-50,in=90] (2,0) to[out=-90,in=0] (0,-2) to[out=180,in=-90] (-2,0) to[out=90,in=-130] (0,4);

	\node at (0,-5) {\( X_1 \)};
    \end{tikzpicture}
    \qquad
    \begin{tikzpicture}[scale=.4]
        \draw (0,0) circle (4)
        (0,0) circle (1);
        \draw[fill] (0,1) circle (0.1)
        (0,4) circle (0.1)
        (4,0) circle (0.1)
        (-4,0) circle (0.1);
        
        \draw[red] (4,0) to[out=-120,in=0] (0,-3) to[out=180,in=-90] (-3,0) to[out=90,in=-155] (0,4);
        \draw[blue] (-4,0) to[out=-60,in=180] (0,-3) to[out=0,in=-90] (3,0) to[out=90,in=-25] (0,4);
        \draw (4,0) to[out=120,in=0] (0,3) to[out=180,in=60] (-4,0);

	\node at (0,-5) {\( X_2 \)};
    \end{tikzpicture}
    \caption{}
    \label{fig:J22 max xs}
\end{figure}

Since \( |X| = 3 \), it must be that \( X= X_2 \) or \( X = h^i (X_1) \) for \( i = 0, 1, 2 \). Note that \( h^3(X_1) = X_1 \). In fact, \( h(X_1) \) is related to \( h^2(X_1) \) by a homeomorphism that projects to \( S \), so we may assume that \( i = 0,1 \). 

Let \( Y_* = \left\{ h^{-3}(y), h^{-2}(y), \dots, h^3(y) \right\} \). Note that every arc in \( Y_* \) intersects at least one other arc in \( Y_* \), except \( y = h^0(y) \). It must be that \( Y = h^i (Y_*) \) for some \( i \in \mathbb{Z} \). However, note that \( h^3 \) projects to a homeomorphism of \( S \), so up to equivalence we may assume that \( i = 0,1,2 \).

We conclude in three steps. 

\textbf{Step 1.} In this step, suppose \( X = X_1 \). 

If \( Y = Y_* \), then we have \( y \in Y \) which intersects no arc in \( X \cup Y \). Therefore, this does not project to a 1-system, since it would contradict the assumption that \( |J| = 2 \). 

If instead \( Y = h(Y_*) \), we may project this choice of \( X, Y \) to a 1-system on \( S \). This is shown in \Cref{fig:J22 X1 syss a}. 

The situation where \( Y = h^2(Y_*) \) can be obtained from the situation where \( Y = h(Y_*) \) by applying a reflection across the vertical axis, which projects to a homeomorphism of \( S \). In particular, this reflection preserves \( X_1 \). So, this projects to a 1-system equivalent to the one above.

\textbf{Step 2.} In this step we assume that \( X = h(X_1) \). 

If \( Y = Y_* \), we may project this choice of \( X, Y \) to a 1-system on \( S \). This is shown in \Cref{fig:J22 X1 syss b}. 

If instead \( Y = h(Y_*) \), then we have \( h(y) \in Y \) which intersects no arc in \( Y = h(Y_*) \), and no arc in \( X = h(X_1) \). Therefore, this choice does not project to a 1-system.

Finally, if \( Y = h^2(Y_{*}) \), we may project this choice to a 1-system on \( S \). This is shown in \Cref{fig:J22 X1 syss c}.

\begin{figure}
	\centering
	\begin{subfigure}{0.45\textwidth} 
		\centering
		\begin{tikzpicture}[baseline=(current bounding box.center),scale=.5]

			\node at (0,0.5) {b};
			\node at (0,4.5) {a};
			\node at (4.5,0) {b};
			\node at (-4.5,0) {b};

			\draw (0,0) circle (4)
			(0,0) circle (1);
			\draw[fill] (0,1) circle (0.1)
			(0,4) circle (0.1)
			(4,0) circle (0.1)
			(-4,0) circle (0.1);

			\draw (-4,0) to[out=50,in=135] (1.2,1.2) to[out=-45,in=45] (1.2,-1.2) to[out=-135,in=-45] (-1.2,-1.2) to[out=135,in=170] (0,1)
			(0,4) to[out=-60,in=45] (1.4,-1.4) to[out=-135,in=-45] (-1.4,-1.4) to[out=135,in=160] (0,1)
			(4,0) to[out=-150,in=-0] (0,-2.2) to[out=180,in=-90] (-2,-0.4) to[out=90,in=150] (0,1);
			\draw (-4,0) to[out=30,in=140] (0,1);
			\draw (0,4) -- (0,1)
			(4,0) to[out=150,in=20] (0,1)
			(-4,0) to[out=-20,in=180] (0,-1.6) to[out=0,in=-45] (1,0.7) to[out=135,in=10] (0,1);

			\draw[dashed] (4,0) to[out=-120,in=0] (0,-3.1) to[out=180,in=-90] (-3,0) to[out=90,in=-155] (0,4)
			(-4,0) to[out=-60,in=180] (0,-2.9) to[out=0,in=-90] (3,0) to[out=90,in=-25] (0,4)
			(0,4) to[out=-50,in=90] (2.6,0) to[out=-90,in=0] (0,-2.6) to[out=180,in=-90] (-2.6,0) to[out=90,in=-130] (0,4);
		\end{tikzpicture}
		\caption{}
		\label{fig:J22 X1 syss a}
	\end{subfigure}
	\begin{subfigure}{0.45\textwidth} 
		\centering
		\begin{tikzpicture}[xscale=-1,baseline=(current bounding box.center),scale=.5]

			\node at (0,0.5) {b};
			\node at (0,4.5) {a};
			\node at (4.5,0) {b};
			\node at (-4.5,0) {b};

			\draw (0,0) circle (4)
			(0,0) circle (1);
			\draw[fill] (0,1) circle (0.1)
			(0,4) circle (0.1)
			(4,0) circle (0.1)
			(-4,0) circle (0.1);

			\draw (0,4) to[out=-60,in=45] (1,-1) to[out=-135,in=-45] (-1,-1) to[out=135,in=160] (0,1)
			(4,0) to[out=-150,in=-0] (0,-2.2) to[out=180,in=-90] (-2,-0.4) to[out=90,in=150] (0,1)
			(-4,0) to[out=30,in=140] (0,1)
			(0,4) -- (0,1)
			(4,0) to[out=150,in=30] (0,1)
			(-4,0) to[out=-20,in=180] (0,-2.3) to[out=0,in=-90] (2,-0.4) to[out=90,in=20] (0,1)
			(0,4) to[out=-120,in=135] (-1.3,-1.3) to[out=-45,in=-135] (1.3,-1.3) to[out=45,in=10] (0,1);

			\draw[dashed] (4,0) to[out=-120,in=0] (0,-3.1) to[out=180,in=-90] (-3,0) to[out=90,in=-155] (0,4)
			(4,0) to[out=120,in=0] (0,2.6) to[out=180,in=90] (-2.6,0) to[out=-90,in=180] (0,-2.6) to[out=0,in=-120] (4,0)
			(4,0) to[out=120,in=0] (0,3) to[out=180,in=60] (-4,0);
		\end{tikzpicture}
		\caption{}
		\label{fig:J22 X1 syss b}
	\end{subfigure}
	\begin{subfigure}{0.45\textwidth} 
		\centering
		\begin{tikzpicture}[xscale=-1,baseline=(current bounding box.center),scale=.5]

			\node at (0,0.5) {b};
			\node at (0,4.5) {a};
			\node at (4.5,0) {b};
			\node at (-4.5,0) {b};

			\draw (0,0) circle (4)
			(0,0) circle (1);
			\draw[fill] (0,1) circle (0.1)
			(0,4) circle (0.1)
			(4,0) circle (0.1)
			(-4,0) circle (0.1);

			\draw (-4,0) to[out=50,in=135] (1.2,1.2) to[out=-45,in=45] (1.2,-1.2) to[out=-135,in=-45] (-1.2,-1.2) to[out=135,in=170] (0,1)
			(0,4) to[out=-60,in=45] (1.4,-1.4) to[out=-135,in=-45] (-1.4,-1.4) to[out=135,in=160] (0,1)
			(4,0) to[out=-150,in=-0] (0,-2.2) to[out=180,in=-90] (-2,-0.4) to[out=90,in=150] (0,1);
			\draw (-4,0) to[out=30,in=140] (0,1);
			\draw (0,4) -- (0,1)
			(4,0) to[out=150,in=20] (0,1)
			(-4,0) to[out=-20,in=180] (0,-1.6) to[out=0,in=-45] (1,0.7) to[out=135,in=10] (0,1);

			\draw[dashed] (4,0) to[out=-120,in=0] (0,-3.1) to[out=180,in=-90] (-3,0) to[out=90,in=-155] (0,4)
			(4,0) to[out=120,in=0] (0,2.6) to[out=180,in=90] (-2.6,0) to[out=-90,in=180] (0,-2.6) to[out=0,in=-120] (4,0)
			(4,0) to[out=120,in=0] (0,3) to[out=180,in=60] (-4,0);
		\end{tikzpicture}
		\caption{}
		\label{fig:J22 X1 syss c}
	\end{subfigure}

	\caption{}
	\label{fig:J22 X1 syss}
\end{figure}
\textbf{Step 3.} In this step we assume that \( X = X_2 \). 

If \( Y = Y_* \), we may lift this choice of \( X, Y \) to a 1-system on \( S \). This is shown in \Cref{fig:J22 X2 syss a}. 

If \( Y = h(Y_*) \), we may lift this choice to another 1-system on \( S \). This is shown in \Cref{fig:J22 X2 syss b}. 

The situation where \( Y = h^2(Y_*) \) can be obtained from the situation where \( Y = h(Y_*) \) by reflection, as in step 1.

\begin{figure}
	\centering
	\begin{subfigure}{0.45\textwidth}
		\centering
		\begin{tikzpicture}[baseline=(current bounding box.center),scale=.5]

			\node at (0,0.5) {b};
			\node at (0,4.5) {a};
			\node at (4.5,0) {b};
			\node at (-4.5,0) {b};

			\draw (0,0) circle (4)
			(0,0) circle (1);
			\draw[fill] (0,1) circle (0.1)
			(0,4) circle (0.1)
			(4,0) circle (0.1)
			(-4,0) circle (0.1);

			\draw (0,4) to[out=-60,in=45] (1,-1) to[out=-135,in=-45] (-1,-1) to[out=135,in=160] (0,1)
			(4,0) to[out=-150,in=-0] (0,-2.2) to[out=180,in=-90] (-2,-0.4) to[out=90,in=150] (0,1)
			(-4,0) to[out=30,in=140] (0,1)
			(0,4) -- (0,1)
			(4,0) to[out=150,in=30] (0,1)
			(-4,0) to[out=-20,in=180] (0,-2.3) to[out=0,in=-90] (2,-0.4) to[out=90,in=20] (0,1)
			(0,4) to[out=-120,in=135] (-1.3,-1.3) to[out=-45,in=-135] (1.3,-1.3) to[out=45,in=10] (0,1);

			\draw[dashed] (4,0) to[out=-120,in=0] (0,-3.1) to[out=180,in=-90] (-3,0) to[out=90,in=-155] (0,4)
			(-4,0) to[out=-60,in=180] (0,-2.9) to[out=0,in=-90] (3,0) to[out=90,in=-25] (0,4)
			(4,0) to[out=120,in=0] (0,3) to[out=180,in=60] (-4,0);
		\end{tikzpicture}
		\caption{}
		\label{fig:J22 X2 syss a}
	\end{subfigure}
	\begin{subfigure}{0.45\textwidth}
		\centering
		\begin{tikzpicture}[baseline=(current bounding box.center),scale=.5]

			\node at (0,0.5) {b};
			\node at (0,4.5) {a};
			\node at (4.5,0) {b};
			\node at (-4.5,0) {b};

			\draw (0,0) circle (4)
			(0,0) circle (1);
			\draw[fill] (0,1) circle (0.1)
			(0,4) circle (0.1)
			(4,0) circle (0.1)
			(-4,0) circle (0.1);

			\draw (-4,0) to[out=50,in=135] (1.2,1.2) to[out=-45,in=45] (1.2,-1.2) to[out=-135,in=-45] (-1.2,-1.2) to[out=135,in=170] (0,1)
			(0,4) to[out=-60,in=45] (1.4,-1.4) to[out=-135,in=-45] (-1.4,-1.4) to[out=135,in=160] (0,1)
			(4,0) to[out=-150,in=-0] (0,-2.2) to[out=180,in=-90] (-2,-0.4) to[out=90,in=150] (0,1);
			\draw (-4,0) to[out=30,in=140] (0,1);
			\draw (0,4) -- (0,1)
			(4,0) to[out=150,in=20] (0,1)
			(-4,0) to[out=-20,in=180] (0,-1.6) to[out=0,in=-45] (1,0.7) to[out=135,in=10] (0,1);

			\draw[dashed] (4,0) to[out=-120,in=0] (0,-3.1) to[out=180,in=-90] (-3,0) to[out=90,in=-155] (0,4)
			(-4,0) to[out=-60,in=180] (0,-2.9) to[out=0,in=-90] (3,0) to[out=90,in=-25] (0,4)
			(4,0) to[out=120,in=0] (0,3) to[out=180,in=60] (-4,0);
		\end{tikzpicture}
		\caption{}
		\label{fig:J22 X2 syss b}
	\end{subfigure}
    \caption{}
    \label{fig:J22 X2 syss}
\end{figure}

\subsection{$J$ contains two non-loops}

In this case, the surface \( S \setminus J \) is an annulus with four marked points, two on each boundary component. This is shown in \Cref{fig:J23 intro a}.

\begin{figure}
	\centering
	\begin{subfigure}[b]{0.6\textwidth}
		\centering
		\begin{tikzpicture}[baseline=(current bounding box.center)]
			\draw (0,0) -- (2,0) -- (2,2) -- (0,2) -- (0,0);

			\draw[fill] (0.7,1.3) circle (0.06)
			(1.3,0.7) circle (0.06);

			\draw (1,0) -- (1.3,0.7) -- (0.7,1.3) -- (1,2);

			\node at (1.5,0.7) {\( a \)};
			\node at (1,1.3) {$b$};
		\end{tikzpicture}
		\qquad
		\begin{tikzpicture}[baseline=(current bounding box.center),scale=.4]
			\draw (0,0) circle (4)
			(0,0) circle (1);
			\draw[fill] (0,1) circle (0.1)
			(0,4) circle (0.1)
			(0,-1) circle (0.1)
			(0,-4) circle (0.1);

			\node at (0,1.6) {$a$};
			\node at (0,3.4) {$a$};
			\node at (0,-1.6) {$b$};
			\node at (0,-3.4) {$b$};
		\end{tikzpicture}
		\caption{}
		\label{fig:J23 intro a}
	\end{subfigure}

	\begin{subfigure}[b]{0.3\textwidth}
		\centering
		\begin{tikzpicture}[scale=.4]
			\draw[blue] (0,0.5) circle (3.5);
			\draw[red] (0,-0.5) circle (1.5);

			\draw[blue] (0,-0.5) circle (3.5);
			\draw[red] (0,0.5) circle (1.5);

			\draw (0,0) circle (4)
			(0,0) circle (1);
			\draw[fill] (0,1) circle (0.1)
			(0,4) circle (0.1)
			(0,-1) circle (0.1)
			(0,-4) circle (0.1);
		\end{tikzpicture}
		\caption{}
		\label{fig:J23 intro b}
	\end{subfigure}
	\begin{subfigure}[b]{0.3\textwidth}
		\begin{tikzpicture}[scale=.4]
			\draw[red] (0,1) -- (0,4) to[out=-45,in=-20] (0,-1) -- (0,-4);

			\node[red] at (-0.4,2.5) {\( c_0 \)};
			\node[red] at (1.8,1.6) {\( c_{1/2} \)};
			\node[red] at (-0.4,-2.5) {\( d_0 \)};

			\draw (0,0) circle (4)
			(0,0) circle (1);
			\draw[fill] (0,1) circle (0.1)
			(0,4) circle (0.1)
			(0,-1) circle (0.1)
			(0,-4) circle (0.1);
		\end{tikzpicture}
		\caption{}
		\label{fig:J23 intro c}
	\end{subfigure}
	\caption{}
	\label{fig:J23 intro}
\end{figure}

Write \( \mathcal{A} \setminus J = X \sqcup Y \), where \( X \) is the set of arcs which have both endpoints on the same boundary component of \( S \setminus J \), and \( Y \) is the set of arcs which have one endpoint on each boundary component.

There are only four arcs which have both endpoints on the same boundary component. These are shown in \Cref{fig:J23 intro b}. We see that the two arcs shown in red intersect twice, and the two arcs shown in blue intersect twice. This means that \( |X| \le 2 \): it may contain at most one of the blue arcs and one of the red arcs.

Now consider the arcs in \( Y \).  
We define \( h \) to be the half Dehn twist about the inner boundary component. Define \( c_0, d_0 \) to be the arcs shown in \Cref{fig:J23 intro c}. Define \( c_{k} = h^{2k} (c_0) \) for \( k \in \frac{1}{2} \mathbb{Z} \) and define \( d_k \) similarly. Then we see that 
\begin{equation}
	i(c_k, c_\ell) = i(d_k, d_\ell) = \lceil |k - \ell| \rceil - 1
	\label{eq:J23 crosses 1}
\end{equation}
\begin{equation}
	i(c_k, d_{\ell}) = \lfloor |k - \ell| \rfloor
	\label{eq:J23 crosses 2}
\end{equation}

In particular, note that every arc in \( Y \) is of the form \( c_k \) or \( d_k \) for some \( k \). Write \( A = \left\{ k \in \frac{1}{2}\mathbb{Z} : c_k \in Y \right\} \) and \( B = \left\{ k \in \frac{1}{2}\mathbb{Z} : d_k \in Y \right\} \).

Since \( |X| \le 2 \), we know \( |Y| \ge 8 \). Therefore we may assume that \( |A| \ge 4 \). In addition, up to applying \( h \), we may assume that \( \min (A) = -1 \). Then, from \Cref{eq:J23 crosses 1}, we see that \( \max( A) \le 1 \). 

\begin{figure}
	\centering
	\begin{subfigure}[b]{0.3\textwidth}
		\centering
		\begin{tikzpicture}[scale=0.4]

			\draw[red]
			(0,4) -- (0,1) 
			(0,4) to[out=-45,in=-20] (0,-1)
			(0,4) to[out=-135,in=-160] (0,-1)
			(0,4) to[out=-30,in=30] (0.9,-1.5) to[out=-150,in=-70] (-1.6,-0.5) to[out=110,in=-110] (-1.6,0.4) to[out=70,in=160] (0,1)
			(0,4) to[out=-150,in=150] (-0.9,-1.5) to[out=-30,in=-110] (1.6,-0.5) to[out=70,in=-70] (1.6,0.4) to[out=110,in=20] (0,1);

			\draw[green] (0,-4) -- (0,-1)
			(0,-4) to[out=45,in=-80] (1.7,0.5) to[out=100,in=30] (0,1)
			(0,-4) to[out=135,in=-100] (-1.7,0.5) to[out=80,in=150] (0,1);

			\draw (0,0) circle (4)
			(0,0) circle (1);
			\draw[fill] (0,1) circle (0.1)
			(0,4) circle (0.1)
			(0,-1) circle (0.1)
			(0,-4) circle (0.1);

		\end{tikzpicture}
		\caption{}
		\label{fig:J23 second a}
	\end{subfigure}
	\begin{subfigure}[b]{0.3\textwidth}
		\centering
		\begin{tikzpicture}[scale=0.4]

			\draw[red] 
			(0,4) to[out=-45,in=20] (0.4,-1.7) to[out=-160,in=-100] (-1,0)
			(0,4) to[out=-60,in=75] (1,0)
			(0,4) to[out=-120,in=105] (-1,0)
			(0,4) to[out=-135,in=160] (-0.4,-1.7) to[out=-20,in=-80] (1,0);

			\draw[green,yscale=-1,xscale=1] 
			(0,4) to[out=-45,in=20] (0.4,-1.7) to[out=-160,in=-100] (-1,0)
			(0,4) to[out=-60,in=75] (1,0)
			(0,4) to[out=-120,in=105] (-1,0)
			(0,4) to[out=-135,in=160] (-0.4,-1.7) to[out=-20,in=-80] (1,0);

			\node[red] at(0.35,2.3) {\( c_0 \)};
			\node[red] at(2.0,2.3) {\( c_{1/2} \)};
			\node[green] at (-0.35,-2.3) {\( d_0 \)};

			\draw (0,0) circle (1) circle (4);
			\draw[fill]
			(1,0) circle (0.1)
			(-1,0) circle (0.1)
			(0,4) circle (0.1)
			(0,-4) circle (0.1);

		\end{tikzpicture}
		\caption{}
		\label{fig:J23 second b}
	\end{subfigure}
	\caption{The two possibilities for \( Y \). A homeomorphism has been applied to the figure on the right to show certain symmetries.}
	\label{fig:J23 second}
\end{figure}
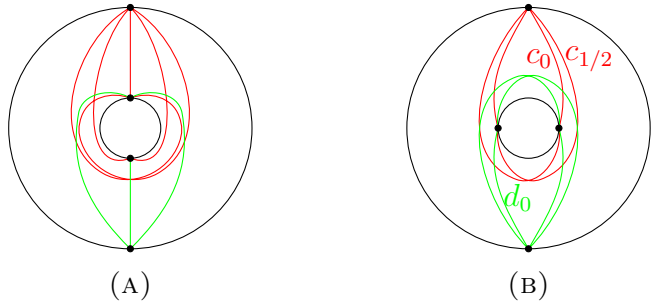

If we have \( c_1 \in Y \), then from \Cref{eq:J23 crosses 2} we have \( B \subset \left\{ -\frac{1}{2}, 0, \frac{1}{2} \right\} \). Since \( |Y| \ge 8 \), in this case it must be that 
\[ Y = \left\{ c_{-1}, c_{-1/2}, c_{0}, c_{1/2}, c_1, d_{-1/2}, d_{0}, d_{1/2} \right\}. \]
This is shown in \Cref{fig:J23 second a}.

If instead \( c_1 \notin Y \), since \( |A| \ge 4 \), it must be that \( A = \left\{ -1, -1/2, 1, {1/2} \right\} \). Then we must have 
\[ Y = \left\{ c_{-1}, c_{-1/2}, c_{0}, c_{1/2}, d_{-1}, d_{-1/2}, d_{0}, d_{1/2} \right\}. \]
This is shown in \Cref{fig:J23 second b}.

For both possibilities we have \( |Y| = 8 \), and so we must have \( |X| = 2 \).

We conclude in two steps.

\textbf{Step 1.} In this step we assume that \( Y \) is the system shown in \Cref{fig:J23 second a}. 

Note that there is one arc in \( Y \) which is disjoint from all the other arcs in \( Y \). So, this arc must intersect at least one of the arcs in \( X \). 

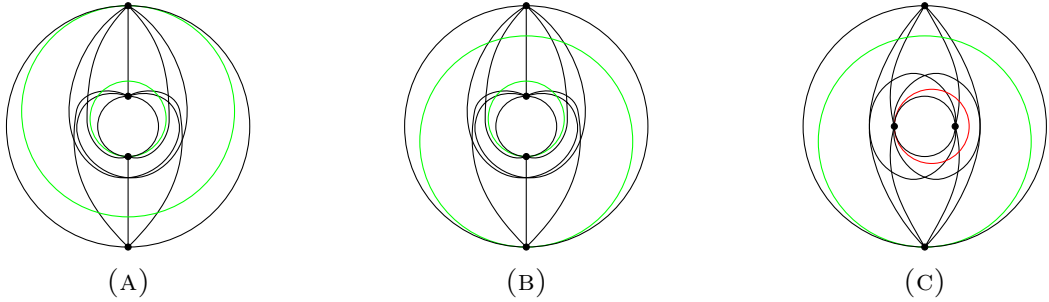
\begin{figure}
	\centering
	\begin{subfigure}[b]{0.3\textwidth}
		\centering
		\begin{tikzpicture}[scale=0.4]

			\draw
			(0,4) -- (0,1) ;
			\draw
	(0,4) to[out=-45,in=90] (1.35,0.2) to[out=-90,in=-20] (0,-1)
	(0,4) to[out=-135,in=90] (-1.35,0.2) to[out=-90,in=-160] (0,-1)
			(0,4) to[out=-30,in=30] (0.9,-1.5) to[out=-150,in=-70] (-1.6,-0.5) to[out=110,in=-110] (-1.6,0.4) to[out=70,in=160] (0,1)
			(0,4) to[out=-150,in=150] (-0.9,-1.5) to[out=-30,in=-110] (1.6,-0.5) to[out=70,in=-70] (1.6,0.4) to[out=110,in=20] (0,1);

			\draw (0,-4) -- (0,-1)
			(0,-4) to[out=45,in=-80] (1.7,0.5) to[out=100,in=30] (0,1)
			(0,-4) to[out=135,in=-100] (-1.7,0.5) to[out=80,in=150] (0,1);

			\draw[green] (0,0.25) circle (1.25)
			(0,0.5) circle (3.5);

			\draw (0,0) circle (4)
			(0,0) circle (1);
			\draw[fill] (0,1) circle (0.1)
			(0,4) circle (0.1)
			(0,-1) circle (0.1)
			(0,-4) circle (0.1);

		\end{tikzpicture}
		\caption{}
		\label{fig:J23 third a}
	\end{subfigure}
	\begin{subfigure}[b]{0.3\textwidth}
		\centering
		\begin{tikzpicture}[scale=0.4]

			\draw
			(0,4) -- (0,1) ;
			\draw
	(0,4) to[out=-45,in=90] (1.35,0.2) to[out=-90,in=-20] (0,-1)
	(0,4) to[out=-135,in=90] (-1.35,0.2) to[out=-90,in=-160] (0,-1)
			(0,4) to[out=-30,in=30] (0.9,-1.5) to[out=-150,in=-70] (-1.6,-0.5) to[out=110,in=-110] (-1.6,0.4) to[out=70,in=160] (0,1)
			(0,4) to[out=-150,in=150] (-0.9,-1.5) to[out=-30,in=-110] (1.6,-0.5) to[out=70,in=-70] (1.6,0.4) to[out=110,in=20] (0,1);

			\draw (0,-4) -- (0,-1)
			(0,-4) to[out=45,in=-80] (1.7,0.5) to[out=100,in=30] (0,1)
			(0,-4) to[out=135,in=-100] (-1.7,0.5) to[out=80,in=150] (0,1);

			\draw[green] (0,0.25) circle (1.25)
			(0,-0.5) circle (3.5);

			\draw (0,0) circle (4)
			(0,0) circle (1);
			\draw[fill] (0,1) circle (0.1)
			(0,4) circle (0.1)
			(0,-1) circle (0.1)
			(0,-4) circle (0.1);

		\end{tikzpicture}
		\caption{}
		\label{fig:J23 third b}
	\end{subfigure}
	\begin{subfigure}[b]{0.3\textwidth}
		\centering
		\begin{tikzpicture}[scale=0.4]

			\draw[] 
			(0,4) to[out=-45,in=20] (0.7,-1.7) to[out=-160,in=-100] (-1,0)
			(0,4) to[out=-60,in=75] (1,0)
			(0,4) to[out=-120,in=105] (-1,0)
			(0,4) to[out=-135,in=160] (-0.7,-1.7) to[out=-20,in=-80] (1,0);

			\draw[yscale=-1,xscale=1] (0,0)
			(0,4) to[out=-45,in=20] (0.7,-1.7) to[out=-160,in=-100] (-1,0)
			(0,4) to[out=-60,in=75] (1,0)
			(0,4) to[out=-120,in=105] (-1,0)
			(0,4) to[out=-135,in=160] (-0.7,-1.7) to[out=-20,in=-80] (1,0);

				\draw[green] (0,-0.5) circle (3.5);
				\draw[red] (0.23,0) circle (1.23);

			\draw (0,0) circle (1) circle (4);
			\draw[fill]
			(1,0) circle (0.1)
			(-1,0) circle (0.1)
			(0,4) circle (0.1)
			(0,-4) circle (0.1);

		\end{tikzpicture}
		\caption{}
		\label{fig:J23 third c}
	\end{subfigure}
	\caption{Three systems on \( S \setminus J \). Each projects to two 1-systems on \( (S, P) \), by precomposing the almost embedding with id or \( h \).}
	\label{fig:J23 third}
\end{figure}

There are three possibilities for \( X \). Two of them result in equivalent systems on \( S \setminus J \). So, there are two non-equivalent 1-systems on \( S \setminus J \). These are shown in \Cref{fig:J23 third a,fig:J23 third b}. Each of these systems projects to a 1-system on \( (S, P) \), and after applying \( h \) to each system we obtain a non-equivalent 1-system on \( (S, P) \). 

So, in this step we obtain four distinct 1-systems. 

\textbf{Step 2.} In this step we assume that \( Y \) is the system shown in \Cref{fig:J23 second b}. 

Note that this configuration is invariant under reflection across the horizontal axis. So, we may assume that the arc shown in green in \Cref{fig:J23 third b} is included in \( X \subset \mathcal{A} \). The configuration is also invariant under reflection across the vertical axis, so we may assume that the arc shown in red in \Cref{fig:J23 third c} is included in \( X \subset \mathcal{A} \). 

The system shown in \Cref{fig:J23 third c} projects to a 1-system on \( (S, P) \), and after applying \( h \) we obtain a non-equivalent 1-system on \( (S, P) \). 

So, in this step we obtain two distinct 1-systems.

\section{$|J| = 1$}
\label{sec:J1}

Write $J = \left\{ u \right\}$. There are two cases: $u$ is a loop or $u$ is a non-loop arc. 

\subsection{$u$ is a loop}

\begin{figure}
	\centering
	\begin{subfigure}[b]{0.2\textwidth}
		\centering
%

		\begin{tikzpicture}
			\draw (0,0) circle (1.2)
			circle (.3);

			\draw[fill] (0,0.3) circle (0.06)
			(0,-0.75) circle (0.06)
			(0,1.2) circle (0.06);

			\draw[dashed] (0,0) circle (0.55) circle (0.95);
		\end{tikzpicture}
		\caption{}
		\label{fig:J11 enum a}
	\end{subfigure}
	\begin{subfigure}[b]{0.7\textwidth}
		\centering
		\begin{tikzpicture}
			\draw (0,0) circle (1.2)
			circle (.3);

			\draw[fill] (0,0.3) circle (0.06)
			(0,-0.75) circle (0.06)
			(0,1.2) circle (0.06);

			\draw (0,0.3) -- (0,1.2);

			\node at (0.4,0.8) {$v$};
		\end{tikzpicture}
		\begin{tikzpicture}
			\draw (0,0) circle (1.2)
			circle (.3);

			\draw[fill] (0,0.3) circle (0.06)
			(0,-0.75) circle (0.06)
			(0,1.2) circle (0.06);

			\draw (0,0.3) to[out=165,in=80] (-0.5,-0.2) to[out=-100,in=155] (0,-0.75);

			\node at (0.3,0.5) {\( w \)};
		\end{tikzpicture}
		
		\begin{tikzpicture}
			\draw (0,0) circle (1.2)
			circle (.3);

			\draw[fill] (0,0.3) circle (0.06)
			(0,-0.75) circle (0.06)
			(0,1.2) circle (0.06);
			
			\draw (0,-.3) circle (0.6);
			\node at (-0.8,0) {\( x \)};
		\end{tikzpicture}
		\begin{tikzpicture}
			\draw (0,0) circle (1.2)
			circle (.3);

			\draw[fill] (0,0.3) circle (0.06)
			(0,-0.75) circle (0.06)
			(0,1.2) circle (0.06);

			\draw (0,0) circle (-0.75);
			\node at (0.9,0) {\( y \)};
		\end{tikzpicture}
		\begin{tikzpicture}
			\draw (0,0) circle (1.2)
			circle (.3);

			\draw[fill] (0,0.3) circle (0.06)
			(0,-0.75) circle (0.06)
			(0,1.2) circle (0.06);

			\draw (0,-0.6) arc (-90:90:0.45)
			(0,-0.9) arc (-90:90:0.6)
			(0,-0.6) arc (90:270:0.15);

			\node at (0.4,0.7) {\( z \)};
		\end{tikzpicture}
		\caption{}
		\label{fig:J11 enum b}
	\end{subfigure}
	\caption{}
	\label{fig:J11 enum}
\end{figure}

In this case, $S \setminus J$ is an annulus with three marked points: one on each boundary component and one in the interior. This is shown in Figure~\ref{fig:J11 enum a}.

Consider all possible arcs on $S \setminus J$ up to homeomorphism. There are five, labelled $v, w, x, y, z$ as shown in Figure~\ref{fig:J11 enum b}. We may write $\mathcal{A}' = V \sqcup W \sqcup X \sqcup Y \sqcup Z$, where $V$ contains exactly the arcs in $\mathcal{A}'$ which can be taken to $v$ by a homeomorphism of $S \setminus J$, and similarly for $W, X, Y, Z$. 

We make the following observations:

\begin{itemize}
	\item
		\( Z = \emptyset \).

		To see this, suppose for contradiction \( Z \ne \emptyset \), so $|Z| \ge 1$. Without loss of generality suppose $z \in Z$. Since \( z \) bounds a disk with one puncture on the interior, we may apply \Cref{cor:help ones} and conclude that the internal arc $g$ is in $J \subset \mathcal{A}$. This contradicts our assumption that $|J| = 1$. So it must be that \( Z = \emptyset \). 

	\item
		$|X| \le 1$

		This follows from the fact that any two distinct elements of $X$ intersect at least twice.
		

	\item
		$|Y| \le 1$

		This follows because $y$ is invariant under any homeomorphism of $S \setminus J$. So, either $Y = \left\{ y \right\}$ or $Y = \emptyset$. 

	\item
		$|X| + |Y| \le 1$

		This follows because $x$ and $y$ intersect twice.

	\item 
		$|V| + |W| \ge 10$

		This follows from the fact that \( |\mathcal{A}'| = 11 \) together with the above observations. 
		
\end{itemize}

\begin{figure}
	\centering
	\begin{subfigure}[b]{0.3\textwidth}
		\centering
		\begin{tikzpicture}
			\draw (0,0) circle (1.2)
			circle (.3);

			\draw[fill] (0,0.3) circle (0.06)
			(0,-0.75) circle (0.06)
			(0,1.2) circle (0.06);

			\draw (0,0.3) -- (0,1.2);
		\end{tikzpicture}
		\caption{$v = v_{00}$}
	\end{subfigure}
	\begin{subfigure}[b]{0.3\textwidth}
		\centering
		\begin{tikzpicture}
			\draw (0,0) circle (1.2)
			circle (.3);

			\draw[fill] (0,0.3) circle (0.06)
			(0,-0.75) circle (0.06)
			(0,1.2) circle (0.06);

			\draw (0,0.3) to[out=160,in=90] (-0.5,0) to[out=-90,in=180] (0,-0.6) to[out=0,in=-90] (0.9,0) to[out=90,in=-40] (0,1.2);
		\end{tikzpicture}
		\caption{\( v_{10} \)}
	\end{subfigure}
	\begin{subfigure}[b]{0.3\textwidth}
		\centering
		\begin{tikzpicture}
			\draw (0,0) circle (1.2)
			circle (.3);

			\draw[fill] (0,0.3) circle (0.06)
			(0,-0.75) circle (0.06)
			(0,1.2) circle (0.06);

			\draw (0,0.3) to[out=160,in=90] (-0.5,0) to[out=-90,in=180] (0,-0.9) to[out=0,in=-90] (0.9,0) to[out=90,in=-40] (0,1.2);
		\end{tikzpicture}
		\caption{\( v_{01} \)}
	\end{subfigure}
	\caption{}
	\label{fig:a01 a10}
\end{figure}

Write $v_{ij}$ for the arc obtained from $v$ by the $i$th power of the Dehn twist about the inner dashed curve of \Cref{fig:J11 enum a} and the $j$th power of the Dehn twist about the outer dashed curve of \Cref{fig:J11 enum a}. See \Cref{fig:a01 a10} for an example. Write $w_{k + 1/2}$ for the arc obtained from $w$ by the $k$th power of the Dehn twist about the inner dashed curve of \Cref{fig:J11 enum a}, and similarly write \( x_k \) for the arc obtained from \( x \) by the \( k \)th power of the Dehn twist about the inner dashed curve of \Cref{fig:J11 enum a}.

Let $h: S\setminus J \to S \setminus J$ be the orientation reversing homeomorphism given by the inversion about the core curve. Note that every arc in $V$ is $v_{jk}$ for some $j,k$; in particular \( h(v_{jk}) = v_{-j,-k} \). We write \( v_{j,k} \) when the indexing is not otherwise clear. We write $W = W_+ \sqcup W_-$, where every arc in $W_+$ is of the form $w_\ell$ for some $\ell$, and every arc in $W_-$ is of the form $h(w_\ell)$ for some $\ell$.

First, we consider the arcs in \( V \). We have the following formula for the intersection number of \( v = v_{00}, v_{jk} \) for \( (j,k) \ne (0,0) \). This is:
\begin{equation}
	i(v, v_{jk}) = \begin{cases}
		|j| + |k| - 2 & \text{when } j\cdot k < 0 \\
		|j| + |k| - 1 & \text{otherwise}
	\end{cases}
	\label{eq:vij ints}
\end{equation}

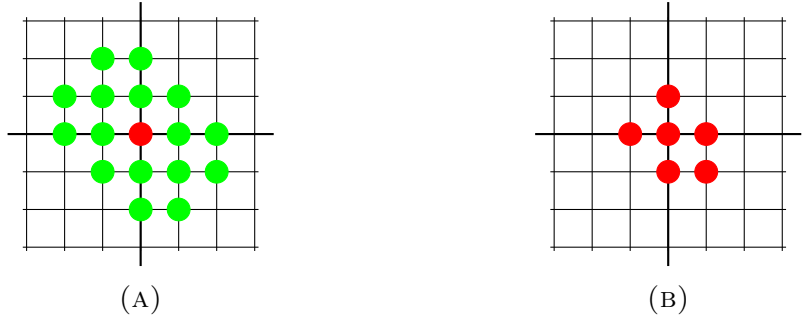
\begin{figure}
	\centering
	\begin{subfigure}[b]{0.4\textwidth}
		\centering
		\begin{tikzpicture}[scale=0.5]

			\foreach\x in {-3,-2,...,3}{
				\draw (\x,-3.1) -- (\x,3.1);
				\draw (-3.1,\x) -- (3.1,\x);
			}
			\draw[thick] (-3.5,0) -- (3.5,0)
			(0,-3.5) -- (0,3.5);

			\draw[fill,red] (0,0) circle (0.3);

			\draw[fill,green] (1,0) circle (0.3)
			(2,0) circle (0.3)
			(-1,0) circle (0.3)
			(-2,0) circle (0.3)
			(0,1) circle (0.3)
			(1,1) circle (0.3)
			(-1,1) circle (0.3)
			(-2,1) circle (0.3)
			(0,2) circle (0.3)
			(-1,2) circle (0.3)
			(1,-1) circle (0.3)
			(2,-1) circle (0.3)
			(0,-1) circle (0.3)
			(-1,-1) circle (0.3)
			(1,-2) circle (0.3)
			(0,-2) circle (0.3);
		\end{tikzpicture}
		\caption{}
		\label{fig:big z squared a}
	\end{subfigure}
	\begin{subfigure}[b]{0.4\textwidth}
		\centering
		\begin{tikzpicture}[scale=0.5]
			\foreach\x in {-3,-2,...,3}{
				\draw (\x,-3.1) -- (\x,3.1);
				\draw (-3.1,\x) -- (3.1,\x);
			}
			\draw[thick] (-3.5,0) -- (3.5,0)
			(0,-3.5) -- (0,3.5);

			\draw[fill, red] (0,0) circle (0.3)
			(1,0) circle (0.3)
			(-1,0) circle (0.3)
			(1,-1) circle (0.3)
			(0,-1) circle (0.3)
			(0,1) circle (0.3);
		\end{tikzpicture}
		\caption{}
		\label{fig:big z squared b}
	\end{subfigure}
	\caption{The vertex \( (j, k) \) corresponds to the arc \( v_{jk} \)}
	\label{fig:big z squared}
\end{figure}

Suppose that \( V \ne \emptyset \). Without loss of generality, we may assume that \( v_{00} \in V \). We note that, based on \Cref{eq:vij ints} and the fact that \( \mathcal{A} \) is a 1-system, only certain other arcs of the form \( v_{jk} \) may be included in \( V \). In \Cref{fig:big z squared a}, the arc \( v = v_{00} \) is shown in red, while the arcs that intersect \( v \) at most once are shown in green. 

Let \( V_* \) be the set of arcs represented by vertices shown in red in \Cref{fig:big z squared b}. 

\begin{claim}
	Up to shifts and reflection across the line \( y = x \), \( V \) must be a subset of \( V_* \). 
\end{claim}
\begin{proof}
	First, suppose \( v_{00}, v_{02} \in V \). Then, using \Cref{eq:vij ints}, we determine that the only other arcs which may be in \( V \) are the ones represented by vertices shown in green in \Cref{fig:smol z squared a}. Since \( v_{1,0}, v_{-1,2} \) intersect twice, they cannot both be in \( V \). If \( v_{10} \notin V \), then \( V \) is contained in \( V_* \) reflected, then shifted up by one. If \( v_{-1,2} \notin V \), then \( V \) is contained in \( V_* \) shifted up one. So, we are done.

	Now suppose \( v_{00}, v_{-1,2} \in V \). If either \( v_{02} \) or \( v_{-1,0} \) are in \( V \), we are done as above. So the only other arcs which may be in \( V \) are the ones represented by vertices shown in green in \Cref{fig:smol z squared b}. Since \( v_{1,1}, v_{-2,1} \) intersect twice, they cannot both be in \( V \). If \( v_{10} \notin V \), then \( V \) is contained in \( V_* \) shifted up one and left one. If \( v_{-1,2} \notin V \), then \( V \) is contained in \( V_* \) reflected, then shifted up one. So, we are done.

	Finally, suppose neither of the above two configurations occurs. Then any two arcs differ by at most one in the vertical direction and at most one in the horizontal direction. Such a configuration must be a subset of the arcs \( \left\{ v_{00}, v_{10}, v_{0,-1}, v_{1,-1} \right\} \subset V_* \). 
\end{proof}

\begin{figure}
	\centering
	\begin{subfigure}[b]{0.4\textwidth}
		\centering
		\begin{tikzpicture}[scale=0.5]

			\foreach\x in {-3,-2,...,3}{
				\draw (\x,-3.1) -- (\x,3.1);
				\draw (-3.1,\x) -- (3.1,\x);
			}
			\draw[thick] (-3.5,0) -- (3.5,0)
			(0,-3.5) -- (0,3.5);

			\draw[fill,red] (0,0) circle (0.3)
			(0,2) circle (0.3);

			\draw[fill,green] (1,0) circle (0.3)
			(0,1) circle (0.3)
			(1,1) circle (0.3)
			(-1,1) circle (0.3)
			(-1,2) circle (0.3);
		\end{tikzpicture}
		\caption{}
		\label{fig:smol z squared a}
	\end{subfigure}
	\begin{subfigure}[b]{0.4\textwidth}
		\centering
		\begin{tikzpicture}[scale=0.5]
			\draw[fill,color=white!70!blue] (0,2) circle (0.3)
			(-1,0) circle (0.3);

			\foreach\x in {-3,-2,...,3}{
				\draw (\x,-3.1) -- (\x,3.1);
				\draw (-3.1,\x) -- (3.1,\x);
			}
			\draw[thick] (-3.5,0) -- (3.5,0)
			(0,-3.5) -- (0,3.5);

			\draw[fill,red] (0,0) circle (0.3)
			(-1,2) circle (0.3);

			\draw[fill,green] 
			(0,1) circle (0.3)
			(1,1) circle (0.3)
			(-1,1) circle (0.3)
			(-2,1) circle (0.3);

		\end{tikzpicture}
		\caption{}
		\label{fig:smol z squared b}
	\end{subfigure}
	\caption{}
	\label{fig:smol z squared}
\end{figure}

Next, we consider arcs in \( W \). We have the following intersection formulas for intersection:
\begin{equation*}
	i(w_j, w_k) = |j-k| - 1\\
\end{equation*}
\begin{equation*}
	i(w_j, h(w_k)) = 0\\
\end{equation*}
From this we can see that $|W_+| \le 3$, and $|W_-| \le 3$. This also allows us to conclude that \( |V| \ge 4 \), since we know that \( |V| + |W| \ge 10 \). 

We also have the following formula:
\begin{equation}
	v_{jk}, w_\ell = \lfloor |j - \ell| \rfloor
	\label{eq:vjk wl ints}
\end{equation}
We obtain a similar formula for the intersection number of $v_{jk}, h(w_\ell)$ by noting that $h^2 = \text{id}$ and $h(v_{jk}) = v_{-k-j}$. 

We now conclude in three steps.

\textbf{Step 1.} In this step, we assume that \( |V| = 6 \). 

By the analysis above, we may assume that:
\[ V = \left\{ v_{00}, v_{-1,0}, v_{0,1}, v_{0,-1}, v_{1,0}, v_{1,-1} \right\} \]
From \Cref{eq:vjk wl ints}, we see that \( W_+ \subset \left\{ w_{-1/2}, w_{1/2} \right\} \) and \( W_- \subset \left\{ h(w_{-1/2}), h(w_{1/2}) \right\} \). Since we know that \( |V| + |W| \ge 10 \), it must be that these containments are equalities. We see that every arc in \( V \cup W \) intersects at least one other arc in \( V \cup W \), except \( v_{00} \). In order to complete this to a 1-system, it must be that either \( |X| = 1 \) or \( |Y| = 1 \). 

First, suppose \( |X| = 1 \). We may assume that \( X = \left\{ x_j \right\} \) for some $j$ since, for the \( V, W \) we have determined, we have \( h(V \cup W) = V \cup W \). 

We present one final formula: 
\begin{equation*}
	i(x_j, w_k) = 2 \lfloor |j-k| \rfloor
\end{equation*}
Using this formula, and the fact that \( w_{-1/2}, w_{1/2} \in W \), the only possible arc that could be in \( X \) is the arc \( x_0 \). However, \( x_0 \) does not intersect \( v_{00} \). So, if \( X \cup Y = \left\{ x_0 \right\} \), this would contradict our assumption that \( |J| = 1 \). 

On the other hand, if \( X = \emptyset \) and \( Y = \left\{ y \right\} \), we obtain a 1-system as shown in \Cref{fig:J11 V6}.

\begin{figure}
	\centering
	\begin{subfigure}[b]{0.4\textwidth}
		\centering
		\begin{tikzpicture}[scale=1.6]

			\draw[green] (0,0) circle (0.75);

			\draw (0,0.3) to[out=175,in=90] (-0.4,-0.2) to[out=-90,in=180] (0,-0.6) to[out=0,in=-90] (0.9,0) to[out=90,in=-40] (0,1.2);
			\draw (0,0.3) to[out=160,in=90] (-0.7,-0.2) to[out=-90,in=180] (0,-0.9) to[out=0,in=-90] (1.05,0) to[out=90,in=-30] (0,1.2);
			\draw(0,0.3) -- (0,1.2);
			\draw (0,0.3) to[out=10,in=90] (0.4,-0.2) to[out=-90,in=0] (0,-0.55) to[out=180,in=-90] (-0.75,0) to[out=90,in=-140] (0,1.2);
			\draw (0,0.3) to[out=20,in=90] (0.6,-0.2) to[out=-90,in=0] (0,-0.95) to[out=180,in=-90] (-1.05,0) to[out=90,in=-155] (0,1.2);
			\draw (0,0.3) to[out=170,in=90] (-0.5,-0.2) to[out=-90,in=180] (0,-0.65) arc(90:-90:0.1) (0,-0.85) to[out=180,in=-90] (-0.95,0) to[out=90,in=-150] (0,1.2);

			\draw[red] (0,0.3) to[out=15,in=90] (0.5,-0.2) to[out=-90,in=20] (0,-0.75)
			(0,0.3) to[out=165,in=90] (-0.6,-0.2) to[out=-90,in=160] (0,-0.75);
			\draw[red] (0,1.2) to[out=-35,in=90] (0.95,0) to[out=-90,in=0] (0,-0.75)
			(0,1.2) to[out=-145,in=90] (-0.85,0) to[out=-90,in=180] (0,-0.75);

			\draw (0,0) circle (1.2)
			circle (.3);

			\draw[fill] (0,0.3) circle (0.06)
			(0,-0.75) circle (0.06)
			(0,1.2) circle (0.06);

		\end{tikzpicture}
		\caption{}
		\label{fig:J11 V6}
	\end{subfigure}
	\begin{subfigure}[b]{0.4\textwidth}
		\centering
		\begin{tikzpicture}[scale=1.6]

			\draw[green] (0,1.2) to[out=-45,in=90] (0.75,0) to[out=-90,in=0] (0,-0.65) to[out=180,in=-90] (-0.75,0) to[out=90,in=-135] (0,1.2);

			\draw (0,0.3) to[out=175,in=90] (-0.4,-0.2) to[out=-90,in=180] (0,-0.6) to[out=0,in=-90] (0.85,0) to[out=90,in=-40] (0,1.2);
			\draw (0,0.3) to[out=160,in=90] (-0.7,-0.2) to[out=-90,in=180] (0,-0.9) to[out=0,in=-90] (1.05,0) to[out=90,in=-30] (0,1.2);
			\draw(0,0.3) -- (0,1.2);
			\draw (0,0.3) to[out=20,in=90] (0.6,-0.2) to[out=-90,in=0] (0,-0.95) to[out=180,in=-90] (-1.05,0) to[out=90,in=-150] (0,1.2);
			\draw (0,0.3) to[out=170,in=90] (-0.5,-0.2) to[out=-90,in=180] (0,-0.65) arc(90:-90:0.1) (0,-0.85) to[out=180,in=-90] (-0.95,0) to[out=90,in=-145] (0,1.2);

			\draw[red] (0,0.3) to[out=15,in=90] (0.5,-0.2) to[out=-90,in=20] (0,-0.75)
			(0,0.3) to[out=165,in=90] (-0.6,-0.2) to[out=-90,in=160] (0,-0.75)
			(0,.3) to[out=175,in=90] (-0.32,0) to[out=-90,in=180] (0,-0.35) to[out=0,in=-90] (0.4,0) to[out=90,in=0] (0,0.4) to[out=180,in=90] (-0.65,-0.2) to[out=-90,in=170] (0,-0.75);
			\draw[red] (0,1.2) to[out=-35,in=90] (0.95,0) to[out=-90,in=0] (0,-0.75)
			(0,1.2) to[out=-140,in=90] (-0.85,0) to[out=-90,in=180] (0,-0.75);

			\draw (0,0) circle (1.2)
			circle (.3);

			\draw[fill] (0,0.3) circle (0.06)
			(0,-0.75) circle (0.06)
			(0,1.2) circle (0.06);

		\end{tikzpicture}
		\caption{}
		\label{fig:J11 V5}
	\end{subfigure}
	\caption{}
	\label{fig:J11 finals}
\end{figure}

\textbf{Step 2.} In this step we assume that \( |V| = 5 \). 

In this step we must have \( |W| \ge 5 \). Without loss of generality, we may assume that \( |W_+| \ge 3 \). We note that if \( v_{-1j}, v_{1k} \in V \), for any \( i,j \), it must be that \( |W_+| \le 2 \), using \Cref{eq:vjk wl ints}. Therefore, it must be that
\[ V = \left\{ v_{00}, v_{0,1}, v_{0,-1}, v_{1,0}, v_{1,-1} \right\}. \]
This lets us determine that
\[ W = W_+ \cup W_- = \left\{ w_{-1/2}, w_{1/2}, w_{3/2} \right\} \cup \left\{ h(w_{-1/2}), h(w_{1/2}) \right\}. \]

We can see that every arc in \( V \cup W \) intersects at least one other arc in \( V \cup W \), except \( w_{1/2} \). The arc \( y \) does not intersect \( w_{1/2} \). Any arc of the form \( x_k \) intersects \( w_{1/2} \) either zero times or at least twice. It must be that \( X \cup Y = h(x_k) \) for some \( k \), and indeed since \( h(w_{-1/2}), h(w_{1/2}) \in W \), it must be that \( X \cup Y = h(x_0) \). We obtain a 1-system as shown in \Cref{fig:J11 V5}.

\textbf{Step 3.} In this step we assume that \( |V| = 4 \). 

In this step we must have \( |W| \ge 6 \), so by the above it must be that \( |W_+| = |W_-| = 3 \). Similar to the above, we also note that if \( v_{i,-1}, v_{j,1} \in V \) for any \( i,j \), then we would have \( W_- \subset \left\{ h(w_{-1/2}), h(w_{1/2}) \right\} \), contradicting the fact that \( |W_-| = 3 \). From this we may conclude that 
\[ V = \left\{ v_{00}, v_{0,-1}, v_{1,0}, v_{1,-1} \right\}. \]
This lets us determine that
\[ W = W_+ \cup W_- = \left\{ w_{-1/2}, w_{1/2}, w_{3/2} \right\} \cup \left\{ h(w_{-1/2}), h(w_{1/2}), h(w_{3/2}) \right\}. \]

We can see that every arc in \( V \cup W \) intersects at least one other arc in \( V \cup W \), except \( w_{1/2} \) and \( h(w_{1/2}) \). The arc \( y \) does not intersect either \( w_{1/2} \) or \( h(w_{1/2}) \). Any arc of the form \( x_k \) intersects \( w_{1/2} \) either 0 times or at least twice, and similarly any arc of the form \( h(x_k) \) intersects \( h(w_{1/2}) \) either zero times or at least twice. There is no arc which intersects \( w_{1/2} \) and \( h(w_{1/2}) \) exactly once each. So, in this step, we do not recover a 1-system.



\subsection{$u$ is a non-loop arc}

In this case, \( S \setminus J \) is a torus with an open disk removed, and two marked points on the single boundary component.

Let \( S' \) be the surface obtained from \( S \) by quotienting the boundary to a single point, and let \( q: S \to S' \) be the quotient map. Let \( P' = q(\partial (S \setminus J)) = \left\{ o \right\} \) be the single point in the image of the the boundary. So, \( (S', P') \) is a torus with one marked point.

Since all arcs on \( (S', P') \) are loops, we have a natural map from arcs on \( (S', P') \) to simple closed curves on \( S' \). Since homotopic arcs are taken under this map to homotopic curves, this descends to a well-defined map of homotopy classes. Every curve can be homotoped to pass through \( o \), so this map is surjective. 

To see that the map is injective, suppose two arcs \( u, v \) are both mapped to the same homotopy class of curves. We may assume that \( u, v \) are in minimal position. However, their images intersect as curves at \( o \). Homotopic curves that intersect form a bigon. However, as arcs, \( u, v \) do not form any bigons or half bigons, since they are in minimal position. This can only occur if the marked point appears in the bigon twice, as shown in \Cref{fig:J12 intro a}. Clearly, this bigon shows that \( u, v \) are homotopic as arcs. Therefore, homotopy classes of arcs on \( (S', P') \) are in bijection with homotopy classes of simple closed curves on \( S' \). 

By the discussion on page 19 of~\cite{farb_2012}, homotopy classes of arcs on \( (S', P') \) are in bijection with \( \mathbb{Q} \cup \left\{ \infty \right\} \), which represents the slope. We also get a formula for the intersection number of a pair of closed curves. Since intersection number between arcs only counts intersections outside of \( P' \), the intersection between two arcs is given by 
\begin{equation}
	i(a/b, c/d) = |ad - bc| - 1.
	\label{eq:q torus ints}
\end{equation}

For any arc \( u \) on \( S \setminus J \), its image under \( q \) is a possibly non-essential arc \( q \circ u \) on \( (S', P') \). We write \( \mathcal{A} \setminus J = X \sqcup Y \), where \( Y \) contains the arcs \( y \) such that \( q \circ y \) is non-essential. We write \( X = \bigsqcup _{\lambda \in \mathbb{Q} \cup \left\{ \infty \right\}} X_\lambda \), where \( X_\lambda \) contains the arcs \( x \) such that \( q \circ x \) is homotopic to the arc with slope number \( \lambda \). 

First, consider \( Y \). There are exactly two arcs on \( S \setminus J \) which map to non-essential arcs on \( S' \), as shown in \Cref{fig:J12 intro b}. We see that these arcs intersect twice, so \( Y = \left\{ y \right\}, Y = \left\{ y' \right\} \), or \( Y = \emptyset \). In particular we have \( |Y| \le 1 \). 

Next, we consider \( X \). For any pair of arcs \( u, v \in X \), we note that \( i(u, v) \ge i(q \circ u, q \circ v) \), since the images might not be in minimal position. This implies that \( \left\{ \lambda | X_\lambda \ne \emptyset \right\} \) is a 1-system on \( (S', P') \). So we first classify 1-systems on \( (S', P') \). 

Let \( B \) be a 1-system on \( (S', P') \). First suppose we have \( u, v \in B \) which intersect. Without loss of generality we may assume that \( u, v \) are represented by \( \frac{1}{1} \) and \( -\frac{1}{1} \), respectively. Applying \Cref{eq:q torus ints}, we see that \( B \subset \left\{ \frac{1}{1}, -\frac{1}{1}, \frac{1}{0}, \frac{0}{1} \right\} \).

The other possibility is that the arcs in \( B \) are pairwise disjoint. Without loss of generality, we may assume \( u, v \in B \) are represented by \( \frac{1}{0}, \frac{0}{1} \) respectively. Applying \Cref{eq:q torus ints}, we see that again \( B \subset \left\{ \frac{1}{1}, -\frac{1}{1}, \frac{1}{0}, \frac{0}{1} \right\} \).

\begin{figure}
	\centering
	\begin{subfigure}{0.3\textwidth}
		\centering
		\begin{tikzpicture}
			\draw (0,0) circle (1.3)
			(0.4,0) circle (1.7);

			\fill[color=white!70!red] (-1.3,0) arc (180:-180:1.3) arc (-180:180:1.7);

			\draw[fill] (-1.3,0) circle (0.05);

			\node at (-1.5,0) {\( o \)};
			\node at (1.1,0) {\( u \)};
			\node at (2.3,0) {\( v \)};
		\end{tikzpicture}
		\caption{}
		\label{fig:J12 intro a}
	\end{subfigure}
	\begin{subfigure}{0.3\textwidth}
		\centering
		\begin{tikzpicture}[scale=1.6]
			\draw[dashed] (0,0) -- (0,2) -- (2,2) -- (2,0) -- (0,0);
			\draw[red] (1.7,1.4) to[out=130,in=20] (1,1.6) to[out=-160,in=160] (1,1.2) to[out=-20,in=-130] (1.7,1.4);

			\draw (1.1,1.4) to[out=35,in=145] (1.7,1.4) to[out=-145,in=-35] (1.1,1.4);
			\draw[fill] (1.1,1.4) circle (0.06)
			(1.7,1.4) circle (0.06);

			\draw (1.1,1.4) to[out=50,in=160] (1.8,1.6) to[out=-20,in=20] (1.8,1.2) to[out=-160,in=-50] (1.1,1.4);

			\node[red] at (1,1) {\( y \)};
			\node at (1.7,1) {\( y' \)};
		\end{tikzpicture}
		\caption{}
		\label{fig:J12 intro b}
	\end{subfigure}
	\caption{}
	\label{fig:J12 intro}
\end{figure}

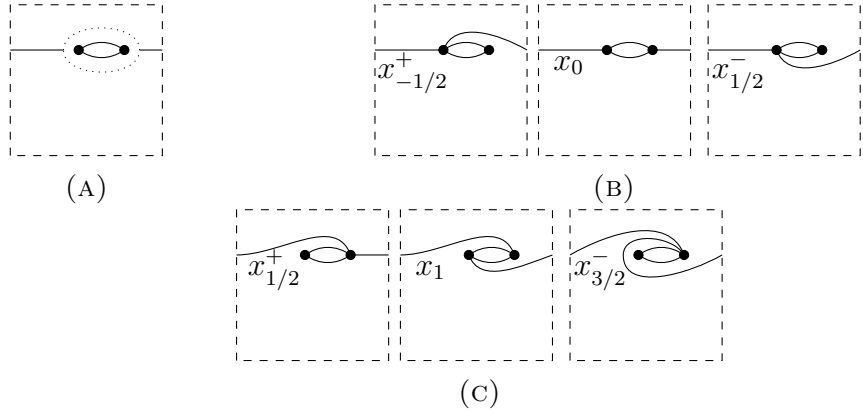
\begin{figure}
	\centering
	\begin{subfigure}{0.2\textwidth}
		\centering
		\begin{tikzpicture}
			\draw[dashed] (0,0) -- (0,2) -- (2,2) -- (2,0) -- (0,0);

			\draw (0.9,1.4) to[out=35,in=145] (1.5,1.4) to[out=-145,in=-35] (0.9,1.4);
			\draw[fill] (0.9,1.4) circle (0.06)
			(1.5,1.4) circle (0.06);

			\draw (0.7,1.4) -- (0,1.4) 
			(2,1.4) -- (1.7,1.4);

			\draw[dotted] (0.7,1.4) to[out=90,in=90] (1.7,1.4) to[out=-90,in=-90] (0.7,1.4);
		\end{tikzpicture}
		\caption{}
		\label{fig:J1 nl x0s a}
	\end{subfigure}
	\begin{subfigure}{0.6\textwidth}
		\centering
		\begin{tikzpicture}
			\node at (0.5,1.1) {\( x_{-1/2}^+ \)};
			\draw[dashed] (0,0) -- (0,2) -- (2,2) -- (2,0) -- (0,0);

			\draw (0.9,1.4) to[out=35,in=145] (1.5,1.4) to[out=-145,in=-35] (0.9,1.4);
			\draw[fill] (0.9,1.4) circle (0.06)
			(1.5,1.4) circle (0.06);

			\draw (0.9,1.4) -- (0,1.4) 
			(2,1.4) to[out=150,in=70] (0.9,1.4);
		\end{tikzpicture}
		\begin{tikzpicture}
			\node at (0.4,1.2) {\( x_0 \)};

			\draw[dashed] (0,0) -- (0,2) -- (2,2) -- (2,0) -- (0,0);

			\draw (0.9,1.4) to[out=35,in=145] (1.5,1.4) to[out=-145,in=-35] (0.9,1.4);
			\draw[fill] (0.9,1.4) circle (0.06)
			(1.5,1.4) circle (0.06);

			\draw (0.9,1.4) -- (0,1.4) 
			(2,1.4) -- (1.5,1.4);
		\end{tikzpicture}
		\begin{tikzpicture}
			\node at (0.4,1.1) {\( x_{1/2} ^- \)};
			\draw[dashed] (0,0) -- (0,2) -- (2,2) -- (2,0) -- (0,0);

			\draw (0.9,1.4) to[out=35,in=145] (1.5,1.4) to[out=-145,in=-35] (0.9,1.4);
			\draw[fill] (0.9,1.4) circle (0.06)
			(1.5,1.4) circle (0.06);

			\draw (0.9,1.4) -- (0,1.4) 
			(2,1.4) to[out=-150,in=-70] (0.9,1.4);
		\end{tikzpicture}
		\caption{}
		\label{fig:J1 nl x0s b}
	\end{subfigure}
	\begin{subfigure}{0.6\textwidth}
		\centering
		\begin{tikzpicture}
			\node at (0.5,1.2) {\( x_{1/2}^+ \)};
			\draw[dashed] (0,0) -- (0,2) -- (2,2) -- (2,0) -- (0,0);

			\draw (0.9,1.4) to[out=35,in=145] (1.5,1.4) to[out=-145,in=-35] (0.9,1.4);
			\draw[fill] (0.9,1.4) circle (0.06)
			(1.5,1.4) circle (0.06);

			\draw (0,1.4) to[out=0,in=110] (1.5,1.4)  
			(2,1.4) -- (1.5,1.4);
		\end{tikzpicture}
		\begin{tikzpicture}
			\node at (0.4,1.2) {\( x_1 \)};

			\draw[dashed] (0,0) -- (0,2) -- (2,2) -- (2,0) -- (0,0);

			\draw (0.9,1.4) to[out=35,in=145] (1.5,1.4) to[out=-145,in=-35] (0.9,1.4);
			\draw[fill] (0.9,1.4) circle (0.06)
			(1.5,1.4) circle (0.06);

			\draw (1.5,1.4) to[out=110,in=0] (0,1.4) 
			(2,1.4) to[out=-160,in=-70] (0.9,1.4);
		\end{tikzpicture}
		\begin{tikzpicture}
			\node at (0.4,1.2) {\( x_{3/2} ^- \)};
			\draw[dashed] (0,0) -- (0,2) -- (2,2) -- (2,0) -- (0,0);

			\draw (0.9,1.4) to[out=35,in=145] (1.5,1.4) to[out=-145,in=-35] (0.9,1.4);
			\draw[fill] (0.9,1.4) circle (0.06)
			(1.5,1.4) circle (0.06);

			\draw (1.5,1.4) to[out=110,in=30] (0,1.4) 
			(2,1.4) to[out=-150,in=-90] (0.7,1.4) to[out=90,in=120] (1.5,1.4);
		\end{tikzpicture}
		\caption{}
		\label{fig:J1 nl x0s c}
	\end{subfigure}
	\caption{Possible elements of \( X \) with slope 0.}
	\label{fig:J1 nl x0s}
\end{figure}

Suppose we have an arc \( x \in X_0 \). This determines \( x \) outside a neighborhood of the boundary component. This is shown in \Cref{fig:J1 nl x0s a}. Let \( h: S \setminus J \to S \setminus J \) be the homeomorphism that rotates the boundary component halfway, like a half Dehn twist about the boundary. Up to applying \( h \), we may fix one end of the arc \( x \). Then, there are three possibilities for the other end. See \Cref{fig:J1 nl x0s b}. 

We see that the arcs in \Cref{fig:J1 nl x0s b} are the only possible arcs in \( X_0 \), up to the action of \( h \). We label them as follows: let \( x_0, x_{-1/2}^+, x_{1/2}^- \) be the arcs shown in \Cref{fig:J1 nl x0s b}. For each \( n \in \mathbb{Z} \) denote by \( x_n \) the arc \( h^n(x_0) \). Similarly, define \( x_{n+1/2}^\pm := h^n(x_{1/2}^\pm) \). For example, see \Cref{fig:J1 nl x0s c}. 

With these definitions, we find the following: \( x_n^{*}, x_m ^{\bullet} \) intersect at most once if and only \( |n - m| \le 1 \) and \( (*, \bullet) \ne (+,+), (-,-) \). 

With this characterization, it is clear that the indices of arcs in \( X_0 \) must lie between \( n \) and \( n+1 \) for some \( n \in \frac{1}{2} \mathbb{Z} \). If \( n \in \mathbb{Z} \), without loss of generality we may assume that \( n=0 \), and we see that \( X_0 \subset \left\{ x_0, x_{1/2}^+, x_{1/2}^-, x_1 \right\} =: A_0 \). If instead \( n \in \mathbb{Z} + \frac{1}{2} \) we may assume that \( n = -\frac{1}{2} \) and without loss of generality we see that \( X_0 \subset \left\{ x_{-1/2}^+, x_0, x_{1/2}^{-} \right\} =: B_0 \). In both cases, we have \( |X_0| \le 4 \).





Define \( A_1 \) to be the set of arcs shown in black in \Cref{fig:J12 3 no 2 a} and define \( B_1 \) to be the set of arcs shown in black in \Cref{fig:J12 3 no 2 b}. Note that \( A_1, B_1 \) are obtained from \( A_0, B_0 \), respectively, by a single Dehn twist.

\begin{claim}
	If \( |X_{\pm 1}| \ge 3 \) then \( |X_{\mp 1}| \le 1 \). 
	\label{claim:3 no 2}
\end{claim}


\begin{proof}
	Suppose \( X_1 \subset B_1, |X_1| \ge 3 \). Consider the arc \( u \), shown in red in \Cref{fig:J12 3 no 2 a}, which intersects every arc of \( A_1 \) exactly once. Note that any other arc with slope number \( -1 \) would intersect at least two arcs of \( A_1 \), and therefore at least one arc of \( X_1 \), at least twice. Therefore, \( X_{-1} \subset \left\{ v \right\} \).  

	Suppose \( X_1 = A_1 \). Consider the arc \( v \), shown in red in \Cref{fig:J12 3 no 2 a}, which intersects every arc in \( X_1 \) exactly once. Note that any other arc with slope number \( -1 \) intersects at least one arc of \( A_1 \) at least twice. Therefore, \( X_{-1} \subset \left\{ v \right\} \). 
\end{proof}


\begin{figure}
	\centering
	\begin{subfigure}[b]{0.4\textwidth}
		\centering
		\begin{tikzpicture}[scale=1.6]
			\draw[dashed] (0,0) rectangle (2,2);

			\draw[red] (1.1,1.4) -- (0.7,2) (0.7,0) -- (0,0.7) (2,0.7) -- (1.7,1.4);

			\draw (1.1,1.4) to[out=35,in=145] (1.7,1.4) to[out=-145,in=-35] (1.1,1.4);
			\draw[fill] (1.1,1.4) circle (0.06)
			(1.7,1.4) circle (0.06);

			\draw (1.1,1.4) -- (0.2,0) (0.2,2) -- (0,1.8) 
			(2,1.8) -- (1.6,1.8) to[out=180,in=80] (1.1,1.4); 
			\draw (1.1,1.4) -- (0.3,0) (0.3,2) -- (0,1.6) (2,1.6) -- (1.7,1.4); 
			\draw (1.7,1.4) -- (0.4,0) 
			(0.4,2) -- (0,1.7) 
			(2,1.7) -- (1.6,1.7) to[out=180,in=60] (1.1,1.4); 
			\draw (1.7,1.4) -- (0.5,0) (0.5,2) -- (0,1.5) (2,1.5) -- (1.7,1.4); 
		\end{tikzpicture}
		\caption{}
		\label{fig:J12 3 no 2 a}
	\end{subfigure}
	\begin{subfigure}[b]{0.4\textwidth}
		\centering
		\begin{tikzpicture}[scale=1.6]
			\draw[dashed] (0,0) rectangle (2,2);

			\draw[red] (1.1,1.4) to[out=-90,in=160] (2,0.6) (0,0.6) -- (0.7,0) (0.7,2) -- (1.1,1.4);

			\draw (1.1,1.4) to[out=35,in=145] (1.7,1.4) to[out=-145,in=-35] (1.1,1.4);
			\draw[fill] (1.1,1.4) circle (0.06)
			(1.7,1.4) circle (0.06);

			\draw (1.1,1.4) -- (0.2,0) (0.2,2) -- (0,1.8) 
			(2,1.8) -- (1.5,1.8) to[out=180,in=80] (1.1,1.4);
			\draw (1.1,1.4) -- (0.3,0) (0.3,2) -- (0,1.7) (2,1.7) -- (1.7,1.4);
			\draw (1.1,1.4) -- (0.4,0) (0.4,2) -- (0,1.6) (2,1.6) to[out=-135,in=0] (1.6,1.2) -- (1.3,1.2) to[out=180,in=-70] (1.1,1.4);
		\end{tikzpicture}
		\caption{}
		\label{fig:J12 3 no 2 b}
	\end{subfigure}
	\caption{}
	\label{fig:J12 3 no 2}
\end{figure}


\begin{claim}
	If \( |X_1| = 4 \) then \( |X_0| \le 3 \). 
	\label{claim:4 no 4}
\end{claim}

\begin{proof}
	Up to homeomorphism, we may assume that \( X_1 = A_1 \). Then, there are five arcs with slope number 1 which intersect no arc in \( X_1 \) more than once. These are shown in red in \Cref{fig:J12 4 no 4}. We see that the pair \( a,c \) intersect twice, as well as the pair \( b,d \). Therefore, \( \mathcal{A} \) contains at most 3 of the five possible arcs. So, we obtain \( |X_0| \le 3 \) as desired. 
\end{proof}

\begin{cor}
	If \( X_1 = B_1 \) and \( |X_0| = 3 \), then \( e \in X_0 \). Furthermore, \( X_0 \) must be one of \( \left\{ a,b,e \right\}, \left\{ a, d, e \right\}, \left\{ c, b, e \right\}, \left\{ c, d, e \right\} \). 
	\label{cor:4 which 3}
\end{cor}

\begin{figure}
	\centering
		\begin{tikzpicture}[scale=1.8]
			\draw[dashed] (0,0) rectangle (2,2);

			\draw[red] (1.1,1.4) -- (0,1.3) (2,1.3) -- (1.7,1.4) 
			(1.7,1.4) to[out=-135,in=0] (0,1.15) (2,1.15) -- (1.7,1.4) 
			(1.1,1.4) -- (0,1.2) (2,1.2) to[out=-170,in=-40] (1.1,1.4) 
			(1.7,1.4) to[out=135,in=10] (0,1.4) (2,1.4) -- (1.7,1.4) 
			(1.1,1.4) -- (0,1.35) (2,1.35) to[out=140,in=50] (1.1,1.4); 

			\draw (1.1,1.4) to[out=35,in=145] (1.7,1.4) to[out=-145,in=-35] (1.1,1.4);
			\draw[fill] (1.1,1.4) circle (0.06)
			(1.7,1.4) circle (0.06);

			\draw (1.1,1.4) -- (0.2,0) (0.2,2) -- (0,1.8) 
			(2,1.8) -- (1.6,1.8) to[out=180,in=80] (1.1,1.4); 
			\draw (1.1,1.4) -- (0.3,0) (0.3,2) -- (0,1.6) (2,1.6) -- (1.7,1.4); 
			\draw (1.7,1.4) -- (0.4,0) 
			(0.4,2) -- (0,1.7) 
			(2,1.7) -- (1.6,1.7) to[out=180,in=60] (1.1,1.4); 
			\draw (1.7,1.4) -- (0.5,0) (0.5,2) -- (0,1.5) (2,1.5) -- (1.7,1.4); 

			\node[red] at (1,1.75) {\( a \)};
			\node[red] at (0.6,1) {\( b \)};
			\node[red] at (1.73,1.6) {\( c \)};
			\node[red] at (1.7,1.05) {\( d \)};
			\node[red] at (-0.1,1.3) {\( e \)};
		\end{tikzpicture}
	\caption{}
	\label{fig:J12 4 no 4}
\end{figure}

Now, we show that if \( Y = \emptyset \), we have \( |X_t| = 4 \) for some \( t \). Suppose instead that \( |X_t| \le 3 \) for all \( t \). We have \( |X_1 \cup X_{-1}| \le 4 \), since if either has three arcs, the other has at most one by \Cref{claim:3 no 2}. Then \( |\mathcal{A} \setminus J| = |X_0| + |X_{\infty}| + |X_1| + |X_{-1}| \le 3 + 3 +4 = 10 \). Since \( |J| = 1 \), this contradicts the assumtion that \( \mathcal{A} \) is maximal. 

We split the proof into three cases. 

\textbf{Case 1.} \( Y = \emptyset, |X_1| = 4 \).

By \Cref{claim:3 no 2}, it must be that \( |X_{-1}| \le 1 \). By \Cref{claim:4 no 4}, it must be that \( |X_0| \le 3, |X_\infty| \le 3 \). We must have \( |X_1| + |X_0| + |X_\infty| + |X_{-1}| = 11 \), because we assume \( Y = \emptyset \), so in fact these inequalities are equalities. 

Up to homeomorphism, \( X_1 \cup X_{-1} \) must be the arcs shown in \Cref{fig:J12 3 no 2 b}. 
Since \( |X_0| = 3 \), it must be one of the sets described in \Cref{cor:4 which 3}.

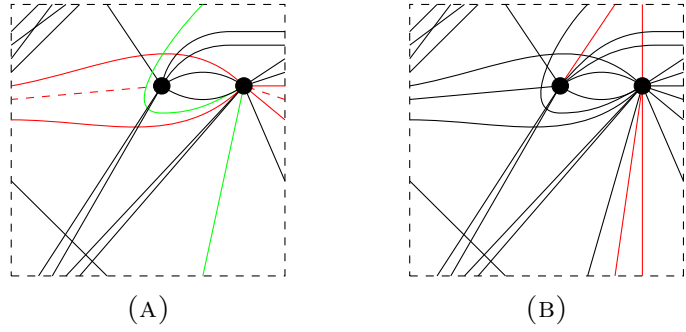
\begin{figure}
	\centering
	\begin{subfigure}[b]{0.3\textwidth}
		\centering
		\begin{tikzpicture}[scale=1.8]
			\draw[dashed] (0,0) rectangle (2,2);

			\draw[dashed,red] (1.1,1.4) -- (0,1.3) (2,1.3) -- (1.7,1.4); 
			\draw[red] (1.7,1.4) to[out=-135,in=0] (0,1.15) (2,1.15) -- (1.7,1.4) 
			(1.7,1.4) to[out=135,in=10] (0,1.4) (2,1.4) -- (1.7,1.4); 

			\draw (1.1,1.4) -- (0.7,2) (0.7,0) -- (0,0.7) (2,0.7) -- (1.7,1.4);

			\draw[green] (1.7,1.4) to[out=-150,in=0] (1.1,1.2) to[out=180,in=-135]  (1.4,2) (1.4,0) -- (1.7,1.4);

			\draw (1.1,1.4) to[out=35,in=145] (1.7,1.4) to[out=-145,in=-35] (1.1,1.4);
			\draw[fill] (1.1,1.4) circle (0.06)
			(1.7,1.4) circle (0.06);

			\draw (1.1,1.4) -- (0.2,0) (0.2,2) -- (0,1.8) 
			(2,1.8) -- (1.6,1.8) to[out=180,in=80] (1.1,1.4); 
			\draw (1.1,1.4) -- (0.3,0) (0.3,2) -- (0,1.6) (2,1.6) -- (1.7,1.4); 
			\draw (1.7,1.4) -- (0.4,0) 
			(0.4,2) -- (0,1.7) 
			(2,1.7) -- (1.6,1.7) to[out=180,in=60] (1.1,1.4); 
			\draw (1.7,1.4) -- (0.5,0) (0.5,2) -- (0,1.5) (2,1.5) -- (1.7,1.4); 

		\end{tikzpicture}
		\caption{}
		\label{fig:J12 final first a}
	\end{subfigure}
	\begin{subfigure}[b]{0.3\textwidth}
		\centering
		\begin{tikzpicture}[scale=1.8]
			\draw[dashed] (0,0) rectangle (2,2);

			\draw (1.1,1.4) -- (0,1.3) (2,1.3) -- (1.7,1.4); 
			\draw (1.7,1.4) to[out=-135,in=0] (0,1.15) (2,1.15) -- (1.7,1.4) 
			(1.7,1.4) to[out=135,in=10] (0,1.4) (2,1.4) -- (1.7,1.4); 

			\draw (1.1,1.4) -- (0.7,2) (0.7,0) -- (0,0.7) (2,0.7) -- (1.7,1.4);

			\draw (1.7,1.4) to[out=-150,in=0] (1.1,1.2) to[out=180,in=-135]  (1.3,2) (1.3,0) -- (1.7,1.4);
			\draw[red] (1.1,1.4) -- (1.5,2) (1.5,0) -- (1.7,1.4)
			(1.7,0) -- (1.7,2);

			\draw (1.1,1.4) to[out=35,in=145] (1.7,1.4) to[out=-145,in=-35] (1.1,1.4);
			\draw[fill] (1.1,1.4) circle (0.06)
			(1.7,1.4) circle (0.06);

			\draw (1.1,1.4) -- (0.2,0) (0.2,2) -- (0,1.8) (2,1.8) -- (1.8,1.8) to[out=180,in=50] (1.1,1.4); 
			\draw (1.1,1.4) -- (0.3,0) (0.3,2) -- (0,1.6) (2,1.6) -- (1.7,1.4); 
			\draw (1.7,1.4) -- (0.4,0) (0.4,2) -- (0,1.7) (2,1.7) -- (1.8,1.7) to[out=180,in=40] (1.1,1.4); 
			\draw (1.7,1.4) -- (0.5,0) (0.5,2) -- (0,1.5) (2,1.5) -- (1.7,1.4); 

		\end{tikzpicture}
		\caption{}
		\label{fig:J12 final first b}
	\end{subfigure}
	\caption{Constructing a 1-system}
	\label{fig:J12 final first}
\end{figure}

If \( X_0 = \left\{ a, b, e \right\} \) shown in red in \Cref{fig:J12 final first a} we note that the arc \( e \), shown in dashed red in \Cref{fig:J12 final first a}, intersects no other arc in \( X_0 \cup X_1 \cup X_{-1} \). This means that there must be some arc in \( X_\infty \) which intersects the arc \( e \). We see that the green arc in \Cref{fig:J12 final first a} is the only arc with slope number \( \infty \) which intersects \( e \), and intersects no arc of \( \mathcal{A} \) twice. So that arc must be included in \( \mathcal{A} \).
Then, there is a unique set of two arcs with slope number \( \infty \) which do not intersect any of the other arcs proved to be in \( \mathcal{A} \) more than once. These are shown in red in \Cref{fig:J12 final first b}. So, we obtain a 1-system.



\begin{figure}
	\centering
	\begin{subfigure}[b]{0.3\textwidth}
		\centering
		\begin{tikzpicture}[scale=1.8]
			\draw[dashed] (0,0) rectangle (2,2);

			\draw[red] (1.1,1.4) -- (0,1.3) (2,1.3) -- (1.7,1.4); 
			\draw[red] (1.1,1.4) -- (0,1.2) (2,1.2) to[out=-170,in=-45] (1.1,1.4) 
			(1.7,1.4) to[out=135,in=10] (0,1.4) (2,1.4) -- (1.7,1.4); 

			\draw (1.1,1.4) -- (0.7,2) (0.7,0) -- (0,0.7) (2,0.7) -- (1.7,1.4);

			\draw (1.1,1.4) to[out=35,in=145] (1.7,1.4) to[out=-145,in=-35] (1.1,1.4);
			\draw[fill] (1.1,1.4) circle (0.06)
			(1.7,1.4) circle (0.06);

			\draw (1.1,1.4) -- (0.2,0) (0.2,2) -- (0,1.8) 
			(2,1.8) -- (1.6,1.8) to[out=180,in=80] (1.1,1.4); 
			\draw (1.1,1.4) -- (0.3,0) (0.3,2) -- (0,1.6) (2,1.6) -- (1.7,1.4); 
			\draw (1.7,1.4) -- (0.4,0) 
			(0.4,2) -- (0,1.7) 
			(2,1.7) -- (1.6,1.7) to[out=180,in=60] (1.1,1.4); 
			\draw (1.7,1.4) -- (0.5,0) (0.5,2) -- (0,1.5) (2,1.5) -- (1.7,1.4); 

		\end{tikzpicture}
		\caption{}
		\label{fig:J12 final first non a}
	\end{subfigure}
	\begin{subfigure}[b]{0.3\textwidth}
		\centering
		\begin{tikzpicture}[scale=1.8]
			\draw[dashed] (0,0) rectangle (2,2);

			\draw[green] (1.1,1.4) -- (0,1.3) (2,1.3) -- (1.7,1.4); 
			\draw (1.1,1.4) -- (0,1.2) (2,1.2) to[out=-170,in=-45] (1.1,1.4) 
			(1.7,1.4) to[out=135,in=10] (0,1.4) (2,1.4) -- (1.7,1.4); 

			\draw (1.1,1.4) -- (0.7,2) (0.7,0) -- (0,0.7) (2,0.7) -- (1.7,1.4);

			\draw[red] (1.1,1.4) -- (1.5,2) (1.5,0) -- (1.7,1.4)
			(1.1,0) -- (1.1,2)
			(1.7,0) -- (1.7,2);

			\draw (1.1,1.4) to[out=35,in=145] (1.7,1.4) to[out=-145,in=-35] (1.1,1.4);
			\draw[fill] (1.1,1.4) circle (0.06)
			(1.7,1.4) circle (0.06);

			\draw (1.1,1.4) -- (0.2,0) (0.2,2) -- (0,1.8) (2,1.8) -- (1.8,1.8) to[out=180,in=50] (1.1,1.4); 
			\draw (1.1,1.4) -- (0.3,0) (0.3,2) -- (0,1.6) (2,1.6) -- (1.7,1.4); 
			\draw (1.7,1.4) -- (0.4,0) (0.4,2) -- (0,1.7) (2,1.7) -- (1.8,1.7) to[out=180,in=40] (1.1,1.4); 
			\draw (1.7,1.4) -- (0.5,0) (0.5,2) -- (0,1.5) (2,1.5) -- (1.7,1.4); 

		\end{tikzpicture}
		\caption{}
		\label{fig:J12 final first non b}
	\end{subfigure}
	\begin{subfigure}{0.3\textwidth}
		\centering
		\begin{tikzpicture}[scale=1.8]
			\draw[dashed] (0,0) rectangle (2,2);

			\draw[green] (1.1,1.4) -- (1.5,2) (1.5,0) -- (1.7,1.4);
			\node[green] at (1.5,2.15) {\( u \)};

			\draw[red] (1.1,1.4) -- (0,1.3) (2,1.3) -- (1.7,1.4) 
			(1.7,1.4) to[out=-135,in=0] (0,1.2) (2,1.2) -- (1.7,1.4) 
			(1.1,1.4) -- (0,1.4) (2,1.4) to[out=140,in=40] (1.1,1.4); 

			\draw (1.1,1.4) -- (0.7,2) (0.7,0) -- (0,0.7) (2,0.7) -- (1.7,1.4);

			\draw (1.1,1.4) to[out=35,in=145] (1.7,1.4) to[out=-145,in=-35] (1.1,1.4);
			\draw[fill] (1.1,1.4) circle (0.06)
			(1.7,1.4) circle (0.06);

			\draw (1.1,1.4) -- (0.2,0) (0.2,2) -- (0,1.8) 
			(2,1.8) -- (1.8,1.8) to[out=180,in=50] (1.1,1.4); 
			\draw (1.1,1.4) -- (0.3,0) (0.3,2) -- (0,1.6) (2,1.6) -- (1.7,1.4); 
			\draw (1.7,1.4) -- (0.4,0) 
			(0.4,2) -- (0,1.7) 
			(2,1.7) -- (1.8,1.7) to[out=180,in=45] (1.1,1.4); 
			\draw (1.7,1.4) -- (0.5,0) (0.5,2) -- (0,1.5) (2,1.5) -- (1.7,1.4); 

		\end{tikzpicture}
		\caption{}
		\label{fig:J12 case 1 third}
	\end{subfigure}
	\caption{Obtaining contradictions}
	\label{fig:J12 final first non}
\end{figure}
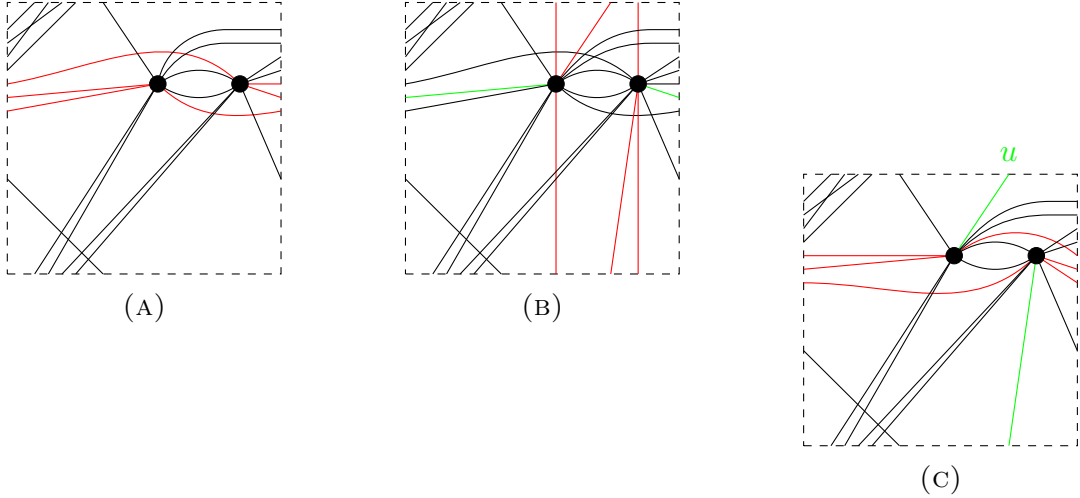

If \( X_0 = \left\{ a, d, e \right\} \), shown in red in \Cref{fig:J12 final first non a}, we again see that arc \( e \), shown in green in \Cref{fig:J12 final first non b}, intersects no other arc with slope number 0, 1, or \( -1 \). Up to a homeomorphism, we appeal to \Cref{claim:4 no 4} to see that there are five arcs with slope number \( \infty \) which do not intersect any arc in \( X_1 \) two or more times, and of these there are three which do not intersect any arc in \( X_0 \) two or more times. These are shown in red in \Cref{fig:J12 final first non b}. None of these arcs intersect the arc \( e \). So, this does not result in a 1-system without contradicting our assumption that \( |J| = 1 \). 

If \( X_0 = \left\{ b, c, e \right\} \), shown in red in \Cref{fig:J12 case 1 third}, consider the arc \( u \) shown in green in \Cref{fig:J12 case 1 third}. Up to applying a homeomorphism, we may appeal to 
\Cref{cor:4 which 3} and we conclude that \( u \) must be included in \( X_\infty \), since we know that \( |X_\infty| = 3 \). We also see that no arc in \( X_1 \cup X_0 \cup X_{-1} \) intersects \( u \), and no arc with slope number \( \infty \) which intersects \( u \) exactly once can be in \( X_\infty \) since \( \mathcal{A} \) is a 1-system. Therefore, we do not obtain a new 1-system.

If \( X_0 = \left\{ c, d, e \right\} \), we may apply a rotation of \( \pi \) and we see that this is equivalent to the case where \( X_0 = \left\{ a, b, e \right\} \). 

\textbf{Case 2.} \( Y = \emptyset, |X_\infty| = 4 \).

By \Cref{claim:4 no 4}, we know that \( |X_1| \le 3, |X_{-1}| \le 3 \). By \Cref{claim:3 no 2}, either \( |X_{\pm 1}| = 3 \) and \( |X_{\mp 1}| \le 1 \), or \( |X_1| = |X_{-1}| = 2 \). This means that \( | X_\infty \cup X_{1} \cup X_{-1}| \le 8 \), so it must be that \( |X_0| = 3 \).

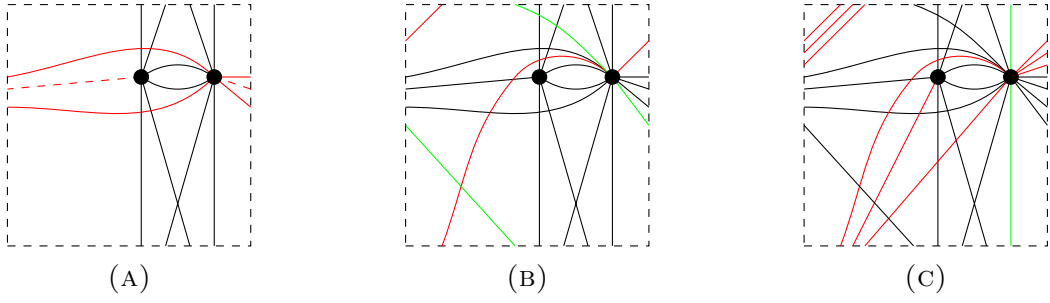
\begin{figure}
	\centering
	\begin{subfigure}[b]{0.3\textwidth}
		\centering
		\begin{tikzpicture}[scale=1.6]
			\draw[dashed] (0,0) rectangle (2,2);

			\draw (1.1,1.4) -- (1.1,2) (1.1,0) -- (1.1,1.4)
			(1.1,1.4) -- (1.3,2) (1.3,0) -- (1.7,1.4)
			(1.7,1.4) -- (1.5,2) (1.5,0) -- (1.1,1.4)
			(1.7,1.4) -- (1.7,2) (1.7,0) -- (1.7,1.4);

			\draw[dashed,red] (1.1,1.4) -- (0,1.3) (2,1.3) -- (1.7,1.4); 
			\draw[red] (1.7,1.4) to[out=-135,in=0] (0,1.15) (2,1.15) -- (1.7,1.4) 
			(1.7,1.4) to[out=135,in=10] (0,1.4) (2,1.4) -- (1.7,1.4); 

			\draw (1.1,1.4) to[out=35,in=145] (1.7,1.4) to[out=-145,in=-35] (1.1,1.4);
			\draw[fill] (1.1,1.4) circle (0.06)
			(1.7,1.4) circle (0.06);

		\end{tikzpicture}
		\caption{}
		\label{fig:J12 case 2 a}
	\end{subfigure}
	\begin{subfigure}[b]{0.3\textwidth}
		\centering
		\begin{tikzpicture}[scale=1.6]
			\draw[dashed] (0,0) rectangle (2,2);

			\draw (1.1,1.4) -- (1.1,2) (1.1,0) -- (1.1,1.4)
			(1.1,1.4) -- (1.3,2) (1.3,0) -- (1.7,1.4)
			(1.7,1.4) -- (1.5,2) (1.5,0) -- (1.1,1.4)
			(1.7,1.4) -- (1.7,2) (1.7,0) -- (1.7,1.4);

			\draw (1.1,1.4) -- (0,1.3) (2,1.3) -- (1.7,1.4) 
			(1.7,1.4) to[out=-135,in=0] (0,1.15) (2,1.15) -- (1.7,1.4) 
			(1.7,1.4) to[out=135,in=10] (0,1.4) (2,1.4) -- (1.7,1.4); 

			\draw[red] (1.7,1.4) -- (2,1.7) (0,1.7) -- (0.3,2) (0.3,0) to[out=70,in=-130] (0.9,1.4) to[out=50,in=135] (1.7,1.4);

			\draw[green] (1.7,1.4) -- (2,1) (0,1) -- (0.9,0) (0.9,2) to[out=-20,in=135] (1.7,1.4);

			\draw (1.1,1.4) to[out=35,in=145] (1.7,1.4) to[out=-145,in=-35] (1.1,1.4);
			\draw[fill] (1.1,1.4) circle (0.06)
			(1.7,1.4) circle (0.06);

		\end{tikzpicture}
		\caption{}
		\label{fig:J12 case 2 b}
	\end{subfigure}
	\begin{subfigure}[b]{0.3\textwidth}
		\centering
		\begin{tikzpicture}[scale=1.6]
			\draw[dashed] (0,0) rectangle (2,2);

			\draw (1.1,1.4) -- (1.1,2) (1.1,0) -- (1.1,1.4)
			(1.1,1.4) -- (1.3,2) (1.3,0) -- (1.7,1.4)
			(1.7,1.4) -- (1.5,2) (1.5,0) -- (1.1,1.4);
			\draw[green] (1.7,0) -- (1.7,2);

			\draw (1.1,1.4) -- (0,1.3) (2,1.3) -- (1.7,1.4) 
			(1.7,1.4) to[out=-135,in=0] (0,1.15) (2,1.15) -- (1.7,1.4) 
			(1.7,1.4) to[out=135,in=10] (0,1.4) (2,1.4) -- (1.7,1.4); 

			\draw[red] (1.7,1.4) -- (2,1.7) (0,1.7) -- (0.3,2) (0.3,0) to[out=70,in=-130] (0.9,1.4) to[out=50,in=135] (1.7,1.4)
			(1.1,1.4) -- (0.4,0) (0.4,2) -- (0,1.6) (2,1.6) -- (1.7,1.4)
			(1.7,1.4) -- (0.5,0) (0.5,2) -- (0,1.5) (2,1.5) -- (1.7,1.4);

			\draw (1.7,1.4) -- (2,1) (0,1) -- (0.9,0) (0.9,2) to[out=-20,in=135] (1.7,1.4);

			\draw (1.1,1.4) to[out=35,in=145] (1.7,1.4) to[out=-145,in=-35] (1.1,1.4);
			\draw[fill] (1.1,1.4) circle (0.06)
			(1.7,1.4) circle (0.06);

		\end{tikzpicture}
		\caption{}
		\label{fig:J12 case 2 c}
	\end{subfigure}
	\caption{Obtaining a contradiction}
	\label{fig:J12 case 2}
\end{figure}

Up to homeomorphism, \( X_\infty \) must be the set of arcs shown in black in \Cref{fig:J12 case 2 a}. Then, from \Cref{cor:4 which 3}, there are two choices for \( X_0 \).

If \( X_0 \) is the set of arcs shown in red in \Cref{fig:J12 case 2 a}, let \( e \) be the arc in dashed red. The arc \( e \) must intersect some other arc in \( \mathcal{A} \setminus J \) by assumption, so up to reflection in the horizontal axis \( e \) intersects some arc in \( X_1 \). The only arc which intersects \( e \), but intersects no arc of \( \mathcal{A} \) twice, is the arc shown in red in \Cref{fig:J12 case 2 b}. Then, we see that there is only one arc with slope number \( -1 \) which intersects no arc of \( \mathcal{A} \) twice, shown in green in \Cref{fig:J12 case 2 b}. This means \( |X_{-1}| \le 1 \), se we must have \( |X_{-1}| = 1, |X_1| = 3 \). There are exactly two more arcs with slope number one which intersect no arc in \( \mathcal{A} \setminus J \) at least twice, shown in red in \Cref{fig:J12 case 2 c}. However, in the resulting 1-system, there is an arc in \( \mathcal{A} \setminus J \), shown in green in \Cref{fig:J12 case 2 c}, which intersects no other arc. This is a contradiction.

\begin{figure}
	\centering
	\begin{subfigure}[b]{0.3\textwidth}
		\centering
		\begin{tikzpicture}[scale=1.6]
			\draw[dashed] (0,0) rectangle (2,2);

			\draw (1.1,1.4) -- (1.1,2) (1.1,0) -- (1.1,1.4)
			(1.1,1.4) -- (1.3,2) (1.3,0) -- (1.7,1.4)
			(1.7,1.4) -- (1.5,2) (1.5,0) -- (1.1,1.4)
			(1.7,1.4) -- (1.7,2) (1.7,0) -- (1.7,1.4);

			\draw[dashed,red] (1.1,1.4) -- (0,1.3) (2,1.3) -- (1.7,1.4); 
			\draw[red](1.1,1.4) -- (0,1.2) (2,1.2) to[out=-170,in=-40] (1.1,1.4) 
			(1.7,1.4) to[out=135,in=10] (0,1.4) (2,1.4) -- (1.7,1.4); 

			\draw (1.1,1.4) to[out=35,in=145] (1.7,1.4) to[out=-145,in=-35] (1.1,1.4);
			\draw[fill] (1.1,1.4) circle (0.06)
			(1.7,1.4) circle (0.06);

		\end{tikzpicture}
		\caption{}
		\label{fig:J12 case 2 other a}
	\end{subfigure}
	\begin{subfigure}[b]{0.3\textwidth}
		\centering
		\begin{tikzpicture}[scale=1.6]
			\draw[dashed] (0,0) rectangle (2,2);

			\draw (1.1,1.4) -- (1.1,2) (1.1,0) -- (1.1,1.4)
			(1.1,1.4) -- (1.3,2) (1.3,0) -- (1.7,1.4)
			(1.7,1.4) -- (1.5,2) (1.5,0) -- (1.1,1.4)
			(1.7,1.4) -- (1.7,2) (1.7,0) -- (1.7,1.4);

			\draw (1.1,1.4) -- (0,1.3) (2,1.3) -- (1.7,1.4) 
			(1.1,1.4) -- (0,1.2) (2,1.2) to[out=-170,in=-40] (1.1,1.4) 
			(1.7,1.4) to[out=135,in=10] (0,1.4) (2,1.4) -- (1.7,1.4); 

			\draw[red] (1.1,1.4) -- (0.4,0) (0.4,2) -- (0,1.6) (2,1.6) to[out=-120,in=-45] (1.1,1.4);

			\draw[green] (1.1,1.4) -- (2,0.7) (0,0.7) -- (0.7,0) (0.7,2) -- (1.1,1.4);

			\draw (1.1,1.4) to[out=35,in=145] (1.7,1.4) to[out=-145,in=-35] (1.1,1.4);
			\draw[fill] (1.1,1.4) circle (0.06)
			(1.7,1.4) circle (0.06);

		\end{tikzpicture}
		\caption{}
		\label{fig:J12 case 2 other b}
	\end{subfigure}
	\begin{subfigure}[b]{0.3\textwidth}
		\centering
		\begin{tikzpicture}[scale=1.6]
			\draw[dashed] (0,0) rectangle (2,2);

			\draw (1.1,1.4) -- (1.1,2) (1.1,0) -- (1.1,1.4)
			(1.1,1.4) -- (1.3,2) (1.3,0) -- (1.7,1.4)
			(1.7,1.4) -- (1.5,2) (1.5,0) -- (1.1,1.4)
			(1.7,1.4) -- (1.7,2) (1.7,0) -- (1.7,1.4);

			\draw (1.1,1.4) -- (0,1.3) (2,1.3) -- (1.7,1.4) 
			(1.1,1.4) -- (0,1.2) (2,1.2) to[out=-170,in=-40] (1.1,1.4) 
			(1.7,1.4) to[out=135,in=10] (0,1.4) (2,1.4) -- (1.7,1.4); 

			\draw[red] (1.1,1.4) -- (0.4,0) (0.4,2) -- (0,1.6) (2,1.6) to[out=-120,in=-45] (1.1,1.4)
			(1.1,1.4) -- (0.35,0) (0.35,2) -- (0,1.65) (2,1.65) -- (1.7,1.4)
			(1.1,1.4) -- (0.3,0) (0.3,2) -- (0,1.7) (2,1.7) to[out=-180,in=45] (1.1,1.4);

			\draw (1.1,1.4) -- (2,0.7) (0,0.7) -- (0.7,0) (0.7,2) -- (1.1,1.4);

			\draw (1.1,1.4) to[out=35,in=145] (1.7,1.4) to[out=-145,in=-35] (1.1,1.4);
			\draw[fill] (1.1,1.4) circle (0.06)
			(1.7,1.4) circle (0.06);

		\end{tikzpicture}
		\caption{}
		\label{fig:J12 case 2 other c}
	\end{subfigure}
	\caption{Constructing a 1-system}
	\label{fig:J12 case 2 other}
\end{figure}
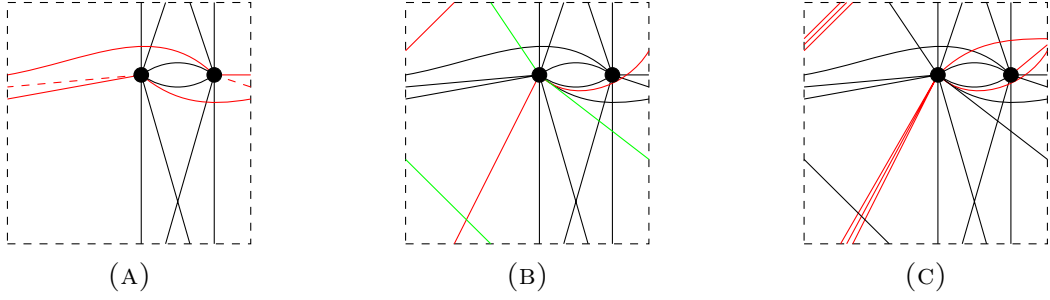

If \( X_0 \) is the set of arcs shown in \Cref{fig:J12 case 2 other a}, let \( e \) again be the arc in dashed red. Note that there is no arc with slope number \( -1 \) which intersects \( e \) but intersects no arc in \( \mathcal{A} \) twice. Therefore, there must be an arc in \( X_{1} \) which intersects \( e \). There are two possibilities, equivalent up to rotation of 180 degrees. One choice is shown in red in \Cref{fig:J12 case 2 other b}. Then, there is only one arc \( u \), shown in green in \Cref{fig:J12 case 2 other b}, with slope number \( -1 \) which intersects no arc in \( \mathcal{A} \) at least twice, so \( X_{-1} \subset \left\{ u \right\} \). Then it must be that \( |X_1| = 3, |X_{-1}| = 1 \). In fact, there is a unique choice for the other two arcs of \( X_1 \). Thus we obtain a 1-system.

\textbf{Case 3.} \( |Y| = 1 \). Without loss of generality, we may assume that \( Y = \left\{ y \right\} \), where \( y \) is the arc shown in red in \Cref{fig:J12 intro a}.

Note that each of \( y, y' \) intersects on arc of \( A_0 \) twice. This means that \( |X_0| \le 3 \). In fact, since \( y \) is invariant under Dehn twists and under \( h \), this means \( |X_\lambda| \le 3 \) for all \( \lambda \). Then, using \Cref{claim:3 no 2} we deduce that \( |X_0| = 3, |X_\infty| = 3 \) and either \( |X_1| = |X_{-1}| = 2 \) or \( |X_{\pm 1}| = 3, |X_{\mp 1}| = 1 \).

We break the proof into steps based on \( X_0 \). 

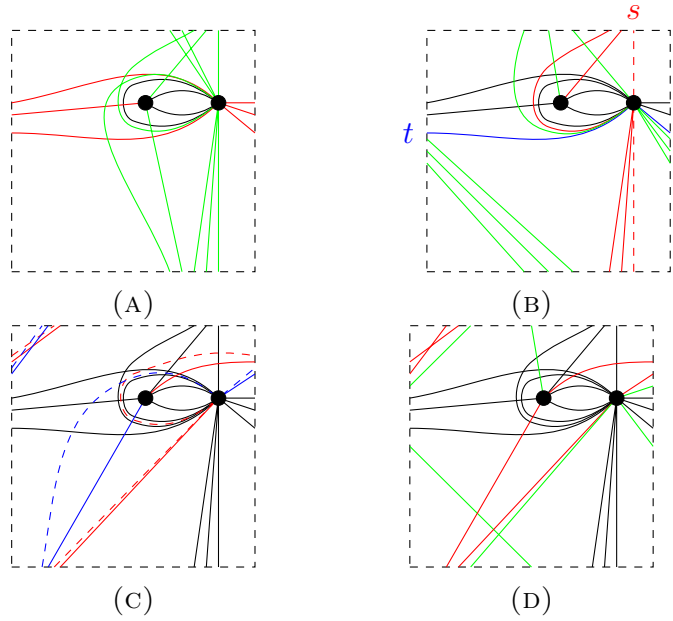
\begin{figure}
	\centering
	\begin{subfigure}[b]{0.3\textwidth}
		\centering
		\begin{tikzpicture}[scale=1.6]
			\draw[dashed] (0,0) rectangle (2,2);

			\draw (1.7,1.4) to[out=140,in=20] (1,1.55) to[out=-160,in=160] (1,1.25) to[out=-20,in=-140] (1.7,1.4);

			\draw[red] (1.1,1.4) -- (0,1.3) (2,1.3) -- (1.7,1.4); 
			\draw[red] (1.7,1.4) to[out=-135,in=0] (0,1.15) (2,1.15) -- (1.7,1.4) 
			(1.7,1.4) to[out=135,in=10] (0,1.4) (2,1.4) -- (1.7,1.4); 

			\draw[green] (1.7,0) -- (1.7,2)
			(1.6,0) -- (1.7,1.4) (1.1,1.4) -- (1.6,2)
			(1.5,0) -- (1.7,1.4) to[out=-140,in=-20] (1,1.2) to[out=160,in=-140] (1,1.7) to[out=40,in=-150] (1.5,2)
			(1.4,0) -- (1.1,1.4) (1.7,1.4) -- (1.4,2)  
			(1.3,0) to[out=110,in=-160] (1,1.6) to[out=20,in=140] (1.7,1.4) -- (1.3,2);

			\draw (1.1,1.4) to[out=35,in=145] (1.7,1.4) to[out=-145,in=-35] (1.1,1.4);
			\draw[fill] (1.1,1.4) circle (0.06)
			(1.7,1.4) circle (0.06);

		\end{tikzpicture}
		\caption{}
		\label{fig:J12 case 3 first a}
	\end{subfigure}
	\begin{subfigure}[b]{0.3\textwidth}
		\centering
		\begin{tikzpicture}[scale=1.6]
			\draw[dashed] (0,0) rectangle (2,2);

			\draw (1.7,1.4) to[out=140,in=20] (1,1.55) to[out=-160,in=160] (1,1.25) to[out=-20,in=-140] (1.7,1.4);

			\draw (1.1,1.4) -- (0,1.3) (2,1.3) -- (1.7,1.4); 
			\draw[blue] (1.7,1.4) to[out=-135,in=0] (0,1.15) (2,1.15) -- (1.7,1.4); 
			\draw (1.7,1.4) to[out=135,in=10] (0,1.4) (2,1.4) -- (1.7,1.4); 

			\draw[dashed,red] (1.7,0) -- (1.7,2);
			\draw[red]
			(1.6,0) -- (1.7,1.4) (1.1,1.4) -- (1.6,2)
			(1.5,0) -- (1.7,1.4) to[out=-140,in=-20] (1,1.2) to[out=160,in=-140] (1,1.7) to[out=40,in=-150] (1.5,2);

			\draw[green] (1.7,1.4) -- (1.2,2) (1.2,0) -- (0,1.1) (2,1.1) -- (1.7,1.4)
			(1.1,1.4) -- (1,2) (1,0) -- (0,1) (2,1) -- (1.7,1.4)
			(1.7,1.4) to[out=-140,in=-20] (0.95,1.18) to[out=160,in=-110] (0.8,2) (0.8,0) -- (0,0.9) (2,0.9) -- (1.7,1.4);

			\draw (1.1,1.4) to[out=35,in=145] (1.7,1.4) to[out=-145,in=-35] (1.1,1.4);
			\draw[fill] (1.1,1.4) circle (0.06)
			(1.7,1.4) circle (0.06);

			\node[red] at (1.7,2.15) {\( s \)};
			\node[blue] at (-0.15,1.15) {\( t \)};

		\end{tikzpicture}
		\caption{}
		\label{fig:J12 case 3 first b}
	\end{subfigure}

	\begin{subfigure}[b]{0.3\textwidth}
		\centering
		\begin{tikzpicture}[scale=1.6]
			\draw[dashed] (0,0) rectangle (2,2);

			\draw (1.7,1.4) to[out=140,in=20] (1,1.55) to[out=-160,in=160] (1,1.25) to[out=-20,in=-140] (1.7,1.4);

			\draw (1.1,1.4) -- (0,1.3) (2,1.3) -- (1.7,1.4) 
			(1.7,1.4) to[out=-135,in=0] (0,1.15) (2,1.15) -- (1.7,1.4); 
			\draw (1.7,1.4) to[out=135,in=10] (0,1.4) (2,1.4) -- (1.7,1.4); 

			\draw (1.7,0) -- (1.7,2)
			(1.6,0) -- (1.7,1.4) (1.1,1.4) -- (1.6,2)
			(1.5,0) -- (1.7,1.4) to[out=-140,in=-20] (1,1.2) to[out=160,in=-140] (1,1.7) to[out=40,in=-150] (1.5,2);


			\draw[dashed,red] (0.35,0) -- (1.7,1.4) to[out=-140,in=-20] (1,1.22) to[out=160,in=-160] (1,1.6) to[out=20,in=170] (2,1.75) (0,1.75) -- (0.35,2);
			\draw[red] (1.1,1.4) to[out=45,in=180] (2,1.7) (0,1.7) -- (0.4,2) (0.4,0) -- (1.7,1.4);
			\draw[blue] (0,1.6) -- (0.3,2) (0.3,0) -- (1.1,1.4) (1.7,1.4) -- (2,1.6);
			\draw[dashed,blue] (0,1.65) -- (0.25,2) (0.25,0) to[out=80,in=-160] (1,1.57) to[out=20,in=140] (1.7,1.4) -- (2,1.65);

			\draw (1.1,1.4) to[out=35,in=145] (1.7,1.4) to[out=-145,in=-35] (1.1,1.4);
			\draw[fill] (1.1,1.4) circle (0.06)
			(1.7,1.4) circle (0.06);

		\end{tikzpicture}
		\caption{}
		\label{fig:J12 case 3 first c}
	\end{subfigure}
	\begin{subfigure}[b]{0.3\textwidth}
		\centering
		\begin{tikzpicture}[scale=1.6]
			\draw[dashed] (0,0) rectangle (2,2);

			\draw (1.7,1.4) to[out=140,in=20] (1,1.55) to[out=-160,in=160] (1,1.25) to[out=-20,in=-140] (1.7,1.4);

			\draw (1.1,1.4) -- (0,1.3) (2,1.3) -- (1.7,1.4) 
			(1.7,1.4) to[out=-135,in=0] (0,1.15) (2,1.15) -- (1.7,1.4); 
			\draw (1.7,1.4) to[out=135,in=10] (0,1.4) (2,1.4) -- (1.7,1.4); 

			\draw (1.7,0) -- (1.7,2)
			(1.6,0) -- (1.7,1.4) (1.1,1.4) -- (1.6,2)
			(1.5,0) -- (1.7,1.4) to[out=-140,in=-20] (1,1.2) to[out=160,in=-140] (1,1.7) to[out=40,in=-150] (1.5,2);

			\draw[green] (1.7,1.4) -- (2,1) (0,1) -- (1,0) (1,2) -- (1.1,1.4);

			\draw[red] (1.1,1.4) to[out=45,in=180] (2,1.7) (0,1.7) -- (0.4,2) (0.4,0) -- (1.7,1.4)
			(0,1.6) -- (0.3,2) (0.3,0) -- (1.1,1.4) (1.7,1.4) -- (2,1.6);
			\draw[green] (0,1.5) -- (0.5,2) (0.5,0) -- (1.7,1.4) -- (2,1.5);

			\draw (1.1,1.4) to[out=35,in=145] (1.7,1.4) to[out=-145,in=-35] (1.1,1.4);
			\draw[fill] (1.1,1.4) circle (0.06)
			(1.7,1.4) circle (0.06);

		\end{tikzpicture}
		\caption{}
		\label{fig:J12 case 3 first d}
	\end{subfigure}
	\caption{Constructing a 1-system}
	\label{fig:J12 case 3 first}
\end{figure}

\textbf{Step 1.} 
First suppose \( X_0 \) is the set of arcs shown in red in \Cref{fig:J12 case 3 first a}. There are exactly five arcs with slope number \( \infty \) which intersect each arc of \( X_0 \cup Y \) at most once. These are shown in green in \Cref{fig:J12 case 3 first a}. Up to reflection, there are two possibilities for \( X_\infty \) given that \( |X_\infty| = 3 \). 

If \( X_\infty \) is the set of arcs shown in red in \Cref{fig:J12 case 3 first b}, there are three arcs with slope number \( -1 \) that intersects each arc proven to be in \( \mathcal{A} \) at most once. These are shown in green in \Cref{fig:J12 case 3 first b}. We see that no arc in \( X_\infty \cup X_{-1} \cup X_{0} \cup Y \) intersects the arc \( s \) shown in dashed red in or the arc \( t \) shown in blue \Cref{fig:J12 case 3 first b}. So, there must be arcs in \( X_1 \) which intersect \( s \) and \( t \). 

Among arcs that may be included in \( \mathcal{A} \), there are two arcs with slope number 1 which intersect \( s \), shown in red in \Cref{fig:J12 case 3 first c}, and two arcs with slope number 1 which intersect \( t \), shown in blue in \Cref{fig:J12 case 3 first c}. At least one of the blue arcs must be in \( X_1 \) and at least one of the red arcs must be in \( X_1 \). However, note that the dashed red arc intersects both blue arcs twice, so it cannot be in \( X_1 \). Similarly, the dashed blue arc cannot be in \( X_1 \). Therefore, \( X_1 \) contains the arcs shown in red in \Cref{fig:J12 case 3 first d}. Then, there is one arc with slope number -1 and one more arc with slope number 1, shown in green in \Cref{fig:J12 case 3 first d}, which intersect each arc proven to be in \( \mathcal{A} \) at most once. So, we obtain a 1-system.




If \( X_\infty \) is the set of arcs shown in red in \Cref{fig:J12 case 3 first e}, note that the arc \( s \) shown in dashed red intersects no other arc in \( Y \cup X_0 \cup X_{\infty} \). Up to reflection across the horizontal axis, we may assume that there is some arc with slope number 1 in \( \mathcal{A} \) which intersects \( s \). However, any arc with slope number 1 which intersects \( s \) intersects some other arc of \( \mathcal{A} \) at least twice. This is a contradiction.

\textbf{Step 2.}
If \( X_0 \) is the set of arcs shown in red in \Cref{fig:J12 case 3 second}, since \( |X_\infty| =  3 \), it must be that \( X_\infty \) is the set of arcs shown in green in \Cref{fig:J12 case 3 second}. By rotating 90 degrees clockwise, we see that this system is equivalent to a system from step 1. Therefore, we do not obtain a new 1-system in this step.

\begin{figure}
	\centering
	\begin{subfigure}[b]{0.3\textwidth}
		\centering
		\begin{tikzpicture}[scale=1.6]
			\draw[dashed] (0,0) rectangle (2,2);

			\draw (1.7,1.4) to[out=140,in=20] (1,1.55) to[out=-160,in=160] (1,1.25) to[out=-20,in=-140] (1.7,1.4);

			\draw (1.1,1.4) -- (0,1.3) (2,1.3) -- (1.7,1.4); 
			\draw (1.7,1.4) to[out=-135,in=0] (0,1.15) (2,1.15) -- (1.7,1.4) 
			(1.7,1.4) to[out=135,in=10] (0,1.4) (2,1.4) -- (1.7,1.4); 

			\draw[dashed,red] (1.7,0) -- (1.7,2);
			\draw[red]
			(1.6,0) -- (1.7,1.4) (1.1,1.4) -- (1.6,2)
			(1.4,0) -- (1.1,1.4) (1.7,1.4) -- (1.4,2)  
			;

			\draw (1.1,1.4) to[out=35,in=145] (1.7,1.4) to[out=-145,in=-35] (1.1,1.4);
			\draw[fill] (1.1,1.4) circle (0.06)
			(1.7,1.4) circle (0.06);

			\node[red] at (1.7,2.15) {\( s \)};

		\end{tikzpicture}
		\caption{}
		\label{fig:J12 case 3 first e}
	\end{subfigure}
	\begin{subfigure}{0.3\textwidth}
		\centering
		\begin{tikzpicture}[scale=1.6]
			\draw[dashed] (0,0) rectangle (2,2);

			\draw (1.7,1.4) to[out=140,in=20] (1,1.55) to[out=-160,in=160] (1,1.25) to[out=-20,in=-140] (1.7,1.4);

			\draw[red] (1.1,1.4) -- (0,1.3) (2,1.3) -- (1.7,1.4) 
			(1.7,1.4) to[out=135,in=20] (0,1.4) (2,1.4) -- (1.7,1.4) 
			(0,1.35) to[out=20,in=135] (1.7,1.4) (1.1,1.4) to[out=-40,in=-160] (2,1.35)
			;

			\draw[green] (1.7,0) -- (1.7,2)
			(1.4,0) -- (1.1,1.4) (1.7,1.4) -- (1.4,2)  
			(1.3,0) to[out=110,in=-160] (1,1.6) to[out=20,in=140] (1.7,1.4) -- (1.3,2)
			;

			\draw (1.1,1.4) to[out=35,in=145] (1.7,1.4) to[out=-145,in=-35] (1.1,1.4);
			\draw[fill] (1.1,1.4) circle (0.06)
			(1.7,1.4) circle (0.06);

		\end{tikzpicture}
		\caption{}
		\label{fig:J12 case 3 second}
	\end{subfigure}
	\caption{Obtaining contradictions}
	\label{fig:J12 case 3 contras}
\end{figure}
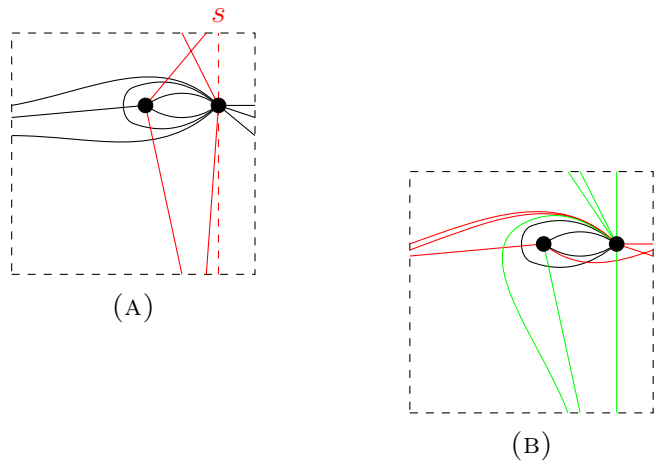

\section{$|J| = 0$}
\label{sec:J0}

In this case, we choose an arc $u \in\mathcal{A}$ which is \emph{minimally intersected} in the sense it minimizes the value \( |\left\{ v \in \mathcal{A} : i(u,v) = 1 \right\}| \) with respect to \( u \). 


We split the proof depending on whether $u$ is a loop arc or $u$ is a non-loop arc.

\subsection{\( u \) is a non-loop arc}

\begin{figure}
	\centering
	\begin{subfigure}{0.24\textwidth}
		\centering
		\begin{tikzpicture}
			\node[inner sep=0pt] at (0.65,1.35) (a) {};
			\node[inner sep=0pt] at (1.35,0.65) (b) {};
			\draw[dashed] (0,0) -- (2,0) -- (2,2) -- (0,2) -- (0,0);
			\draw[fill] (a) circle (0.05) (b) circle (0.05);

			\draw (a) -- (b);
			\draw[blue]
			(a) -- (0.75,2) (0.75,0) to[out=80,in=-135] (1,1) to[out=45,in=190] (2,1.25) (0,1.25) -- (a);
		\end{tikzpicture}
		\caption{a loop arc}
		\label{fig:nl intro a}
	\end{subfigure}
	\begin{subfigure}{0.24\textwidth}
		\centering
		\begin{tikzpicture}
			\fill[color=white!70!red] (0.75,1.25) -- (0.8,2) -- (1.2,2) -- (1,1);
			\fill[color=white!70!red] (1.25,0.75) -- (1.2,0) -- (0.8,0) -- (1,1);
			\draw[dashed] (0,0) -- (2,0) -- (2,2) -- (0,2) -- (0,0);
			\draw[fill] (0.75,1.25) circle (0.05) (1.25,0.75) circle (0.05);

			\draw (0.75, 1.25) -- (1.25,0.75);
			\draw[blue]
			(0.75,1.25) -- (0.8,2) (0.8,0) -- (1.2,2) (1.2,0) -- (1.25,0.75);
		\end{tikzpicture}
		\caption{a non-loop arc}
		\label{fig:nl intro b}
	\end{subfigure}
	\begin{subfigure}{0.24\textwidth}
		\centering
		\begin{tikzpicture}
			\draw[dashed] (0,0) -- (2,0) -- (2,2) -- (0,2) -- (0,0);

			\fill[color=white!70!red] (0.75,1.25) -- (0,1.1) -- (0,0.9);
			\fill[color=white!70!red] (2,1.1) to[out=-180,in=90] (1.1,0.75) to[out=-90,in=-150] (2,0.9);

			\draw[blue] (0.75,1.25) -- (0,1.1) (2,1.1) to[out=-180,in=90] (1.1,0.75) to[out=-90,in=-150] (2,0.9) (0,0.9) -- (0.75,1.25);

			\draw (0.75, 1.25) -- (1.25,0.75);

			\draw[fill] (0.75,1.25) circle (0.05) (1.25,0.75) circle (0.05);
		\end{tikzpicture}
		\caption{a loop arc}
		\label{fig:nl intro c}
	\end{subfigure}
	\caption{}
	\label{fig:nl intro}
\end{figure}

By assumption, \( J \) is empty, so there is some arc \( v \) which intersects \( u \). Up to homeomorphism, \( v \) must be one of the arcs shown in blue in \Cref{fig:nl intro}. 

We claim that \( v \) must be the arc shown in blue in \Cref{fig:nl intro a}. Suppose instead that \( v \) is the arc shown in \Cref{fig:nl intro b} and consider the almost embedded disk shown in red in Figure~\ref{fig:nl intro b}. We may apply Corollary~\ref{cor:twoss} to this disk, and we find that the internal arc of the disk must be included in $\mathcal{A}$. In addition, by applying Lemma~\ref{lem:twos} we can see that any arc which intersects this arc must also intersect the arc $u$, which contradicts the assumption that $u$ is minimally intersected. 

Similarly, suppose \( v \) is the arc shown in \Cref{fig:nl intro c}. We apply \Cref{cor:help ones} to the disk shown in red in \Cref{fig:nl intro c}. This would imply that the internal arc of that disk is in \( J \), which contradicts our assumption that \( |J| = 0 \). 

We conclude that \( v \) is the arc shown in \Cref{fig:nl intro a}. 

\begin{figure}
	\centering
	\begin{subfigure}{0.24\textwidth}
		\centering
		\begin{tikzpicture}
			\fill[color=white!70!red] (1,1) to[out=-100,in=25] (0,0) to[out=75,in=-135] (0.7,1.3);
			\draw[blue] (0,0) -- (1,1);
			\draw[dashed] (0,0) -- (2,0) -- (2,2) -- (0,2) -- (0,0);
			\draw[fill] (0,0) circle (0.05) 
			(2,2) circle (0.05)
			(0,2) circle (0.05)
			(2,0) circle (0.05)
			(1,1) circle (0.05);

			\draw (0,2) -- (1,1);
			\draw (0,0) to[out=75,in=-135] (0.7,1.3) to[out=45,in=-165] (2,2);

			\node at (0.5,1.7) {\( u \)};
			\node at (0.9,1.7) {\( v \)};
		\end{tikzpicture}
		\caption{}
		\label{fig:nl l first a}
	\end{subfigure}
	\begin{subfigure}{0.24\textwidth}
		\centering
		\begin{tikzpicture}
			\draw[blue] (0,0) -- (1,1);
			\draw[red] (0.3,1.1) -- (0.7,0.4) (0.9,0.6) to[out=135,in=-145] (1,1.3);
			\draw[dashed] (0,0) -- (2,0) -- (2,2) -- (0,2) -- (0,0);
			\draw[fill] (0,0) circle (0.05) 
			(2,2) circle (0.05)
			(0,2) circle (0.05)
			(2,0) circle (0.05)
			(1,1) circle (0.05);

			\draw (0,2) -- (1,1);
			\draw (0,0) to[out=75,in=-135] (0.7,1.3) to[out=45,in=-165] (2,2);
		\end{tikzpicture}
		\caption{}
		\label{fig:nl l first b}
	\end{subfigure}
	\begin{subfigure}{0.24\textwidth}
		\centering
		\begin{tikzpicture}
			\draw[blue] (0,0) -- (1,1);
			\draw[red] (0,1.3) to[out=-45,in=160] (2,0) to[out=145,in=-40] (1,0.7) to[out=140,in=-150] (1,1.3) to[out=30,in=160] (2,1.3);
			\draw[dashed] (0,0) -- (2,0) -- (2,2) -- (0,2) -- (0,0);

			\draw[fill] (0,0) circle (0.05) 
			(2,2) circle (0.05)
			(0,2) circle (0.05)
			(2,0) circle (0.05)
			(1,1) circle (0.05);

			\draw (0,2) -- (1,1);
			\draw (0,0) to[out=75,in=-135] (0.7,1.3) to[out=45,in=-165] (2,2);
		\end{tikzpicture}
		\caption{}
		\label{fig:nl l first c}
	\end{subfigure}
	\begin{subfigure}{0.24\textwidth}
		\centering
		\begin{tikzpicture}
			\draw[blue] (0,0) -- (1,1);
			\draw[red] (2,2) to[out=-155,in=150] (0.9,0.75) to[out=-30,in=-160] (2,0.7) (0,0.7) -- (2,0); 
			\draw[dashed] (0,0) -- (2,0) -- (2,2) -- (0,2) -- (0,0);

			\draw[fill] (0,0) circle (0.05) 
			(2,2) circle (0.05)
			(0,2) circle (0.05)
			(2,0) circle (0.05)
			(1,1) circle (0.05);

			\draw (0,2) -- (1,1);
			\draw (0,0) to[out=75,in=-135] (0.7,1.3) to[out=45,in=-165] (2,2);
		\end{tikzpicture}
		\caption{}
		\label{fig:nl l first d}
	\end{subfigure}

	\begin{subfigure}{0.24\textwidth}
		\centering
		\begin{tikzpicture}
			\fill[color=white!70!red] (0,0) to[out=75,in=-125] (0.385,0.935) to[out=-40,in=160] (2,0);
			\fill[color=white!70!red] (0,2) -- (0.7,1.3) to[out=45,in=-165] (2,2);

			\draw[blue] (0,0) -- (1,1);
			\draw[red] (0,1.3) to[out=-45,in=160] (2,0) to[out=145,in=-40] (1,0.7) to[out=140,in=-150] (1,1.3) to[out=30,in=160] (2,1.3);
			\draw[dashed] (2,0) -- (2,2) (0,2) -- (0,0);
			\draw[green] (0,0) -- (2,0) (0,2) -- (2,2);

			\draw[fill] (0,0) circle (0.05) 
			(2,2) circle (0.05)
			(0,2) circle (0.05)
			(2,0) circle (0.05)
			(1,1) circle (0.05);

			\draw (0,2) -- (1,1);
			\draw (0,0) to[out=75,in=-135] (0.7,1.3) to[out=45,in=-165] (2,2);
		\end{tikzpicture}
		\caption{}
		\label{fig:nl l first e}
	\end{subfigure}
	\begin{subfigure}{0.24\textwidth}
		\centering
		\begin{tikzpicture}
			\fill[color=white!70!red] (0,0) to[out=75,in=-135] (0.7,1.3) -- (0.87,1.14) to[out=-105,in=140] (1,0.7) to[out=-40,in=145] (2,0);
			\fill[color=white!70!red] (0,2) -- (0.7,1.3) to[out=45,in=-165] (2,2);

			\draw[green] (1,0.6) -- (1,0) (1,2) -- (1,1.4);

			\draw[blue] (0,0) -- (1,1);
			\draw[red] (0,1.3) to[out=-45,in=160] (2,0) to[out=145,in=-40] (1,0.7) to[out=140,in=-150] (1,1.3) to[out=30,in=160] (2,1.3);
			\draw[dashed] (2,0) -- (2,2) (0,2) -- (0,0);
			\draw[blue] (0,0) -- (2,0) (0,2) -- (2,2);

			\draw[fill] (0,0) circle (0.05) 
			(2,2) circle (0.05)
			(0,2) circle (0.05)
			(2,0) circle (0.05)
			(1,1) circle (0.05);

			\draw (0,2) -- (1,1);
			\draw (0,0) to[out=75,in=-135] (0.7,1.3) to[out=45,in=-165] (2,2);
		\end{tikzpicture}
		\caption{}
		\label{fig:nl l first f}
	\end{subfigure}
	\begin{subfigure}{0.24\textwidth}
		\centering
		\begin{tikzpicture}
			\fill[color=white!70!red] (0,0) to[out=75,in=-117] (0.21,0.62) -- (2,0);
			\fill[color=white!70!red] (0,2) -- (0.7,1.3) to[out=45,in=-165] (2,2);

			\draw[blue] (0,0) -- (1,1);
			\draw[red] (2,2) to[out=-155,in=150] (0.9,0.75) to[out=-30,in=-160] (2,0.7) (0,0.7) -- (2,0); 
			\draw[dashed] (2,0) -- (2,2) (0,2) -- (0,0);
			\draw[green] (0,0) -- (2,0) (0,2) -- (2,2);

			\draw[fill] (0,0) circle (0.05) 
			(2,2) circle (0.05)
			(0,2) circle (0.05)
			(2,0) circle (0.05)
			(1,1) circle (0.05);

			\draw (0,2) -- (1,1);
			\draw (0,0) to[out=75,in=-135] (0.7,1.3) to[out=45,in=-165] (2,2);
		\end{tikzpicture}
		\caption{}
		\label{fig:nl l first g}
	\end{subfigure}
	\begin{subfigure}{0.24\textwidth}
		\centering
		\begin{tikzpicture}
			\fill[color=white!70!red] (0,0) to[out=75,in=-117] (0.21,0.62) -- (2,0);
			\fill[color=white!70!red] (0,2) -- (0.82,1.18) to[out=50,in=-155] (2,2);

			\draw[green] (1,0.6) -- (1,0) (1,2) -- (1,1.5);

			\draw[blue] (0,0) -- (1,1);
			\draw[red] (2,2) to[out=-155,in=150] (0.9,0.75) to[out=-30,in=-160] (2,0.7) (0,0.7) -- (2,0); 
			\draw[dashed] (2,0) -- (2,2) (0,2) -- (0,0);
			\draw[blue] (0,0) -- (2,0) (0,2) -- (2,2);

			\draw[fill] (0,0) circle (0.05) 
			(2,2) circle (0.05)
			(0,2) circle (0.05)
			(2,0) circle (0.05)
			(1,1) circle (0.05);

			\draw (0,2) -- (1,1);
			\draw (0,0) to[out=75,in=-135] (0.7,1.3) to[out=45,in=-165] (2,2);
		\end{tikzpicture}
		\caption{}
		\label{fig:nl l first h}
	\end{subfigure}
	\caption{}
	\label{fig:nl l first}
\end{figure}

Let \( w \) be the arc shown in blue in \Cref{fig:nl l first a}. 

\begin{claim}
	\label{claim:sec5 some arcs}
	$w \in \mathcal{A}$
\end{claim}

\begin{proof}
	Suppose for contradiction that $w \notin \mathcal{A}$. Since $\mathcal{A}$ is saturated, there must be some $w' \in \mathcal{A}$ which intersects $w$ at least twice. We may apply Lemma~\ref{lem:twos} to the disk shown in red in \Cref{fig:nl l first a} and we determine two subarcs of \( w' \). These are shown in red in Figure~\ref{fig:nl l first b}. There are only two possiblilities for $w'$, as shown in Figure~\ref{fig:nl l first c} and Figure~\ref{fig:nl l first d}. 

	In both cases, we can apply Corollary~\ref{cor:twoss} to a certain disk, shown in red in \Cref{fig:nl l first e,fig:nl l first g}, and we conclude that an additional arc \( a \), shown in green in \Cref{fig:nl l first e,fig:nl l first g}, must be included in $\mathcal{A}$. However, we will show that any arc which intersects \( a \) also intersects the arc \( u \), contradicting the assumption that \( u \) is minimally intersected.

	To see this, first suppose \( w' \) is the arc shown in red in \Cref{fig:nl l first c} and suppose there was some arc \( a' \in \mathcal{A} \) which intersects \( a \) but not \( u \). We apply \Cref{lem:twos} to the disk shown in red in \Cref{fig:nl l first e} to determine a subarc of \( a' \), shown in green in \Cref{fig:nl l first f}. Now, consider the disc shown in red in \Cref{fig:nl l first f}. The subarc we determined intersects the boundary of this disk in exactly one point. However, every point in the boundary of the disk is either contained in \( u \), which by assumption does not intersect \( a' \), or contained in \( v \) or \( w' \), both of which are intersected by the subarc we determined above. This contradicts \Cref{lem:twos}, so no such arc \( a' \) can be in \( \mathcal{A} \). 

	This contradicts the assumption that \( u \) is minimally intersected. So, it cannot be that \( w' \in \mathcal{A} \). 

	We have a similar argument if \( w' \) is the arc shown in \Cref{fig:nl l first d}, using the disks shown in red in \Cref{fig:nl l first g,fig:nl l first h}. So, it must be that \( w \in \mathcal{A} \). 
\end{proof}

\begin{figure}
	\centering
	\begin{subfigure}{0.3\textwidth}
		\centering
		\begin{tikzpicture}
			\draw[blue] (0,0) -- (2,2);
			\draw[dashed] (0,0) -- (2,0) -- (2,2) -- (0,2) -- (0,0);
			\draw[fill] (0,0) circle (0.05) 
			(2,2) circle (0.05)
			(0,2) circle (0.05)
			(2,0) circle (0.05)
			(1,1) circle (0.05);

			\draw (0,2) -- (1,1);
			\draw (0,0) to[out=75,in=-135] (0.7,1.3) to[out=45,in=-165] (2,2);
			\draw[fill,white] (0.5,0.5) circle (0.15)
			(1.5,1.5) circle (0.15);

			\node[blue] at (0.5,0.5) {$w$};
			\node[blue] at (1.5,1.5) {$x$};
		\end{tikzpicture}
		\caption{}
		\label{fig:nl l second a}
	\end{subfigure}
	\begin{subfigure}{0.3\textwidth}
		\centering
		\begin{tikzpicture}
			\fill[color=white!70!red] (1,1) to[out=-100,in=25] (0,0) to[out=75,in=-135] (0.7,1.3);

			\draw[red] (0.2,1) -- (0.9,0.5);

			\draw[blue] (0,0) -- (2,2);
			\draw[dashed] (0,0) -- (2,0) -- (2,2) -- (0,2) -- (0,0);
			\draw[fill] (0,0) circle (0.05) 
			(2,2) circle (0.05)
			(0,2) circle (0.05)
			(2,0) circle (0.05)
			(1,1) circle (0.05);

			\draw (0,2) -- (1,1);
			\draw (0,0) to[out=75,in=-135] (0.7,1.3) to[out=45,in=-165] (2,2);
		\end{tikzpicture}
		\caption{}
		\label{fig:nl l second b}
	\end{subfigure}
	\begin{subfigure}{0.3\textwidth}
		\centering
		\begin{tikzpicture}

			\draw[red] (2,0) -- (0,0.7) (2,0.7) -- (1,1);

			\draw[blue] (0,0) -- (2,2);
			\draw[dashed] (0,0) -- (2,0) -- (2,2) -- (0,2) -- (0,0);
			\draw[fill] (0,0) circle (0.05) 
			(2,2) circle (0.05)
			(0,2) circle (0.05)
			(2,0) circle (0.05)
			(1,1) circle (0.05);

			\draw (0,2) -- (1,1);
			\draw (0,0) to[out=75,in=-135] (0.7,1.3) to[out=45,in=-165] (2,2);

			\draw[fill,white] (1,0.4) circle (0.2);
			\node[red] at (1,0.4) {$w'_1$};
		\end{tikzpicture}
		\caption{}
		\label{fig:nl l second c}
	\end{subfigure}
	\begin{subfigure}{0.3\textwidth}
		\centering
		\begin{tikzpicture}

			\draw[red] (2,0) to[out=165,in=-75] (0,2);

			\draw[blue] (0,0) -- (2,2);
			\draw[dashed] (0,0) -- (2,0) -- (2,2) -- (0,2) -- (0,0);
			\draw[fill] (0,0) circle (0.05) 
			(2,2) circle (0.05)
			(0,2) circle (0.05)
			(2,0) circle (0.05)
			(1,1) circle (0.05);

			\draw (0,2) -- (1,1);
			\draw (0,0) to[out=75,in=-135] (0.7,1.3) to[out=45,in=-165] (2,2);

			\draw[fill,white] (1,0.4) circle (0.2);
			\node[red] at (1,0.4) {$w'_2$};
		\end{tikzpicture}
		\caption{}
		\label{fig:nl l second d}
	\end{subfigure}
	\begin{subfigure}{0.3\textwidth}
		\centering
		\begin{tikzpicture}

			\draw[red] (2,2) to[out=-115,in=0] (1,0.8) to[out=180,in=-65] (0,2);

			\draw[blue] (0,0) -- (2,2);
			\draw[dashed] (0,0) -- (2,0) -- (2,2) -- (0,2) -- (0,0);
			\draw[fill] (0,0) circle (0.05) 
			(2,2) circle (0.05)
			(0,2) circle (0.05)
			(2,0) circle (0.05)
			(1,1) circle (0.05);

			\draw (0,2) -- (1,1);
			\draw (0,0) to[out=75,in=-135] (0.7,1.3) to[out=45,in=-165] (2,2);

			\draw[fill,white] (1.4,1) circle (0.2);
			\node[red] at (1.4,1) {$w'_3$};
		\end{tikzpicture}
		\caption{}
		\label{fig:nl l second e}
	\end{subfigure}
	\begin{subfigure}{0.3\textwidth}
		\centering
		\begin{tikzpicture}
			\draw[red] (2,2) to[out=-115,in=20] (1.1,0.8) to[out=-160,in=0] (0,0.6) (2,0.6) -- (0,0);

			\draw[blue] (0,0) -- (2,2);
			\draw[dashed] (0,0) -- (2,0) -- (2,2) -- (0,2) -- (0,0);
			\draw[fill] (0,0) circle (0.05) 
			(2,2) circle (0.05)
			(0,2) circle (0.05)
			(2,0) circle (0.05)
			(1,1) circle (0.05);

			\draw (0,2) -- (1,1);
			\draw (0,0) to[out=75,in=-135] (0.7,1.3) to[out=45,in=-165] (2,2);

			\draw[fill,white] (1.4,1) circle (0.2);
			\node[red] at (1.4,1) {$w'_4$};
		\end{tikzpicture}
		\caption{}
		\label{fig:nl l second f}
	\end{subfigure}
	\caption{}
	\label{fig:nl l second}
\end{figure}

By symmetry, it must be that the arc $x$, shown in Figure~\ref{fig:nl l second a}, is also included in $\mathcal{A}$. 

By the assumption that $u$ is minimally intersected, there must be some arc $w' \in \mathcal{A}$ which intersects $w$ but not $u$. We may apply Lemma~\ref{lem:twos} to the disk shown in red in \Cref{fig:nl l second b} in order to determine a subarc of \( w' \). This subarc is shown in red in Figure~\ref{fig:nl l second b}. We see that \( w' \) must be one of the four arcs shown in \Crefrange{fig:nl l second c}{fig:nl l second f}. Call these arcs \( w'_1, w'_2, w'_3, w'_4 \) as labelled in the figure.

Again by symmetry, there must be some arc $x' \in \mathcal{A}$ which intersects $x$ but not $u$. For each $i$, let $x'_i$ be the arc obtained from $w'_i$ by reflection over the line containing the arc \( u \) in \Cref{fig:nl l second}.

Up to reflection, we can assume that $w'_i, x'_j \in \mathcal{A}$ with $i \le j$.

Now note the following: $w'_1$ intersects at least twice each $x'_1$, $x'_3$, and $x'_4$. $w'_2$ intersects twice $x'_4$, and $w'_4$ intersecs twice $x'_4$. 

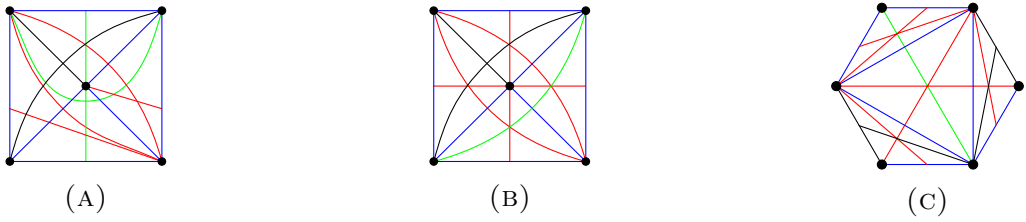
\begin{figure}
	\centering
	\begin{subfigure}{0.24\textwidth}
		\centering
		\begin{tikzpicture}
			\draw[red] (2,0) to[out=155,in=-45] (0.7,0.7) to[out=135,in=-75] (0,2);
			\draw[green] (0,2) to[out=-65,in=180] (1,0.8) to[out=0,in=-105] (2,2);
			\draw[green] (1,0) -- (1,2);

			\draw[red] (2,0) -- (0,0.7) (2,0.7) -- (1,1);
			\draw[red] (2,0) to[out=105,in=-15] (0,2);

			\draw[blue] (0,0) -- (2,2);
			\draw[blue] (0,0) -- (2,0) -- (2,2) -- (0,2) -- (0,0);
			\draw[blue](1,1) -- (2,0);

			\draw[fill] (0,0) circle (0.05) 
			(2,2) circle (0.05)
			(0,2) circle (0.05)
			(2,0) circle (0.05)
			(1,1) circle (0.05);

			\draw (0,2) -- (1,1);
			\draw (0,0) to[out=75,in=-135] (0.7,1.3) to[out=45,in=-165] (2,2);
		\end{tikzpicture}
		\caption{}
	\end{subfigure}
	\begin{subfigure}{0.4\textwidth}
		\centering
		\begin{tikzpicture}
			\draw[red] (0,1) -- (2,1)
			(1,0) -- (1,2);
			\draw[red] (2,0) to[out=105,in=-15] (0,2);
			\draw[red] (2,0) to[out=165,in=-75] (0,2);

			\draw[blue] (0,0) -- (2,2);
			\draw[blue] (0,0) -- (2,0) -- (2,2) -- (0,2) -- (0,0);
			\draw[blue] (2,0) -- (1,1);

			\draw[green] (0,0) to[out=15,in=-105] (2,2);

			\draw[fill] (0,0) circle (0.05) 
			(2,2) circle (0.05)
			(0,2) circle (0.05)
			(2,0) circle (0.05)
			(1,1) circle (0.05);

			\draw (0,2) -- (1,1);
			\draw (0,0) to[out=75,in=-135] (0.7,1.3) to[out=45,in=-165] (2,2);
		\end{tikzpicture}
		\caption{}
	\end{subfigure}
	\begin{subfigure}{0.23\textwidth}
		\centering
		\begin{tikzpicture}[scale=1.2]
			\draw[green] (-60:1) -- (120:1);

			\draw[red] (60:1) -- (-120:1)
			(0:1) -- (180:1);

			\draw[blue] (-60:1) -- (180:1) -- (60:1) -- (-60:1);

			\draw[blue] (0:1) -- (-60:1) -- (-120:1) (180:1) -- (120:1) -- (60:1);
			\draw[red] (150:0.87) -- (60:1) -- (-30:0.87)
			(-90:0.87) -- (180:1) -- (90:0.87);

			\draw (0:1) -- (60:1) (180:1) -- (-120:1);
			\draw (30:0.87) -- (-60:1) -- (-150:0.87);

			\foreach \x in {0,60,...,300} {
				\draw[fill] (\x:1) circle (0.05);
			}
		\end{tikzpicture}
		\caption{}
	\end{subfigure}
	\caption{The 3 1-systems obtained when \( |J| = 0 \).}
	\label{fig:J0 finals}
\end{figure}

So, we split the proof in 5 steps. We will obtain three non-equivalent 1-systems, shown in \Cref{fig:J0 finals}.

\textbf{Step 1.} In this step, we assume $w_1', x_2' \in \mathcal{A}$. This is shown in Figure~\ref{fig:nl l step 1 z}.

\begin{figure}
	\centering
	\begin{subfigure}{0.3\textwidth}
		\centering
		\begin{tikzpicture}
			\draw[red] (2,0) -- (0,0.7) (2,0.7) -- (1,1);
			\draw[red] (2,0) to[out=105,in=-15] (0,2);

			\draw[blue] (0,0) -- (2,2);
			\draw[dashed] (0,0) -- (2,0) -- (2,2) -- (0,2) -- (0,0);
			\draw[fill] (0,0) circle (0.05) 
			(2,2) circle (0.05)
			(0,2) circle (0.05)
			(2,0) circle (0.05)
			(1,1) circle (0.05);

			\draw (0,2) -- (1,1);
			\draw (0,0) to[out=75,in=-135] (0.7,1.3) to[out=45,in=-165] (2,2);
		\end{tikzpicture}
		\caption{}
		\label{fig:nl l step 1 z}
	\end{subfigure}
	\begin{subfigure}{0.6\textwidth}
		\centering
		\begin{tikzpicture}
			\fill[color=white!70!red] (0,0) to[out=75,in=-117] (0.21,0.62) -- (2,0);
			\fill[color=white!70!red] (0,2) -- (0.7,1.3) to[out=45,in=-165] (2,2);

			\draw[red] (2,0) -- (0,0.7) (2,0.7) -- (1,1);
			\draw[red] (2,0) to[out=105,in=-15] (0,2);

			\draw[blue] (0,0) -- (2,2);
			\draw[dashed] (0,0) -- (2,0) -- (2,2) -- (0,2) -- (0,0);
			\draw[fill] (0,0) circle (0.05) 
			(2,2) circle (0.05)
			(0,2) circle (0.05)
			(2,0) circle (0.05)
			(1,1) circle (0.05);

			\draw (0,2) -- (1,1);
			\draw (0,0) to[out=75,in=-135] (0.7,1.3) to[out=45,in=-165] (2,2);
		\end{tikzpicture}
		\begin{tikzpicture}
			\fill[color=white!70!red] (2,0) -- (0.21,0.62) to[out=63,in=-140] (0.98,1.55) to[out=-36,in=105] (2,0);

			\draw[red] (2,0) -- (0,0.7) (2,0.7) -- (1,1);
			\draw[red] (2,0) to[out=105,in=-15] (0,2);

			\draw[blue] (0,0) -- (2,2);
			\draw[dashed] (0,0) -- (2,0) -- (2,2) -- (0,2) -- (0,0);
			\draw[fill] (0,0) circle (0.05) 
			(2,2) circle (0.05)
			(0,2) circle (0.05)
			(2,0) circle (0.05)
			(1,1) circle (0.05);

			\draw (0,2) -- (1,1);
			\draw (0,0) to[out=75,in=-135] (0.7,1.3) to[out=45,in=-165] (2,2);
		\end{tikzpicture}
		\begin{tikzpicture}
			\fill[color=white!70!red] (0,0) to[out=75,in=-135] (0.7,1.3) -- (0,2);
			\fill[color=white!70!red] (2,0) to[out=105,in=-38] (0.98,1.55) to[out=38,in=-165] (2,2);

			\draw[red] (2,0) -- (0,0.7) (2,0.7) -- (1,1);
			\draw[red] (2,0) to[out=105,in=-15] (0,2);

			\draw[blue] (0,0) -- (2,2);
			\draw[dashed] (0,0) -- (2,0) -- (2,2) -- (0,2) -- (0,0);
			\draw[fill] (0,0) circle (0.05) 
			(2,2) circle (0.05)
			(0,2) circle (0.05)
			(2,0) circle (0.05)
			(1,1) circle (0.05);

			\draw (0,2) -- (1,1);
			\draw (0,0) to[out=75,in=-135] (0.7,1.3) to[out=45,in=-165] (2,2);
		\end{tikzpicture}
		\caption{}
		\label{fig:nl l step 1 a}
	\end{subfigure}
	\caption{}
	\label{fig:nl l step 1}
\end{figure}

We may apply \Cref{cor:twoss} to the three disks shown in red in \Cref{fig:nl l step 1 a} to obtain three arcs which must be in \( \mathcal{A} \). Call these arcs $a, b, c$ as shown in \Cref{fig:nl l step 1 second z}. 

\begin{figure}
	\centering
	\begin{subfigure}{0.3\textwidth}
		\centering
		\begin{tikzpicture}
			\draw[red] (2,0) -- (0,0.7) (2,0.7) -- (1,1);
			\draw[red] (2,0) to[out=105,in=-15] (0,2);

			\draw[blue] (0,0) -- (2,2);
			\draw[blue] (0,0) -- (2,0) -- (2,2) -- (0,2) -- (0,0);
			\draw[blue](1,1) -- (2,0);

			\draw[fill] (0,0) circle (0.05) 
			(2,2) circle (0.05)
			(0,2) circle (0.05)
			(2,0) circle (0.05)
			(1,1) circle (0.05);

			\draw (0,2) -- (1,1);
			\draw (0,0) to[out=75,in=-135] (0.7,1.3) to[out=45,in=-165] (2,2);

			\node[blue] at (1,0.1) {$a$};
			\draw[fill,white] (1.5,0.5) circle (0.15);
			\node[blue] at (1.5,0.5) {$b$};
			\draw[fill,white] (2,1.2) circle (0.15);
			\node[blue] at (2,1.2) {$c$};
		\end{tikzpicture}
		\caption{}
		\label{fig:nl l step 1 second z}
	\end{subfigure}
	\begin{subfigure}{0.3\textwidth}
		\centering
		\begin{tikzpicture}
			\fill[color=white!70!red] (0,0) to[out=75,in=-117] (0.21,0.62) -- (2,0);
			\fill[color=white!70!red] (0,2) -- (0.7,1.3) to[out=45,in=-165] (2,2);

			\draw[green] (1.3,1.5) -- (1.2,2) (1.2,0) -- (1.1,0.5);

			\draw[red] (2,0) -- (0,0.7) (2,0.7) -- (1,1);
			\draw[red] (2,0) to[out=105,in=-15] (0,2);

			\draw[blue] (0,0) -- (2,2);
			\draw[blue] (0,0) -- (2,0) -- (2,2) -- (0,2) -- (0,0);
			\draw[blue](1,1) -- (2,0);

			\draw[fill] (0,0) circle (0.05) 
			(2,2) circle (0.05)
			(0,2) circle (0.05)
			(2,0) circle (0.05)
			(1,1) circle (0.05);

			\draw (0,2) -- (1,1);
			\draw (0,0) to[out=75,in=-135] (0.7,1.3) to[out=45,in=-165] (2,2);
		\end{tikzpicture}
		\caption{}
		\label{fig:nl l step 1 second a}
	\end{subfigure}
	\begin{subfigure}{0.3\textwidth}
		\centering
		\begin{tikzpicture}
			\draw[green] (1,0) -- (1,2);

			\draw[red] (2,0) -- (0,0.7) (2,0.7) -- (1,1);
			\draw[red] (2,0) to[out=105,in=-15] (0,2);

			\draw[blue] (0,0) -- (2,2);
			\draw[blue] (0,0) -- (2,0) -- (2,2) -- (0,2) -- (0,0);
			\draw[blue](1,1) -- (2,0);

			\draw[fill] (0,0) circle (0.05) 
			(2,2) circle (0.05)
			(0,2) circle (0.05)
			(2,0) circle (0.05)
			(1,1) circle (0.05);

			\draw (0,2) -- (1,1);
			\draw (0,0) to[out=75,in=-135] (0.7,1.3) to[out=45,in=-165] (2,2);
		\end{tikzpicture}
		\caption{}
		\label{fig:nl l step 1 second b}
	\end{subfigure}
	\caption{}
	\label{fig:nl l step 1 second}
\end{figure}

By assumption, $u$ is minimally intersected, so there is some arc $a' \in \mathcal{A}$ which intersects $a$ but not $u$. We apply Lemma~\ref{lem:twos} to determine a subarc of $a'$, shown in green in Figure~\ref{fig:nl l step 1 second a}. Then, \( a' \) must be the arc shown in \Cref{fig:nl l step 1 second b}. 

\begin{figure}
	\centering
	\begin{subfigure}{0.24\textwidth}
		\centering
		\begin{tikzpicture}
			\fill[color=white!70!red] (2,0) -- (0.21,0.62) to[out=63,in=-140] (0.98,1.55) to[out=-36,in=105] (2,0);

			\draw[red] (2,0) -- (0,0.7) (2,0.7) -- (1,1);
			\draw[red] (2,0) to[out=105,in=-15] (0,2);

			\draw[blue] (0,0) -- (2,2);
			\draw[blue] (0,0) -- (2,0) -- (2,2) -- (0,2) -- (0,0);
			\draw[blue](1,1) -- (2,0);

			\draw[fill] (0,0) circle (0.05) 
			(2,2) circle (0.05)
			(0,2) circle (0.05)
			(2,0) circle (0.05)
			(1,1) circle (0.05);

			\draw (0,2) -- (1,1);
			\draw (0,0) to[out=75,in=-135] (0.7,1.3) to[out=45,in=-165] (2,2);
		\end{tikzpicture}
		\caption{}
		\label{fig:nl l step 1 third a}
	\end{subfigure}
	\begin{subfigure}{0.24\textwidth}
		\centering
		\begin{tikzpicture}
			\draw[green] (1,0.2) -- (1.9,0.5);
			\draw[red] (2,0) -- (0,0.7) (2,0.7) -- (1,1);
			\draw[red] (2,0) to[out=105,in=-15] (0,2);

			\draw[blue] (0,0) -- (2,2);
			\draw[blue] (0,0) -- (2,0) -- (2,2) -- (0,2) -- (0,0);
			\draw[blue](1,1) -- (2,0);

			\draw[fill] (0,0) circle (0.05) 
			(2,2) circle (0.05)
			(0,2) circle (0.05)
			(2,0) circle (0.05)
			(1,1) circle (0.05);

			\draw (0,2) -- (1,1);
			\draw (0,0) to[out=75,in=-135] (0.7,1.3) to[out=45,in=-165] (2,2);
		\end{tikzpicture}
		\caption{}
		\label{fig:nl l step 1 third b}
	\end{subfigure}
	\begin{subfigure}{0.24\textwidth}
		\centering
		\begin{tikzpicture}
			\draw[green] (0.15,0.9) -- (1.9,0.5);
			\draw[red] (2,0) -- (0,0.7) (2,0.7) -- (1,1);
			\draw[red] (2,0) to[out=105,in=-15] (0,2);

			\draw[blue] (0,0) -- (2,2);
			\draw[blue] (0,0) -- (2,0) -- (2,2) -- (0,2) -- (0,0);
			\draw[blue](1,1) -- (2,0);

			\draw[fill] (0,0) circle (0.05) 
			(2,2) circle (0.05)
			(0,2) circle (0.05)
			(2,0) circle (0.05)
			(1,1) circle (0.05);

			\draw (0,2) -- (1,1);
			\draw (0,0) to[out=75,in=-135] (0.7,1.3) to[out=45,in=-165] (2,2);
		\end{tikzpicture}
		\caption{}
		\label{fig:nl l step 1 third c}
	\end{subfigure}
	\begin{subfigure}{0.24\textwidth}
		\centering
		\begin{tikzpicture}
			\draw[green] (0,2) to[out=-75,in=180] (1,0.8) to[out=0,in=-105] (2,2);

			\draw[red] (2,0) -- (0,0.7) (2,0.7) -- (1,1);
			\draw[red] (2,0) to[out=105,in=-15] (0,2);

			\draw[blue] (0,0) -- (2,2);
			\draw[blue] (0,0) -- (2,0) -- (2,2) -- (0,2) -- (0,0);
			\draw[blue](1,1) -- (2,0);

			\draw[fill] (0,0) circle (0.05) 
			(2,2) circle (0.05)
			(0,2) circle (0.05)
			(2,0) circle (0.05)
			(1,1) circle (0.05);

			\draw (0,2) -- (1,1);
			\draw (0,0) to[out=75,in=-135] (0.7,1.3) to[out=45,in=-165] (2,2);
		\end{tikzpicture}
		\caption{}
		\label{fig:nl l step 1 third d}
	\end{subfigure}
	\caption{}
	\label{fig:nl l step 1 third}
\end{figure}

Similarly, there must be some arc $b' \in \mathcal{A}$ which intersecs $b$ but not $u$. 
Applying \Cref{lem:twos} to the disk shown in red in \Cref{fig:nl l step 1 third a}, we determine that either the path shown in green in \Cref{fig:nl l step 1 third b} or the path shown in green in \Cref{fig:nl l step 1 third c} must be a subarc of \( b' \). We find that the path shown in \Cref{fig:nl l step 1 third b} cannot be the subarc of an arc belonging to \( \mathcal{A} \), while the path shown in \Cref{fig:nl l step 1 third c} is uniquely a subarc of the arc shown in green in \Cref{fig:nl l step 1 third d}. 


There must also be an arc in $\mathcal{A}$ which intersects $c$ but not $u$, but $w_1 \in \mathcal{A}$ satisfies this.

\begin{figure}
	\centering
	\begin{subfigure}{0.24\textwidth}
		\centering
		\begin{tikzpicture}
			\fill[color=white!70!red] (0,2) to[out=-75,in=180] (1,0.8) -- (1.18,0.82) -- (2,0) -- (0,0.7);
			\fill[color=white!70!red] (2,0.7) -- (1.36,0.88) to[out=36,in=-105] (2,2);

			\draw[green] (0,2) to[out=-75,in=180] (1,0.8) to[out=0,in=-105] (2,2);
			\draw[green] (1,0) -- (1,2);

			\draw[red] (2,0) -- (0,0.7) (2,0.7) -- (1,1);
			\draw[red] (2,0) to[out=105,in=-15] (0,2);

			\draw[blue] (0,0) -- (2,2);
			\draw[blue] (0,0) -- (2,0) -- (2,2) -- (0,2) -- (0,0);
			\draw[blue](1,1) -- (2,0);

			\draw[fill] (0,0) circle (0.05) 
			(2,2) circle (0.05)
			(0,2) circle (0.05)
			(2,0) circle (0.05)
			(1,1) circle (0.05);

			\draw (0,2) -- (1,1);
			\draw (0,0) to[out=75,in=-135] (0.7,1.3) to[out=45,in=-165] (2,2);
		\end{tikzpicture}
		\caption{}
		\label{fig:nl l step 1 fourth a}
	\end{subfigure}
	\begin{subfigure}{0.24\textwidth}
		\centering
		\begin{tikzpicture}
			\draw[red] (2,0) to[out=155,in=-45] (0.7,0.7) to[out=135,in=-75] (0,2);
			\draw[green] (0,2) to[out=-65,in=180] (1,0.8) to[out=0,in=-105] (2,2);
			\draw[green] (1,0) -- (1,2);

			\draw[red] (2,0) -- (0,0.7) (2,0.7) -- (1,1);
			\draw[red] (2,0) to[out=105,in=-15] (0,2);

			\draw[blue] (0,0) -- (2,2);
			\draw[blue] (0,0) -- (2,0) -- (2,2) -- (0,2) -- (0,0);
			\draw[blue](1,1) -- (2,0);

			\draw[fill] (0,0) circle (0.05) 
			(2,2) circle (0.05)
			(0,2) circle (0.05)
			(2,0) circle (0.05)
			(1,1) circle (0.05);

			\draw (0,2) -- (1,1);
			\draw (0,0) to[out=75,in=-135] (0.7,1.3) to[out=45,in=-165] (2,2);
		\end{tikzpicture}
		\caption{}
		\label{fig:nl l step 1 fourth b}
	\end{subfigure}
	\caption{}
	\label{fig:nl l step 1 fourth}
\end{figure}

We have now determined 11 of the 12 arcs in $\mathcal{A}$. We may apply Corollary~\ref{cor:twoss} to the disk shown in red in \Cref{fig:nl l step 1 fourth a}, and we obtain the final arc which must be in $\mathcal{A}$. This is shown in Figure~\ref{fig:nl l step 1 fourth b}.

\begin{figure}
	\centering
	\begin{subfigure}{0.22\textwidth}
		\centering
		\begin{tikzpicture}
			\fill[color=white!70!red] (2,2) to[out=-165,in=45] (0.7,1.3) -- (0,2);
			\fill[color=white!70!red] (2,0) to[out=165,in=-52] (0.45,1.02) to[out=-128,in=75] (0,0);

			\draw[red] (2,0) to[out=105,in=-15] (0,2);
			\draw[red] (2,0) to[out=165,in=-75] (0,2);

			\draw[blue] (0,0) -- (2,2);
			\draw[dashed] (0,0) -- (2,0) -- (2,2) -- (0,2) -- (0,0);
			\draw[fill] (0,0) circle (0.05) 
			(2,2) circle (0.05)
			(0,2) circle (0.05)
			(2,0) circle (0.05)
			(1,1) circle (0.05);

			\draw (0,2) -- (1,1);
			\draw (0,0) to[out=75,in=-135] (0.7,1.3) to[out=45,in=-165] (2,2);
		\end{tikzpicture}
		\caption{}
		\label{fig:nl l step 2 a}
	\end{subfigure}
	\begin{subfigure}{0.22\textwidth}
		\centering
		\begin{tikzpicture}
			\fill[color=white!70!red] (2,0) to[out=105,in=-38] (0.98,1.55) to[out=-142,in=52] (0.45,1.02) to[out=-52,in=165] (2,0);

			\draw[red] (2,0) to[out=105,in=-15] (0,2);
			\draw[red] (2,0) to[out=165,in=-75] (0,2);

			\draw[blue] (0,0) -- (2,2);
			\draw[dashed] (0,0) -- (2,0) -- (2,2) -- (0,2) -- (0,0);
			\draw[fill] (0,0) circle (0.05) 
			(2,2) circle (0.05)
			(0,2) circle (0.05)
			(2,0) circle (0.05)
			(1,1) circle (0.05);

			\draw (0,2) -- (1,1);
			\draw (0,0) to[out=75,in=-135] (0.7,1.3) to[out=45,in=-165] (2,2);
		\end{tikzpicture}
		\caption{}
		\label{fig:nl l step 2 b}
	\end{subfigure}
	\begin{subfigure}{0.22\textwidth}
		\centering
		\begin{tikzpicture}
			\fill[color=white!70!red] (0,0) to[out=75,in=-135] (0.7,1.3) -- (0,2);
			\fill[color=white!70!red] (2,0) to[out=105,in=-38] (0.98,1.55) to[out=38,in=-165] (2,2);

			\draw[red] (2,0) to[out=105,in=-15] (0,2);
			\draw[red] (2,0) to[out=165,in=-75] (0,2);

			\draw[blue] (0,0) -- (2,2);
			\draw[dashed] (0,0) -- (2,0) -- (2,2) -- (0,2) -- (0,0);
			\draw[fill] (0,0) circle (0.05) 
			(2,2) circle (0.05)
			(0,2) circle (0.05)
			(2,0) circle (0.05)
			(1,1) circle (0.05);

			\draw (0,2) -- (1,1);
			\draw (0,0) to[out=75,in=-135] (0.7,1.3) to[out=45,in=-165] (2,2);
		\end{tikzpicture}
		\caption{}
		\label{fig:nl l step 2 c}
	\end{subfigure}
	\begin{subfigure}{0.3\textwidth}
		\centering
		\begin{tikzpicture}
			\draw[red] (2,0) to[out=105,in=-15] (0,2);
			\draw[red] (2,0) to[out=165,in=-75] (0,2);

			\draw[blue] (0,0) -- (2,2);
			\draw[blue] (0,0) -- (2,0) -- (2,2) -- (0,2) -- (0,0);
			\draw[blue] (2,0) -- (1,1);

			\draw[fill] (0,0) circle (0.05) 
			(2,2) circle (0.05)
			(0,2) circle (0.05)
			(2,0) circle (0.05)
			(1,1) circle (0.05);

			\draw (0,2) -- (1,1);
			\draw (0,0) to[out=75,in=-135] (0.7,1.3) to[out=45,in=-165] (2,2);

			\node[blue] at (1,0.1) {$a$};
			\draw[fill,white] (1.5,0.5) circle (0.15);
			\node[blue] at (1.5,0.5) {$b$};
			\draw[fill,white] (2,1.2) circle (0.15);
			\node[blue] at (2,1.2) {$c$};
		\end{tikzpicture}
		\caption{}
		\label{fig:nl l step 2 d}
	\end{subfigure}
	\caption{}
	\label{fig:nl l step 2}
\end{figure}

\textbf{Step 2.} In this step, we assume $w_2, x_2 \in \mathcal{A}$. This is shown in Figure~\ref{fig:nl l step 2}. We apply Corollary~\ref{cor:twoss} to three disks, shown in red in \Cref{fig:nl l step 2 a,fig:nl l step 2 b,fig:nl l step 2 c}. We conclude that three more arcs \( a, b, c \) must be included in $\mathcal{A}$, as shown in Figure~\ref{fig:nl l step 2 d}. 

By assumption, $u$ is minimally intersected, so there must be an arc $c' \in \mathcal{A}$ which intersects $c$ but not $u$. We apply Lemma~\ref{lem:twos} to the disk shown in red in Figure~\ref{fig:nl l step 2 c}, to determine a subarc of $c'$. This is shown in Figure~\ref{fig:nl l step 2 second a}. 

\begin{figure}
	\centering
	\begin{subfigure}{0.24\textwidth}
		\centering
		\begin{tikzpicture}
			\draw[green] (1.4,0.9) -- (2,0.8) (0,0.8) -- (0.5,0.7);
			\draw[red] (2,0) to[out=105,in=-15] (0,2);
			\draw[red] (2,0) to[out=165,in=-75] (0,2);

			\draw[blue] (0,0) -- (2,2);
			\draw[blue] (0,0) -- (2,0) -- (2,2) -- (0,2) -- (0,0);
			\draw[blue] (2,0) -- (1,1);

			\draw[fill] (0,0) circle (0.05) 
			(2,2) circle (0.05)
			(0,2) circle (0.05)
			(2,0) circle (0.05)
			(1,1) circle (0.05);

			\draw (0,2) -- (1,1);
			\draw (0,0) to[out=75,in=-135] (0.7,1.3) to[out=45,in=-165] (2,2);
		\end{tikzpicture}
		\caption{}
		\label{fig:nl l step 2 second a}
	\end{subfigure}
	\begin{subfigure}{0.24\textwidth}
		\centering
		\begin{tikzpicture}
			\draw[green] (1,1) -- (2,0.6) (0,0.6) -- (2,0);
			\draw[red] (2,0) to[out=105,in=-15] (0,2);
			\draw[red] (2,0) to[out=165,in=-75] (0,2);

			\draw[blue] (0,0) -- (2,2);
			\draw[blue] (0,0) -- (2,0) -- (2,2) -- (0,2) -- (0,0);
			\draw[blue] (2,0) -- (1,1);

			\draw[fill] (0,0) circle (0.05) 
			(2,2) circle (0.05)
			(0,2) circle (0.05)
			(2,0) circle (0.05)
			(1,1) circle (0.05);

			\draw (0,2) -- (1,1);
			\draw (0,0) to[out=75,in=-135] (0.7,1.3) to[out=45,in=-165] (2,2);
		\end{tikzpicture}
		\caption{}
		\label{fig:nl l step 2 second b}
	\end{subfigure}
	\begin{subfigure}{0.24\textwidth}
		\centering
		\begin{tikzpicture}
			\draw[green] (0,1) -- (2,1);
			\draw[red] (2,0) to[out=105,in=-15] (0,2);
			\draw[red] (2,0) to[out=165,in=-75] (0,2);

			\draw[blue] (0,0) -- (2,2);
			\draw[blue] (0,0) -- (2,0) -- (2,2) -- (0,2) -- (0,0);
			\draw[blue] (2,0) -- (1,1);

			\draw[fill] (0,0) circle (0.05) 
			(2,2) circle (0.05)
			(0,2) circle (0.05)
			(2,0) circle (0.05)
			(1,1) circle (0.05);

			\draw (0,2) -- (1,1);
			\draw (0,0) to[out=75,in=-135] (0.7,1.3) to[out=45,in=-165] (2,2);
		\end{tikzpicture}
		\caption{}
		\label{fig:nl l step 2 second c}
	\end{subfigure}
	\begin{subfigure}{0.24\textwidth}
		\centering
		\begin{tikzpicture}
			\draw[green] (0,1) -- (2,1)
			(1,0) -- (1,2);
			\draw[red] (2,0) to[out=105,in=-15] (0,2);
			\draw[red] (2,0) to[out=165,in=-75] (0,2);

			\draw[blue] (0,0) -- (2,2);
			\draw[blue] (0,0) -- (2,0) -- (2,2) -- (0,2) -- (0,0);
			\draw[blue] (2,0) -- (1,1);

			\draw[fill] (0,0) circle (0.05) 
			(2,2) circle (0.05)
			(0,2) circle (0.05)
			(2,0) circle (0.05)
			(1,1) circle (0.05);

			\draw (0,2) -- (1,1);
			\draw (0,0) to[out=75,in=-135] (0.7,1.3) to[out=45,in=-165] (2,2);
		\end{tikzpicture}
		\caption{}
		\label{fig:nl l step 2 second d}
	\end{subfigure}
	\caption{}
	\label{fig:nl l step 2 second}
\end{figure}

We see that \( c' \) must be one of two arcs, shown in green in \Cref{fig:nl l step 2 second b,fig:nl l step 2 second c}. Note that the arc shown in green \Cref{fig:nl l step 2 second b} is \( w_1' \), so up to the analysis in Step 1 we may assume that \( c' \) is the arc shown in green in \Cref{fig:nl l step 2 second c}. By a symmetric argument we see that the other arc shown in green in \Cref{fig:nl l step 2 second d} must be in \( \mathcal{A} \). 


\begin{figure}
	\centering
	\begin{subfigure}{0.4\textwidth}
		\centering
		\begin{tikzpicture}
			\draw[red] (0,1) -- (2,1)
			(1,0) -- (1,2);
			\draw[red] (2,0) to[out=105,in=-15] (0,2);
			\draw[red] (2,0) to[out=165,in=-75] (0,2);

			\draw[blue] (0,0) -- (2,2);
			\draw[blue] (0,0) -- (2,0) -- (2,2) -- (0,2) -- (0,0);
			\draw[blue] (2,0) -- (1,1);

			\draw[green] (0.5,1.3) to[out=-60,in=150] (1,0.7) to[out=-30,in=-150] (1.8,0.8);

			\draw[fill] (0,0) circle (0.05) 
			(2,2) circle (0.05)
			(0,2) circle (0.05)
			(2,0) circle (0.05)
			(1,1) circle (0.05);

			\draw (0,2) -- (1,1);
			\draw (0,0) to[out=75,in=-135] (0.7,1.3) to[out=45,in=-165] (2,2);
		\end{tikzpicture}
		\begin{tikzpicture}
			\draw[red] (0,1) -- (2,1)
			(1,0) -- (1,2);
			\draw[red] (2,0) to[out=105,in=-15] (0,2);
			\draw[red] (2,0) to[out=165,in=-75] (0,2);

			\draw[blue] (0,0) -- (2,2);
			\draw[blue] (0,0) -- (2,0) -- (2,2) -- (0,2) -- (0,0);
			\draw[blue] (2,0) -- (1,1);

			\draw[green] (1.2,0.2) -- (1.8,0.8);

			\draw[fill] (0,0) circle (0.05) 
			(2,2) circle (0.05)
			(0,2) circle (0.05)
			(2,0) circle (0.05)
			(1,1) circle (0.05);

			\draw (0,2) -- (1,1);
			\draw (0,0) to[out=75,in=-135] (0.7,1.3) to[out=45,in=-165] (2,2);
		\end{tikzpicture}
		\caption{}
		\label{fig:nl l step 2 third a}
	\end{subfigure}
	\begin{subfigure}{0.4\textwidth}
		\centering
		\begin{tikzpicture}
			\draw[red] (0,1) -- (2,1)
			(1,0) -- (1,2);
			\draw[red] (2,0) to[out=105,in=-15] (0,2);
			\draw[red] (2,0) to[out=165,in=-75] (0,2);

			\draw[blue] (0,0) -- (2,2);
			\draw[blue] (0,0) -- (2,0) -- (2,2) -- (0,2) -- (0,0);
			\draw[blue] (2,0) -- (1,1);

			\draw[green] (0,0) to[out=15,in=-105] (2,2);

			\draw[fill] (0,0) circle (0.05) 
			(2,2) circle (0.05)
			(0,2) circle (0.05)
			(2,0) circle (0.05)
			(1,1) circle (0.05);

			\draw (0,2) -- (1,1);
			\draw (0,0) to[out=75,in=-135] (0.7,1.3) to[out=45,in=-165] (2,2);
		\end{tikzpicture}
		\caption{}
		\label{fig:nl l step 2 third b}
	\end{subfigure}
	\caption{}
	\label{fig:nl l step 2 third}
\end{figure}

By assumption, \( u \) is minimally intersected, so there must be an arc \( b' \in \mathcal{A} \) which intersects \( b \) but not \( u \). We apply \Cref{lem:twos} to the disk shown in red in \Cref{fig:nl l step 2 b} and we conclude that, up to reflection, one of the two paths shown in green in \Cref{fig:nl l step 2 third a} must be a subarc of \( b' \). We find that in fact \( b' \) must be the arc shown in green in \Cref{fig:nl l step 2 third b}. This is a maximal 1-system. 


\textbf{Step 3.} In this step we assume \( w_3', x_3' \in \mathcal{A}\). This is shown in \Cref{fig:nl l step 4 first a}. We redraw this on the surface \( S \setminus \left\{ u, w, x \right\} \), as shown in \Cref{fig:nl l step 4 first b}. 

\begin{figure}
	\centering
	\begin{subfigure}[c]{0.23\textwidth}
		\centering
		\begin{tikzpicture}
			\draw[red] (0,0) to[out=25,in=-90] (1.2,1) to[out=90,in=-25] (0,2);
			\draw[red] (2,2) to[out=-115,in=0] (1,0.8) to[out=180,in=-65] (0,2);

			\draw[blue] (0,0) -- (2,2);
			\draw[dashed] (0,0) -- (2,0) -- (2,2) -- (0,2) -- (0,0);
			\draw[fill] (0,0) circle (0.05) 
			(2,2) circle (0.05)
			(0,2) circle (0.05)
			(2,0) circle (0.05)
			(1,1) circle (0.05);

			\draw (0,2) -- (1,1);
			\draw (0,0) to[out=75,in=-135] (0.7,1.3) to[out=45,in=-165] (2,2);

			\draw[fill,white] (1.7,1) circle (0.2)
			(1,0.3) circle (0.2);
			\node[red] at (1,0.3) {\( x'_3 \)};
			\node[red] at (1.7,1) {\( w'_3 \)};

			\draw[fill,white] (1.5,1.5) circle (0.2)
			(0.5,0.5) circle (0.2);
			\node[blue] at (0.5,0.5) {\( w \)};
			\node[blue] at (1.5,1.5) {\( x \)};
		\end{tikzpicture}
		\caption{}
		\label{fig:nl l step 4 first a}
	\end{subfigure}
	\begin{subfigure}[c]{0.23\textwidth}
		\centering
		\begin{tikzpicture}[scale=1.2]
			\fill[color=white!70!red] (-60:1) to[out=120,in=0] (180:1) -- (-120:.67);
			\draw[green] (-60:1) -- (180:1);

			\draw[blue] (0:1) -- (-60:1) -- (-120:1) (180:1) -- (120:1) -- (60:1);
			\draw[red] (150:0.87) -- (60:1) -- (-30:0.87)
			(-90:0.87) -- (180:1) -- (90:0.87);

			\draw (0:1) -- (60:1) (180:1) -- (-120:1);
			\draw (30:0.87) -- (-60:1) -- (-150:0.87);

			\foreach \x in {0,60,...,300} {
				\draw[fill] (\x:1) circle (0.05);
			}

			\node at (1,0.5) {\( u \)};
			\node at (-1,-0.5) {\( u \)};

			\node[blue] at (1,-0.5) {\( w \)};
			\node[blue] at (-1,0.5) {\( w \)};

			\node[blue] at (0,1.1) {\( x \)};
			\node[blue] at (0,-1.1) {\( x \)};
		\end{tikzpicture}
		\caption{}
		\label{fig:nl l step 4 first b}
	\end{subfigure}
	\begin{subfigure}[c]{0.23\textwidth}
		\centering
		\begin{tikzpicture}[scale=1.2]
			\fill[color=white!70!red] (-60:1) to[out=120,in=0] (180:1) -- (-150:0.87) -- (-120:.67) -- (-90:0.87);
			\draw[red] (-90:0.78) -- (-60:0.3)
			(-150:0.78) -- (180:0.3);

			\draw[blue] (0:1) -- (-60:1) -- (-120:1) (180:1) -- (120:1) -- (60:1);
			\draw[red] (150:0.87) -- (60:1) -- (-30:0.87)
			(-90:0.87) -- (180:1) -- (90:0.87);

			\draw (0:1) -- (60:1) (180:1) -- (-120:1);
			\draw (30:0.87) -- (-60:1) -- (-150:0.87);

			\foreach \x in {0,60,...,300} {
				\draw[fill] (\x:1) circle (0.05);
			}
		\end{tikzpicture}
		\caption{}
		\label{fig:nl l step 4 first c}
	\end{subfigure}
	\begin{subfigure}[c]{0.23\textwidth}
		\centering
		\begin{tikzpicture}[scale=1.2]
			\fill[color=white!70!red] (-60:1) to[out=120,in=0] (180:1) -- (-150:0.87) -- (-120:.67) -- (-90:0.87);
			\fill[color=white!70!red] (60:1) -- (120:0.67) -- (90:0.87);
			\fill[color=white!70!red] (60:1) -- (0:0.67) -- (30:0.87);
			\draw[red] (80:0.82) -- (80:0.88) (-80:0.88) -- (-60:0.3)
			(40:0.82) -- (40:0.88) (-160:0.88) -- (180:0.3);

			\draw[blue] (0:1) -- (-60:1) -- (-120:1) (180:1) -- (120:1) -- (60:1);
			\draw[red] (150:0.87) -- (60:1) -- (-30:0.87)
			(-90:0.87) -- (180:1) -- (90:0.87);

			\draw (0:1) -- (60:1) (180:1) -- (-120:1);
			\draw (30:0.87) -- (-60:1) -- (-150:0.87);

			\foreach \x in {0,60,...,300} {
				\draw[fill] (\x:1) circle (0.05);
			}
		\end{tikzpicture}
		\caption{}
		\label{fig:nl l step 4 first d}
	\end{subfigure}
	\caption{}
	\label{fig:nl l step 4 first}
\end{figure}

Let \( a \) be the arc shown in green in \Cref{fig:nl l step 4 first b}. 

\begin{claim}
	\( a \in \mathcal{A} \). 
\end{claim}

\begin{proof}
	Suppose for contradiction that \( a \notin \mathcal{A} \). Since \( \mathcal{A} \) is saturated, there must be some \( a' \in \mathcal{A} \) which intersects \( a \) at least twice. We apply \Cref{lem:twos} to the disk shown in red in \Cref{fig:nl l step 4 first b} and we determine two subarcs of \( a' \). These are shown in red in \Cref{fig:nl l step 4 first c}. 

	Now we apply \Cref{lem:twos} to the disk shown in red in \Cref{fig:nl l step 4 first c}, and we determine longer paths which are subarcs of \( a' \). Finally, we apply \Cref{lem:twos} to the disk shown in red in \Cref{fig:nl l step 4 first d}, and we see that in fact \( a' \) must intersect some arc of \( \mathcal{A} \) at least twice. This is a contradiction. So, it must be that \( a \in \mathcal{A} \).
\end{proof}

We invoke the symmetry of the surface, and we conclude that two more arcs \( b,c \) must be in \( \mathcal{A} \), as shown in \Cref{fig:nl l step 4 second a}.

\begin{figure}
	\centering
	\begin{subfigure}{0.23\textwidth}
		\centering
		\begin{tikzpicture}[scale=1.2]
			\draw[blue] (-60:1) -- (180:1) -- (60:1) -- (-60:1);

			\draw[blue] (0:1) -- (-60:1) -- (-120:1) (180:1) -- (120:1) -- (60:1);
			\draw[red] (150:0.87) -- (60:1) -- (-30:0.87)
			(-90:0.87) -- (180:1) -- (90:0.87);

			\draw (0:1) -- (60:1) (180:1) -- (-120:1);
			\draw (30:0.87) -- (-60:1) -- (-150:0.87);

			\foreach \x in {0,60,...,300} {
				\draw[fill] (\x:1) circle (0.05);
			}

			\draw[fill,white] (0:0.5) circle (0.12)
			(120:0.5) circle (0.12)
			(-120:0.5) circle (0.12);

			\node[blue] at (0:0.5) {$c$};
			\node[blue] at (120:0.5) {$b$};
			\node[blue] at (-120:0.5) {$a$};
		\end{tikzpicture}
		\caption{}
		\label{fig:nl l step 4 second a}
	\end{subfigure}
	\begin{subfigure}{0.23\textwidth}
		\centering
		\begin{tikzpicture}[scale=1.2]
			\fill[color=white!70!red] (-60:1) to[out=120,in=0] (180:1) -- (-150:0.87) -- (-120:0.67)-- (-90:.87);
			\fill[color=white!70!red] (60:1) -- (90:0.87) -- (120:.67);
			\draw[blue] (-60:1) -- (180:1);

			\draw[blue] (0:1) -- (-60:1) -- (-120:1) (180:1) -- (120:1) -- (60:1);
			\draw[red] (150:0.87) -- (60:1) -- (-30:0.87)
			(-90:0.87) -- (180:1) -- (90:0.87);

			\draw (0:1) -- (60:1) (180:1) -- (-120:1);
			\draw (30:0.87) -- (-60:1) -- (-150:0.87);

			\foreach \x in {0,60,...,300} {
				\draw[fill] (\x:1) circle (0.05);
			}
		\end{tikzpicture}
		\caption{}
		\label{fig:nl l step 4 second b}
	\end{subfigure}
	\caption{}
	\label{fig:nl l step 4 second part 1}
\end{figure}
\begin{figure}
	\centering
	\begin{subfigure}{0.23\textwidth}
		\centering
		\begin{tikzpicture}[scale=1.2]
			\draw[green] (0,-0.3) -- (-80:0.88) (80:0.88) -- (75:0.6);

			\draw[blue] (-60:1) -- (180:1);

			\draw[blue] (0:1) -- (-60:1) -- (-120:1) (180:1) -- (120:1) -- (60:1);
			\draw[red] (150:0.87) -- (60:1) -- (-30:0.87)
			(-90:0.87) -- (180:1) -- (90:0.87);

			\draw (0:1) -- (60:1) (180:1) -- (-120:1);
			\draw (30:0.87) -- (-60:1) -- (-150:0.87);

			\foreach \x in {0,60,...,300} {
				\draw[fill] (\x:1) circle (0.05);
			}
		\end{tikzpicture}
		\caption{}
		\label{fig:nl l step 4 second c}
	\end{subfigure}
	\begin{subfigure}{0.23\textwidth}
		\centering
		\begin{tikzpicture}[scale=1.2]
			\draw[green] (0,-0.3) -- (-80:0.88) (80:0.88) -- (-60:1);

			\draw[blue] (-60:1) -- (180:1);

			\draw[blue] (0:1) -- (-60:1) -- (-120:1) (180:1) -- (120:1) -- (60:1);
			\draw[red] (150:0.87) -- (60:1) -- (-30:0.87)
			(-90:0.87) -- (180:1) -- (90:0.87);

			\draw (0:1) -- (60:1) (180:1) -- (-120:1);
			\draw (30:0.87) -- (-60:1) -- (-150:0.87);

			\foreach \x in {0,60,...,300} {
				\draw[fill] (\x:1) circle (0.05);
			}
		\end{tikzpicture}
		\caption{}
		\label{fig:nl l step 4 second d}
	\end{subfigure}
	\begin{subfigure}{0.23\textwidth}
		\centering
		\begin{tikzpicture}[scale=1.2]
			\draw[green] (60:1) -- (-80:0.88) (80:0.88) -- (180:1);

			\draw[blue] (-60:1) -- (180:1);

			\draw[blue] (0:1) -- (-60:1) -- (-120:1) (180:1) -- (120:1) -- (60:1);
			\draw[red] (150:0.87) -- (60:1) -- (-30:0.87)
			(-90:0.87) -- (180:1) -- (90:0.87);

			\draw (0:1) -- (60:1) (180:1) -- (-120:1);
			\draw (30:0.87) -- (-60:1) -- (-150:0.87);

			\foreach \x in {0,60,...,300} {
				\draw[fill] (\x:1) circle (0.05);
			}
		\end{tikzpicture}
		\caption{}
		\label{fig:nl l step 4 second e}
	\end{subfigure}
	\caption{}
	\label{fig:nl l step 4 second part 2}
\end{figure}
\begin{figure}
	\centering
	\begin{subfigure}{0.23\textwidth}
		\centering
		\begin{tikzpicture}[scale=1.2]
			\draw[green] (0,0) -- (-120:0.8);

			\draw[blue] (-60:1) -- (180:1);

			\draw[blue] (0:1) -- (-60:1) -- (-120:1) (180:1) -- (120:1) -- (60:1);
			\draw[red] (150:0.87) -- (60:1) -- (-30:0.87)
			(-90:0.87) -- (180:1) -- (90:0.87);

			\draw (0:1) -- (60:1) (180:1) -- (-120:1);
			\draw (30:0.87) -- (-60:1) -- (-150:0.87);

			\foreach \x in {0,60,...,300} {
				\draw[fill] (\x:1) circle (0.05);
			}
		\end{tikzpicture}
		\caption{}
		\label{fig:nl l step 4 second f}
	\end{subfigure}
	\begin{subfigure}{0.23\textwidth}
		\centering
		\begin{tikzpicture}[scale=1.2]
			\draw[green] (60:1) -- (-100:0.88)
			(100:0.88) -- (120:0.86);

			\draw[blue] (-60:1) -- (180:1);

			\draw[blue] (0:1) -- (-60:1) -- (-120:1) (180:1) -- (120:1) -- (60:1);
			\draw[red] (150:0.87) -- (60:1) -- (-30:0.87)
			(-90:0.87) -- (180:1) -- (90:0.87);

			\draw (0:1) -- (60:1) (180:1) -- (-120:1);
			\draw (30:0.87) -- (-60:1) -- (-150:0.87);

			\foreach \x in {0,60,...,300} {
				\draw[fill] (\x:1) circle (0.05);
			}
		\end{tikzpicture}
		\caption{}
		\label{fig:nl l step 4 second h}
	\end{subfigure}
	\begin{subfigure}{0.23\textwidth}
		\centering
		\begin{tikzpicture}[scale=1.2]
			\draw[green] (-120:1) -- (80:0.88);
			\draw[green] (-80:0.88) -- (-85:0.82);

			\draw[blue] (-60:1) -- (180:1);

			\draw[blue] (0:1) -- (-60:1) -- (-120:1) (180:1) -- (120:1) -- (60:1);
			\draw[red] (150:0.87) -- (60:1) -- (-30:0.87)
			(-90:0.87) -- (180:1) -- (90:0.87);

			\draw (0:1) -- (60:1) (180:1) -- (-120:1);
			\draw (30:0.87) -- (-60:1) -- (-150:0.87);

			\foreach \x in {0,60,...,300} {
				\draw[fill] (\x:1) circle (0.05);
			}
		\end{tikzpicture}
		\caption{}
		\label{fig:nl l step 4 second g}
	\end{subfigure}
	\begin{subfigure}{0.23\textwidth}
		\centering
		\begin{tikzpicture}[scale=1.2]
			\draw[green] (60:1) -- (-120:1);

			\draw[blue] (-60:1) -- (180:1);

			\draw[blue] (0:1) -- (-60:1) -- (-120:1) (180:1) -- (120:1) -- (60:1);
			\draw[red] (150:0.87) -- (60:1) -- (-30:0.87)
			(-90:0.87) -- (180:1) -- (90:0.87);

			\draw (0:1) -- (60:1) (180:1) -- (-120:1);
			\draw (30:0.87) -- (-60:1) -- (-150:0.87);

			\foreach \x in {0,60,...,300} {
				\draw[fill] (\x:1) circle (0.05);
			}
		\end{tikzpicture}
		\caption{}
		\label{fig:nl l step 4 second i}
	\end{subfigure}
	\caption{}
	\label{fig:nl l step 4 second}
\end{figure}

By assumption, \( u \) is minimally intersected, so there must be some arc \( a' \in \mathcal{A} \) which intersects \( a \) but not \( u \). We apply \Cref{lem:twos} to the disk shown in red in \Cref{fig:nl l step 4 second b}, and we see that either the path shown in green in \Cref{fig:nl l step 4 second c} or the path shown in green in \Cref{fig:nl l step 4 second f} must be a subarc of \( a' \). 

If the path shown in \Cref{fig:nl l step 4 second c} is a subarc, we consider how this may be completed to an arc, by the behaviour of the subarc at the top of the subfigure. If it is completed as in \Cref{fig:nl l step 4 second d}, the other end of the arc cannot be completed. Otherwise, it could be the arc shown in \Cref{fig:nl l step 4 second e}. Call this arc \( a'_1 \). 

If instead the path shown in \Cref{fig:nl l step 4 second f} is a subarc of \( a' \), consider the following. Becuase arcs cannot intersect twice, and we assume that \( a' \) does not intersect \( u \), we see that \( a' \) must contain one of the paths shown in green in \Cref{fig:nl l step 4 second g,fig:nl l step 4 second h,fig:nl l step 4 second i}. However, we see that the paths in \Cref{fig:nl l step 4 second g,fig:nl l step 4 second h} cannot be completed to arcs. So, let \( a'_2 \) be the arc shown in \Cref{fig:nl l step 4 second i}. 

We claim that \( a'_2 \in \mathcal{A} \). Indeed, suppose instead that \( a'_2 \notin \mathcal{A} \). Then we have \( a'_1 \in \mathcal{A} \), and there must be some arc \( \hat{a} \in \mathcal{A} \) which intersects \( a'_2 \) at least twice. We apply \Cref{lem:twos} to the disk shown in red in \Cref{fig:nl l step 4 second j} to determine two paths which must be subarcs of \( \hat{a} \). This is shown in green in \Cref{fig:nl l step 4 second j}. 

\begin{figure}
	\centering
	\begin{subfigure}{0.23\textwidth}
		\centering
		\begin{tikzpicture}[scale=1.2]
			\fill[color=white!70!red] (-120:1) to[out=80,in=-140] (60:1) -- (-80:0.88);

			\draw[green] (120:0.4) -- (-30:0.5) 
			(-132:0.8) -- (-100:0.88) (100:0.88) -- (120:0.8);
			\draw[red] (60:1) -- (-80:0.88) (80:0.88) -- (180:1);

			\draw[blue] (-60:1) -- (180:1);

			\draw[blue] (0:1) -- (-60:1) -- (-120:1) (180:1) -- (120:1) -- (60:1);
			\draw[red] (150:0.87) -- (60:1) -- (-30:0.87)
			(-90:0.87) -- (180:1) -- (90:0.87);

			\draw (0:1) -- (60:1) (180:1) -- (-120:1);
			\draw (30:0.87) -- (-60:1) -- (-150:0.87);

			\foreach \x in {0,60,...,300} {
				\draw[fill] (\x:1) circle (0.05);
			}
		\end{tikzpicture}
		\caption{}
		\label{fig:nl l step 4 second j}
	\end{subfigure}
	\begin{subfigure}{0.23\textwidth}
		\centering
		\begin{tikzpicture}[scale=1.2]
			\draw[green] (120:0.2) -- (-30:0.5) (-120:0.8) -- (-100:0.88) (100:0.88) -- (180:1);
			\draw[red] (60:1) -- (-80:0.88) (80:0.88) -- (180:1);

			\draw[blue] (-60:1) -- (180:1);

			\draw[blue] (0:1) -- (-60:1) -- (-120:1) (180:1) -- (120:1) -- (60:1);
			\draw[red] (150:0.87) -- (60:1) -- (-30:0.87)
			(-90:0.87) -- (180:1) -- (90:0.87);

			\draw (0:1) -- (60:1) (180:1) -- (-120:1);
			\draw (30:0.87) -- (-60:1) -- (-150:0.87);

			\foreach \x in {0,60,...,300} {
				\draw[fill] (\x:1) circle (0.05);
			}
		\end{tikzpicture}
		\caption{}
		\label{fig:nl l step 4 second k}
	\end{subfigure}
	\begin{subfigure}{0.23\textwidth}
		\centering
		\begin{tikzpicture}[scale=1.2]
			\draw[green] (120:0.2) -- (-30:0.5) 
			(180:1) -- (-100:0.88) (100:0.88) -- (140:0.88) (-20:0.88) -- (0:0.8);

			\draw[red] (60:1) -- (-80:0.88) (80:0.88) -- (180:1);

			\draw[blue] (-60:1) -- (180:1);

			\draw[blue] (0:1) -- (-60:1) -- (-120:1) (180:1) -- (120:1) -- (60:1);
			\draw[red] (150:0.87) -- (60:1) -- (-30:0.87)
			(-90:0.87) -- (180:1) -- (90:0.87);

			\draw (0:1) -- (60:1) (180:1) -- (-120:1);
			\draw (30:0.87) -- (-60:1) -- (-150:0.87);

			\foreach \x in {0,60,...,300} {
				\draw[fill] (\x:1) circle (0.05);
			}
		\end{tikzpicture}
		\caption{}
		\label{fig:nl l step 4 second l}
	\end{subfigure}
	\begin{subfigure}{0.23\textwidth}
		\centering
		\begin{tikzpicture}[scale=1.2]
			\draw[green] (120:0.2) -- (-30:0.5) 
			(-120:0.84) -- (-100:0.88) (100:0.88) -- (140:0.88) (-20:0.88) -- (60:1);

			\draw[red] (60:1) -- (-80:0.88) (80:0.88) -- (180:1);

			\draw[blue] (-60:1) -- (180:1);

			\draw[blue] (0:1) -- (-60:1) -- (-120:1) (180:1) -- (120:1) -- (60:1);
			\draw[red] (150:0.87) -- (60:1) -- (-30:0.87)
			(-90:0.87) -- (180:1) -- (90:0.87);

			\draw (0:1) -- (60:1) (180:1) -- (-120:1);
			\draw (30:0.87) -- (-60:1) -- (-150:0.87);

			\foreach \x in {0,60,...,300} {
				\draw[fill] (\x:1) circle (0.05);
			}
		\end{tikzpicture}
		\caption{}
		\label{fig:nl l step 4 second m}
	\end{subfigure}
	\caption{}
	\label{fig:nl l step 4 second part 4}
\end{figure}

We see that in fact there is no way to extend these paths to make an arc. See \Cref{fig:nl l step 4 second part 4}. So, there is no arc that intersects \( a'_2 \) twice which could be in \( \mathcal{A} \). This contradicts the assumption that \( \mathcal{A} \) is maximal. So, it must be that \( a'_2 \in \mathcal{A} \). 

\begin{figure}
	\centering
	\begin{subfigure}{0.23\textwidth}
		\centering
		\begin{tikzpicture}[scale=1.2]
			\fill[color=white!70!red] (90:0.87) -- (120:0.67) -- (150:0.87) -- (-1,0) to[out=10,in=-130] (60:1);

			\draw[green] (60:1) -- (-120:1)
			(0:1) -- (180:1);

			\draw[blue] (-60:1) -- (180:1) -- (60:1) -- (-60:1);

			\draw[blue] (0:1) -- (-60:1) -- (-120:1) (180:1) -- (120:1) -- (60:1);
			\draw[red] (150:0.87) -- (60:1) -- (-30:0.87)
			(-90:0.87) -- (180:1) -- (90:0.87);

			\draw (0:1) -- (60:1) (180:1) -- (-120:1);
			\draw (30:0.87) -- (-60:1) -- (-150:0.87);

			\foreach \x in {0,60,...,300} {
				\draw[fill] (\x:1) circle (0.05);
			}
		\end{tikzpicture}
		\caption{}
		\label{fig:nl l step 4 third a}
	\end{subfigure}
	\begin{subfigure}{0.23\textwidth}
		\centering
		\begin{tikzpicture}[scale=1.2]
			\draw[green] (120:0.3) -- (160:0.88) (-40:0.88) -- (-40:0.82);

			\draw[red] (60:1) -- (-120:1)
			(0:1) -- (180:1);

			\draw[blue] (-60:1) -- (180:1) -- (60:1) -- (-60:1);

			\draw[blue] (0:1) -- (-60:1) -- (-120:1) (180:1) -- (120:1) -- (60:1);
			\draw[red] (150:0.87) -- (60:1) -- (-30:0.87)
			(-90:0.87) -- (180:1) -- (90:0.87);

			\draw (0:1) -- (60:1) (180:1) -- (-120:1);
			\draw (30:0.87) -- (-60:1) -- (-150:0.87);

			\foreach \x in {0,60,...,300} {
				\draw[fill] (\x:1) circle (0.05);
			}
		\end{tikzpicture}
		\caption{}
		\label{fig:nl l step 4 third b}
	\end{subfigure}
	\begin{subfigure}{0.23\textwidth}
		\centering
		\begin{tikzpicture}[scale=1.2]
			\draw[green] (30:0.3) -- (160:0.88) (-40:0.88) -- (-40:0.82);

			\draw[red] (60:1) -- (-120:1)
			(0:1) -- (180:1);

			\draw[blue] (-60:1) -- (180:1) -- (60:1) -- (-60:1);

			\draw[blue] (0:1) -- (-60:1) -- (-120:1) (180:1) -- (120:1) -- (60:1);
			\draw[red] (150:0.87) -- (60:1) -- (-30:0.87)
			(-90:0.87) -- (180:1) -- (90:0.87);

			\draw (0:1) -- (60:1) (180:1) -- (-120:1);
			\draw (30:0.87) -- (-60:1) -- (-150:0.87);

			\foreach \x in {0,60,...,300} {
				\draw[fill] (\x:1) circle (0.05);
			}
		\end{tikzpicture}
		\caption{}
		\label{fig:nl l step 4 third b5}
	\end{subfigure}
	\begin{subfigure}{0.23\textwidth}
		\centering
		\begin{tikzpicture}[scale=1.2]
			\draw[green] (30:0.3) -- (160:0.88)
			(-40:0.88) -- (60:1);

			\draw[red] (60:1) -- (-120:1)
			(0:1) -- (180:1);

			\draw[blue] (-60:1) -- (180:1) -- (60:1) -- (-60:1);

			\draw[blue] (0:1) -- (-60:1) -- (-120:1) (180:1) -- (120:1) -- (60:1);
			\draw[red] (150:0.87) -- (60:1) -- (-30:0.87)
			(-90:0.87) -- (180:1) -- (90:0.87);

			\draw (0:1) -- (60:1) (180:1) -- (-120:1);
			\draw (30:0.87) -- (-60:1) -- (-150:0.87);

			\foreach \x in {0,60,...,300} {
				\draw[fill] (\x:1) circle (0.05);
			}
		\end{tikzpicture}
		\caption{}
		\label{fig:nl l step 4 third c}
	\end{subfigure}
	\begin{subfigure}{0.23\textwidth}
		\centering
		\begin{tikzpicture}[scale=1.2]
			\draw[green] (120:0.3) -- (120:0.9);

			\draw[red] (60:1) -- (-120:1)
			(0:1) -- (180:1);

			\draw[blue] (-60:1) -- (180:1) -- (60:1) -- (-60:1);

			\draw[blue] (0:1) -- (-60:1) -- (-120:1) (180:1) -- (120:1) -- (60:1);
			\draw[red] (150:0.87) -- (60:1) -- (-30:0.87)
			(-90:0.87) -- (180:1) -- (90:0.87);

			\draw (0:1) -- (60:1) (180:1) -- (-120:1);
			\draw (30:0.87) -- (-60:1) -- (-150:0.87);

			\foreach \x in {0,60,...,300} {
				\draw[fill] (\x:1) circle (0.05);
			}
		\end{tikzpicture}
		\caption{}
		\label{fig:nl l step 4 third d}
	\end{subfigure}
	\begin{subfigure}{0.23\textwidth}
		\centering
		\begin{tikzpicture}[scale=1.2]
			\draw[green] (-60:0.4) -- (120:1);

			\draw[red] (60:1) -- (-120:1)
			(0:1) -- (180:1);

			\draw[blue] (-60:1) -- (180:1) -- (60:1) -- (-60:1);

			\draw[blue] (0:1) -- (-60:1) -- (-120:1) (180:1) -- (120:1) -- (60:1);
			\draw[red] (150:0.87) -- (60:1) -- (-30:0.87)
			(-90:0.87) -- (180:1) -- (90:0.87);

			\draw (0:1) -- (60:1) (180:1) -- (-120:1);
			\draw (30:0.87) -- (-60:1) -- (-150:0.87);

			\foreach \x in {0,60,...,300} {
				\draw[fill] (\x:1) circle (0.05);
			}
		\end{tikzpicture}
		\caption{}
		\label{fig:nl l step 4 third e}
	\end{subfigure}
	\begin{subfigure}{0.23\textwidth}
		\centering
		\begin{tikzpicture}[scale=1.2]
			\draw[green] (80:0.83) -- (80:0.88) (-80:0.88) -- (120:1);

			\draw[red] (60:1) -- (-120:1)
			(0:1) -- (180:1);

			\draw[blue] (-60:1) -- (180:1) -- (60:1) -- (-60:1);

			\draw[blue] (0:1) -- (-60:1) -- (-120:1) (180:1) -- (120:1) -- (60:1);
			\draw[red] (150:0.87) -- (60:1) -- (-30:0.87)
			(-90:0.87) -- (180:1) -- (90:0.87);

			\draw (0:1) -- (60:1) (180:1) -- (-120:1);
			\draw (30:0.87) -- (-60:1) -- (-150:0.87);

			\foreach \x in {0,60,...,300} {
				\draw[fill] (\x:1) circle (0.05);
			}
		\end{tikzpicture}
		\caption{}
		\label{fig:nl l step 4 third f}
	\end{subfigure}
	\begin{subfigure}{0.23\textwidth}
		\centering
		\begin{tikzpicture}[scale=1.2]
			\draw[green] (-60:1) -- (120:1);

			\draw[red] (60:1) -- (-120:1)
			(0:1) -- (180:1);

			\draw[blue] (-60:1) -- (180:1) -- (60:1) -- (-60:1);

			\draw[blue] (0:1) -- (-60:1) -- (-120:1) (180:1) -- (120:1) -- (60:1);
			\draw[red] (150:0.87) -- (60:1) -- (-30:0.87)
			(-90:0.87) -- (180:1) -- (90:0.87);

			\draw (0:1) -- (60:1) (180:1) -- (-120:1);
			\draw (30:0.87) -- (-60:1) -- (-150:0.87);

			\foreach \x in {0,60,...,300} {
				\draw[fill] (\x:1) circle (0.05);
			}
		\end{tikzpicture}
		\caption{}
		\label{fig:nl l step 4 third g}
	\end{subfigure}
	\caption{Constructing a 1-system}
	\label{fig:nl l step 4 third}
\end{figure}
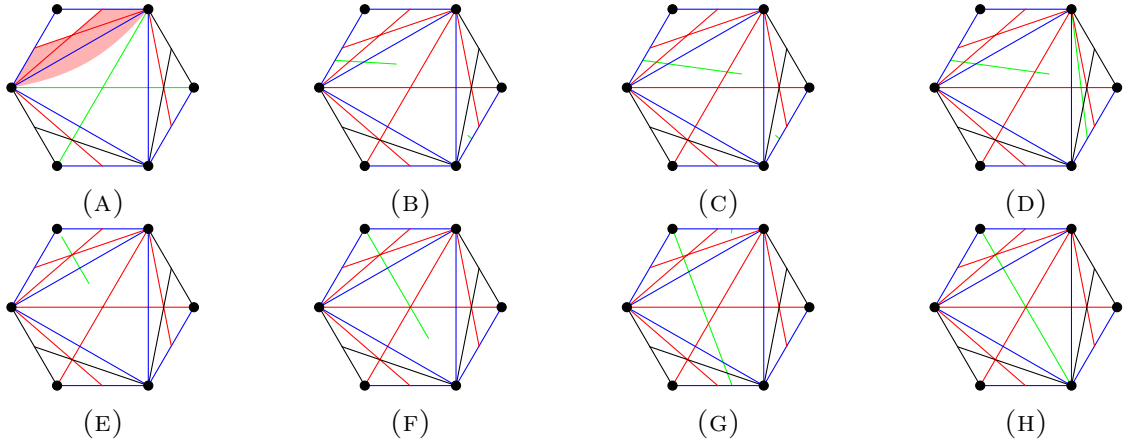

We have \( a' = a'_2 \in \mathcal{A} \). By symmetry, we find that another arc \( c' \) must be in \( \mathcal{A} \), as shown in green in \Cref{fig:nl l step 4 third a}. 

Finally, there must be some arc \( b' \in \mathcal{A} \) which intersects \( b \) but not \( u \). We apply \Cref{lem:twos} to the disk shown in red in \Cref{fig:nl l step 4 third a}, and up to symmetry we see that either the path shown in green in \Cref{fig:nl l step 4 third b} or in \Cref{fig:nl l step 4 third d} must be a subarc of \( b' \). 

First, if \( b' \) contains the path shown in \Cref{fig:nl l step 4 third b}, we try to extend this path to an arc. We uniquely extend one end as shown in \Cref{fig:nl l step 4 third b5}, then we uniquely extend the other end as shown in \Cref{fig:nl l step 4 third c}. We see that this path cannot be extended to an arc without intersecting some arc of \( \mathcal{A} \) twice.

So, it must be that \( b' \) contains the path shown in \Cref{fig:nl l step 4 third d}. We extend this path as shown in \Cref{fig:nl l step 4 third e}. Then, it cannot be that \( b' \) contains the path shown in \Cref{fig:nl l step 4 third f}, so it must be that \( b' \) is the arc shown in \Cref{fig:nl l step 4 third g}. 

Thus, we obtain a maximal 1-system.

\textbf{Step 4.} In this step, we assume \( w_2', x_3' \in \mathcal{A} \). This is shown in \Cref{fig:nl l step 3 a}. Let \( b \) be the arc shown in green in \Cref{fig:nl l step 3 b}. 

\begin{figure}
	\centering
	\begin{subfigure}{0.23\textwidth}
		\centering
		\begin{tikzpicture}
			\draw[red] (2,0) to[out=165,in=-75] (0,2);
			\draw[red] (0,0) to[out=25,in=-90] (1.2,1) to[out=90,in=-25] (0,2);

			\draw[blue] (0,0) -- (2,2);
			\draw[dashed] (0,0) -- (2,0) -- (2,2) -- (0,2) -- (0,0);
			\draw[fill] (0,0) circle (0.05) 
			(2,2) circle (0.05)
			(0,2) circle (0.05)
			(2,0) circle (0.05)
			(1,1) circle (0.05);

			\draw (0,2) -- (1,1);
			\draw (0,0) to[out=75,in=-135] (0.7,1.3) to[out=45,in=-165] (2,2);
		\end{tikzpicture}
		\caption{}
		\label{fig:nl l step 3 a}
	\end{subfigure}
	\begin{subfigure}{0.23\textwidth}
		\centering
		\begin{tikzpicture}
			\draw[green] (2,0) -- (1,1);

			\draw[red] (2,0) to[out=165,in=-75] (0,2);
			\draw[red] (0,0) to[out=25,in=-90] (1.2,1) to[out=90,in=-25] (0,2);

			\draw[blue] (0,0) -- (2,2);
			\draw[dashed] (0,0) -- (2,0) -- (2,2) -- (0,2) -- (0,0);
			\draw[fill] (0,0) circle (0.05) 
			(2,2) circle (0.05)
			(0,2) circle (0.05)
			(2,0) circle (0.05)
			(1,1) circle (0.05);

			\draw (0,2) -- (1,1);
			\draw (0,0) to[out=75,in=-135] (0.7,1.3) to[out=45,in=-165] (2,2);
		\end{tikzpicture}
		\caption{}
		\label{fig:nl l step 3 b}
	\end{subfigure}
	\begin{subfigure}{0.23\textwidth}
		\centering
		\begin{tikzpicture}
			\draw[green] (0.7,0) to[out=45,in=-110] (2,2) to[out=-125,in=50] (1.3,1) to[out=-130,in=-60] (0.7,1) to[out=120,in=-110] (0.7,2);

			\draw[red] (2,0) to[out=165,in=-75] (0,2);
			\draw[red] (0,0) to[out=25,in=-90] (1.2,1) to[out=90,in=-25] (0,2);

			\draw[blue] (0,0) -- (2,2);
			\draw[dashed] (0,0) -- (2,0) -- (2,2) -- (0,2) -- (0,0);
			\draw[fill] (0,0) circle (0.05) 
			(2,2) circle (0.05)
			(0,2) circle (0.05)
			(2,0) circle (0.05)
			(1,1) circle (0.05);

			\draw (0,2) -- (1,1);
			\draw (0,0) to[out=75,in=-135] (0.7,1.3) to[out=45,in=-165] (2,2);
		\end{tikzpicture}
		\caption{}
		\label{fig:nl l step 3 c}
	\end{subfigure}
	\begin{subfigure}{0.23\textwidth}
		\centering
		\begin{tikzpicture}
			\draw[green] (0,2) to[out=-65,in=240] (1.25,0.9) to[out=60,in=-70] (1.3,2) (1.3,0) -- (2,2); 

			\draw[red] (2,0) to[out=165,in=-75] (0,2);
			\draw[red] (0,0) to[out=25,in=-90] (1.2,1) to[out=90,in=-25] (0,2);

			\draw[blue] (0,0) -- (2,2);
			\draw[dashed] (0,0) -- (2,0) -- (2,2) -- (0,2) -- (0,0);
			\draw[fill] (0,0) circle (0.05) 
			(2,2) circle (0.05)
			(0,2) circle (0.05)
			(2,0) circle (0.05)
			(1,1) circle (0.05);

			\draw (0,2) -- (1,1);
			\draw (0,0) to[out=75,in=-135] (0.7,1.3) to[out=45,in=-165] (2,2);
		\end{tikzpicture}
		\caption{}
		\label{fig:nl l step 3 d}
	\end{subfigure}
	\caption{}
	\label{fig:nl l step 3}
\end{figure}

\begin{claim}
	\( b \in \mathcal{A} \).
\end{claim}

\begin{proof}
	From the proof of \Cref{claim:sec5 some arcs}, we see that the only two arcs which could be in \( \mathcal{A} \) and which intersect \( b \) twice are the arcs shown in green in \Cref{fig:nl l step 3 c} and \Cref{fig:nl l step 3 d}. However, the arc shown in \Cref{fig:nl l step 3 c} intersects \( x'_3 \) twice, and the arc shown in \Cref{fig:nl l step 3 d} intersects \( v \) twice. Neither of these arcs are in \( \mathcal{A} \), and \( \mathcal{A} \) is saturated, so it must be that \( b \in \mathcal{A} \). 
\end{proof}

We apply \Cref{cor:twoss} to the disk shown in red in \Cref{fig:nl l step 3 second a} and we obtain an arc \( a \) which must be in \( \mathcal{A} \), as shown in green in \Cref{fig:nl l step 3 second b}. 

\begin{figure}
	\centering
	\begin{subfigure}{0.23\textwidth}
		\centering
		\begin{tikzpicture}
			\fill[color=white!70!red] (0,2) -- (0.7,1.3) to[out=45,in=-165] (2,2);
			\fill[color=white!70!red] (0,0) to[out=75,in=-128] (0.45,1.02) to[out=-52,in=165] (2,0);

			\draw[red] (2,0) to[out=165,in=-75] (0,2);
			\draw[red] (0,0) to[out=25,in=-90] (1.2,1) to[out=90,in=-25] (0,2);

			\draw[blue] (0,0) -- (2,2);
			\draw[blue] (2,0) -- (1,1);

			\draw[dashed] (0,0) -- (2,0) -- (2,2) -- (0,2) -- (0,0);
			\draw[fill] (0,0) circle (0.05) 
			(2,2) circle (0.05)
			(0,2) circle (0.05)
			(2,0) circle (0.05)
			(1,1) circle (0.05);

			\draw (0,2) -- (1,1);
			\draw (0,0) to[out=75,in=-135] (0.7,1.3) to[out=45,in=-165] (2,2);
		\end{tikzpicture}
		\caption{}
		\label{fig:nl l step 3 second a}
	\end{subfigure}
	\begin{subfigure}{0.23\textwidth}
		\centering
		\begin{tikzpicture}
			\draw[green] (0,0) -- (2,0) (0,2) -- (2,2);

			\draw[red] (2,0) to[out=165,in=-75] (0,2);
			\draw[red] (0,0) to[out=25,in=-90] (1.2,1) to[out=90,in=-25] (0,2);

			\draw[blue] (0,0) -- (2,2);
			\draw[blue] (2,0) -- (1,1);

			\draw[dashed] (2,0) -- (2,2) (0,2) -- (0,0);
			\draw[fill] (0,0) circle (0.05) 
			(2,2) circle (0.05)
			(0,2) circle (0.05)
			(2,0) circle (0.05)
			(1,1) circle (0.05);

			\draw (0,2) -- (1,1);
			\draw (0,0) to[out=75,in=-135] (0.7,1.3) to[out=45,in=-165] (2,2);

			\node[green] at (1,0.15) {\( a \)};
		\end{tikzpicture}
		\caption{}
		\label{fig:nl l step 3 second b}
	\end{subfigure}
	\begin{subfigure}{0.23\textwidth}
		\centering
		\begin{tikzpicture}
			\draw[green] (1.3,1.5) -- (1.2,2) (1.2,0) -- (1.4,0.4);

			\draw[blue] (0,0) -- (2,0) (0,2) -- (2,2);

			\draw[red] (2,0) to[out=165,in=-75] (0,2);
			\draw[red] (0,0) to[out=25,in=-90] (1.2,1) to[out=90,in=-25] (0,2);

			\draw[blue] (0,0) -- (2,2);
			\draw[blue] (2,0) -- (1,1);

			\draw[dashed] (2,0) -- (2,2) (0,2) -- (0,0);
			\draw[fill] (0,0) circle (0.05) 
			(2,2) circle (0.05)
			(0,2) circle (0.05)
			(2,0) circle (0.05)
			(1,1) circle (0.05);

			\draw (0,2) -- (1,1);
			\draw (0,0) to[out=75,in=-135] (0.7,1.3) to[out=45,in=-165] (2,2);
		\end{tikzpicture}
		\caption{}
		\label{fig:nl l step 3 second c}
	\end{subfigure}
	\begin{subfigure}{0.23\textwidth}
		\centering
		\begin{tikzpicture}
			\draw[green] (2,0) -- (1.2,2) (1.2,0) -- (1,1);

			\draw[blue] (0,0) -- (2,0) (0,2) -- (2,2);

			\draw[red] (2,0) to[out=165,in=-75] (0,2);
			\draw[red] (0,0) to[out=25,in=-90] (1.2,1) to[out=90,in=-25] (0,2);

			\draw[blue] (0,0) -- (2,2);
			\draw[blue] (2,0) -- (1,1);

			\draw[dashed] (2,0) -- (2,2) (0,2) -- (0,0);
			\draw[fill] (0,0) circle (0.05) 
			(2,2) circle (0.05)
			(0,2) circle (0.05)
			(2,0) circle (0.05)
			(1,1) circle (0.05);

			\draw (0,2) -- (1,1);
			\draw (0,0) to[out=75,in=-135] (0.7,1.3) to[out=45,in=-165] (2,2);
		\end{tikzpicture}
		\caption{}
		\label{fig:nl l step 3 second d}
	\end{subfigure}
	\begin{subfigure}{0.23\textwidth}
		\centering
		\begin{tikzpicture}
			\draw[green] (1,1) -- (1.4,2) (1.4,0) -- (2,2);

			\draw[blue] (0,0) -- (2,0) (0,2) -- (2,2);

			\draw[red] (2,0) to[out=165,in=-75] (0,2);
			\draw[red] (0,0) to[out=25,in=-90] (1.2,1) to[out=90,in=-25] (0,2);

			\draw[blue] (0,0) -- (2,2);
			\draw[blue] (2,0) -- (1,1);

			\draw[dashed] (2,0) -- (2,2) (0,2) -- (0,0);
			\draw[fill] (0,0) circle (0.05) 
			(2,2) circle (0.05)
			(0,2) circle (0.05)
			(2,0) circle (0.05)
			(1,1) circle (0.05);

			\draw (0,2) -- (1,1);
			\draw (0,0) to[out=75,in=-135] (0.7,1.3) to[out=45,in=-165] (2,2);
		\end{tikzpicture}
		\caption{}
		\label{fig:nl l step 3 second e}
	\end{subfigure}
	\begin{subfigure}{0.23\textwidth}
		\centering
		\begin{tikzpicture}
			\draw[red] (1,1) -- (1.4,2) (1.4,0) -- (2,2);
			\draw[red] (0,2) to[out=-65,in=180] (1,0.8) to[out=0,in=-115] (2,2);

			\draw[blue] (0,0) -- (2,0) (0,2) -- (2,2);

			\draw[red] (2,0) to[out=165,in=-75] (0,2);
			\draw[red] (0,0) to[out=25,in=-90] (1.2,1) to[out=90,in=-25] (0,2);

			\draw[blue] (0,0) -- (2,2);
			\draw[blue] (2,0) -- (1,1);

			\draw[dashed] (2,0) -- (2,2) (0,2) -- (0,0);
			\draw[fill] (0,0) circle (0.05) 
			(2,2) circle (0.05)
			(0,2) circle (0.05)
			(2,0) circle (0.05)
			(1,1) circle (0.05);

			\draw (0,2) -- (1,1);
			\draw (0,0) to[out=75,in=-135] (0.7,1.3) to[out=45,in=-165] (2,2);
		\end{tikzpicture}
		\caption{}
		\label{fig:nl l step 3 second f}
	\end{subfigure}
	\caption{}
	\label{fig:nl l step 3 second}
\end{figure}
By assumption, \( u \) is minimally intersected, so there must be some arc \( a' \in \mathcal{A} \) which intersects \( a \) but not \( u \). Applying \Cref{lem:twos} to the disk shown in red in \Cref{fig:nl l step 3 second a}, we find a path which must be a subarc of \( a' \). This is shown in green in \Cref{fig:nl l step 3 second c}. Then, we see that \( a' \) must be the arc shown in green either in \Cref{fig:nl l step 3 second d} or in \Cref{fig:nl l step 3 second e}.

We see that the arc shown in green in \Cref{fig:nl l step 3 second d} is in fact \( x'_1 \). So, up to a reflection switching \( w \) and \( x \), by Step 1 we assume that \( c' \) is not this arc.

So, \( c' \) must be the arc shown in \Cref{fig:nl l step 3 second e}. After a rotation of 90 degrees, we see that \( \mathcal{A} \) contains the arcs shown in \Cref{fig:nl l step 1 a}. Therefore we assume that \( w \) is not minimally intersected. 

There must be at least one other arc in \( \mathcal{A} \) which intersects \( w \) but not \( u \), in addition to \( w'_2 \). We see that \( w'_1 \) intersects twice \( c' \), and \( w'_4 \) intersects twice \( w'_2 \). So, it must be that \( w'_3 \in \mathcal{A} \). This is shown in \Cref{fig:nl l step 3 second f}. 

Now we have shown that \( x'_3, w'_3 \in \mathcal{A} \). So, we are done by Step 3.

\textbf{Step 5.} In this step we assume \( w_3', x_4' \in \mathcal{A} \). This is shown in \Cref{fig:nl l step 5 first a}. 

\begin{figure}
	\centering
	\begin{subfigure}{0.3\textwidth}
		\centering
		\begin{tikzpicture}
			\draw[red] (0,0) to[out=25,in=-110] (1.2,0.9) to[out=70,in=-90] (1.4,2) (1.4,0) -- (2,2);
			\draw[red] (2,2) to[out=-115,in=0] (1,0.8) to[out=180,in=-65] (0,2);

			\draw[blue] (0,0) -- (2,2);
			\draw[dashed] (0,0) -- (2,0) -- (2,2) -- (0,2) -- (0,0);
			\draw[fill] (0,0) circle (0.05) 
			(2,2) circle (0.05)
			(0,2) circle (0.05)
			(2,0) circle (0.05)
			(1,1) circle (0.05);

			\draw (0,2) -- (1,1);
			\draw (0,0) to[out=75,in=-135] (0.7,1.3) to[out=45,in=-165] (2,2);

		\end{tikzpicture}
		\caption{}
		\label{fig:nl l step 5 first a}
	\end{subfigure}
	\begin{subfigure}{0.3\textwidth}
		\centering
		\begin{tikzpicture}
			\draw[red] (0,0) to[out=25,in=-115] (2,2);
			\draw[red] (0,2) to[out=-65,in=155] (2,0);

			\draw[blue] (0,0) -- (1,1) -- (2,0);
			\draw[dashed] (0,0) -- (2,0) -- (2,2) -- (0,2) -- (0,0);
			\draw[fill] (0,0) circle (0.05) 
			(2,2) circle (0.05)
			(0,2) circle (0.05)
			(2,0) circle (0.05)
			(1,1) circle (0.05);

			\draw (0,2) -- (1,1);
			\draw (0,0) to[out=65,in=-180] (1,1.2) to[out=0,in=115] (2,0);

		\end{tikzpicture}
		\caption{}
		\label{fig:nl l step 5 first b}
	\end{subfigure}
	\caption{}
	\label{fig:nl l step 5 first}
\end{figure}


After applying a Dehn twist, we see that this system is equivalent to a subset of the system obtained in step 1. See \Cref{fig:nl l step 5 first b}. So, we may assume that the arc \( w \) is not minimally intersected. 

However, we see that \( w_1', w_2',w_4' \) intersect twice \( x'_4 \). So, we do not obtain a new 1-system in this step.

\subsection{\( u \) is a loop}

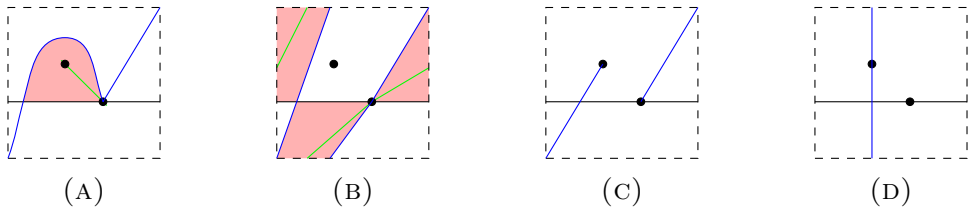
\begin{figure}
    \centering
    \begin{subfigure}{0.2\textwidth}
	    \centering
	    \begin{tikzpicture}
		\fill[color=white!70!red] (1.25,0.75) -- (0.2,0.75) to[out=75,in=180] (0.75,1.6) to[out=0,in=110] (1.25,0.75);
	    
		\draw[dashed] (0,0) -- (2,0) -- (2,2) -- (0,2) -- cycle;

		\draw[green] (0.7,1.3) -- (1.3,0.7);
		
		\draw[fill] (0.75,1.25) circle (0.05) (1.25,0.75) circle (0.05);
		
		\draw (0,0.75) -- (2,0.75);

		\draw[blue] (1.25,0.75) -- (2,2) (0,0) to[out=70,in=-105] (0.2,0.75) to[out=75,in=180] (0.75,1.6) to[out=0,in=110] (1.25,0.75);
	    \end{tikzpicture}
	    \caption{}
	    \label{fig:ack1}
    \end{subfigure}
    \begin{subfigure}{0.2\textwidth}
	    \centering
	    \begin{tikzpicture}
		\fill[color=white!70!red] (1.25,0.75) -- (2,0.75) -- (2,2) -- (1.25,0.75);
		
		\fill[color=white!70!red] (0,0.75) -- (0.2625,0.75) -- (0.7,2) -- (0,2) -- (0,0.75);
		
		\fill[color=white!70!red] (0,0) -- (0.7,0) -- (1.25,0.75) -- (0.2625,0.75) -- (0,0);

		\draw[green] (1.25,0.75) -- (2,1.2) (0,1.2) -- (0.4,2) (0.4,0) -- (1.25,0.75);

		\draw[dashed] (0,0) -- (2,0) -- (2,2) -- (0,2) -- cycle;
		
		\draw[fill] (0.75,1.25) circle (0.05) (1.25,0.75) circle (0.05);
		
		\draw (0,0.75) -- (2,0.75);
		
		\draw[blue] (1.25,0.75) -- (2,2) (0,0) -- (0.7,2) (0.7,0) -- (1.25,0.75);
	    \end{tikzpicture}
	    \caption{}
	    \label{fig:ack2}
    \end{subfigure}
    \begin{subfigure}{0.2\textwidth}
	    \centering
	    \begin{tikzpicture}
		\draw[dashed] (0,0) -- (2,0) -- (2,2) -- (0,2) -- cycle;
		
		\draw[fill] (0.75,1.25) circle (0.05) (1.25,0.75) circle (0.05);
		
		\draw (0,0.75) -- (2,0.75);
		
		\draw[blue] (1.25,0.75) -- (2,2) (0,0) -- (0.75,1.25);
	    \end{tikzpicture}
	    \caption{}
	    \label{fig:ack3}
    \end{subfigure}
    \begin{subfigure}{0.2\textwidth}
	    \centering
	    \begin{tikzpicture}
		\draw[dashed] (0,0) -- (2,0) -- (2,2) -- (0,2) -- cycle;
		
		\draw[fill] (0.75,1.25) circle (0.05) (1.25,0.75) circle (0.05);
		
		\draw (0,0.75) -- (2,0.75);
		
		\draw[blue] (0.75,0) -- (0.75,2);
	    \end{tikzpicture}
	    \caption{}
	    \label{fig:ack4}
    \end{subfigure}
    \caption{drawings of the torus}
    \label{fig:ack}
\end{figure}
Now, suppose $u$ is a loop, and assume furthermore that no non-loop arc is minimally intersected. Let $a$ denote the marked point at which $u$ is based and let $b$ denote the other marked point. Since $J = \emptyset$, there must be an arc $v \in \mathcal{A}$ that intersects $u$. There are four cases up to homeomorphism, as described in Figure  \ref{fig:ack}.


\textbf{Case 1.} 
In this case, we derive a contradiction. We apply Corollary~\ref{cor:twoss} to the disk shown in red in \Cref{fig:ack1} and we conclude that another arc $w$, shown in green in \Cref{fig:ack1}, must be included in $\mathcal{A}$. Furthermore, applying Lemma~\ref{lem:twos}, we see that any arc which intersects $w$ must also intersect $u$. However, this contradicts the assumption that $u$ is minimally intersected. Therefore, we can exclude the case shown in \Cref{fig:ack1}.


\textbf{Case 2.} In this case, we also derive a contradiction. We apply Corollary~\ref{cor:twoss} to the disk shown in red in \Cref{fig:ack2} and we conclude that another arc \( w \), shown in green in \Cref{fig:ack2}, must be included in $\mathcal{A}$. Then, by applying Lemma~\ref{lem:twos}, we see that any arc which intersects $w$ must also intersect $u$. This contradicts the assumption that $u$ is intersected minimally. Therefore, we can exclude the second case shown in the figure.


\textbf{Case 3.} In this case, it takes us a little longer to derive a contradiction. 

\begin{figure}
	\centering
	\begin{subfigure}{0.24\textwidth}
		\centering
		\begin{tikzpicture}
			\draw[dashed] (0,0) -- (2,0) -- (2,2) -- (0,2) -- (0,0);
			\draw[fill] (0,0) circle (0.05) 
			(2,2) circle (0.05)
			(0,2) circle (0.05)
			(2,0) circle (0.05)
			(1,1) circle (0.05);

			\draw[blue] (0,2) -- (1,1);
			\draw (0,0) to[out=75,in=-135] (0.7,1.3) to[out=45,in=-165] (2,2);

			\node[blue] at (0.5,1.7) {\( v \)};
			\node at (0.5,0.6) {\( u \)};
		\end{tikzpicture}
		\caption{}
		\label{fig:J013 intro a}
	\end{subfigure}
	\begin{subfigure}{0.24\textwidth}
		\centering
		\begin{tikzpicture}
			\draw[blue] (0,0) -- (1,1);
			\draw[dashed] (0,0) -- (2,0) -- (2,2) -- (0,2) -- (0,0);
			\draw[fill] (0,0) circle (0.05) 
			(2,2) circle (0.05)
			(0,2) circle (0.05)
			(2,0) circle (0.05)
			(1,1) circle (0.05);

			\draw (0,2) -- (1,1);
			\draw (0,0) to[out=75,in=-135] (0.7,1.3) to[out=45,in=-165] (2,2);
		\end{tikzpicture}
		\caption{}
		\label{fig:J013 intro b}
	\end{subfigure}
	\begin{subfigure}{0.24\textwidth}
		\centering
		\begin{tikzpicture}
			\draw[blue] (0,0) -- (1,1);
			\draw[red] (0,1.3) to[out=-45,in=160] (2,0) to[out=145,in=-40] (1,0.7) to[out=140,in=-150] (1,1.3) to[out=30,in=160] (2,1.3);
			\draw[dashed] (0,0) -- (2,0) -- (2,2) -- (0,2) -- (0,0);

			\draw[fill] (0,0) circle (0.05) 
			(2,2) circle (0.05)
			(0,2) circle (0.05)
			(2,0) circle (0.05)
			(1,1) circle (0.05);

			\draw (0,2) -- (1,1);
			\draw (0,0) to[out=75,in=-135] (0.7,1.3) to[out=45,in=-165] (2,2);
		\end{tikzpicture}
		\caption{}
		\label{fig:J013 intro c}
	\end{subfigure}
	\begin{subfigure}{0.24\textwidth}
		\centering
		\begin{tikzpicture}
			\draw[blue] (0,0) -- (1,1);
			\draw[red] (2,2) to[out=-155,in=150] (0.9,0.75) to[out=-30,in=-160] (2,0.7) (0,0.7) -- (2,0); 
			\draw[dashed] (0,0) -- (2,0) -- (2,2) -- (0,2) -- (0,0);

			\draw[fill] (0,0) circle (0.05) 
			(2,2) circle (0.05)
			(0,2) circle (0.05)
			(2,0) circle (0.05)
			(1,1) circle (0.05);

			\draw (0,2) -- (1,1);
			\draw (0,0) to[out=75,in=-135] (0.7,1.3) to[out=45,in=-165] (2,2);
		\end{tikzpicture}
		\caption{}
		\label{fig:J013 intro d}
	\end{subfigure}

	\begin{subfigure}{0.24\textwidth}
		\centering
		\begin{tikzpicture}
			\fill[color=white!70!red] (2,0) to[out=145,in=-40] (1,0.7) to[out=140,in=-150] (1,1.3) to[out=30,in=160] (2,1.3);
			\fill[color=white!70!red] (0,0) to[out=75,in=-120] (0.38,0.93) to[out=135,in=-45] (0,1.3);

			\draw[blue] (2,0) -- (1,1);

			\draw[red] (0,1.3) to[out=-45,in=160] (2,0) to[out=145,in=-40] (1,0.7) to[out=140,in=-150] (1,1.3) to[out=30,in=160] (2,1.3);
			\draw[dashed] (2,0) -- (2,2) (0,2) -- (0,0);
			\draw[dashed] (0,0) -- (2,0) (0,2) -- (2,2);

			\draw[fill] (0,0) circle (0.05) 
			(2,2) circle (0.05)
			(0,2) circle (0.05)
			(2,0) circle (0.05)
			(1,1) circle (0.05);

			\draw (0,2) -- (1,1);
			\draw (0,0) to[out=75,in=-135] (0.7,1.3) to[out=45,in=-165] (2,2);
		\end{tikzpicture}
		\caption{}
		\label{fig:J013 intro e}
	\end{subfigure}
	\begin{subfigure}{0.24\textwidth}
		\centering
		\begin{tikzpicture}
			\fill[color=white!70!red] (0,0) to[out=75,in=-117] (0.21,0.62) -- (0,0.7);
			\fill[color=white!70!red] (2,2) to[out=-155,in=150] (0.9,0.75) to[out=-30,in=-160] (2,0.7) -- (2,0) -- (0.21,0.62) to[out=63,in=-165] (2,2);

			\draw[blue] (2,0) to[out=155,in=-50] (0.8,0.7) to[out=130,in=-130] (0.8,1.3) to[out=50,in=-155] (2,2);

			\draw[red] (2,2) to[out=-155,in=150] (0.9,0.75) to[out=-30,in=-160] (2,0.7) (0,0.7) -- (2,0); 
			\draw[dashed] (2,0) -- (2,2) (0,2) -- (0,0);
			\draw[dashed] (0,0) -- (2,0) (0,2) -- (2,2);

			\draw[fill] (0,0) circle (0.05) 
			(2,2) circle (0.05)
			(0,2) circle (0.05)
			(2,0) circle (0.05)
			(1,1) circle (0.05);

			\draw (0,2) -- (1,1);
			\draw (0,0) to[out=75,in=-135] (0.7,1.3) to[out=45,in=-165] (2,2);
		\end{tikzpicture}
		\caption{}
		\label{fig:J013 intro f}
	\end{subfigure}
	\caption{}
	\label{fig:J013 intro}
\end{figure}

Let $w$ be the arc shown in blue in \Cref{fig:J013 intro b}. 

\begin{claim}
	\( w \in \mathcal{A} \)
\end{claim}

\begin{proof}

	Suppose for contradiction that \( w \notin \mathcal{A} \). Since \( \mathcal{A} \) is saturated, there must be some arc \( w' \in \mathcal{A} \) which intersects \( w \) at least twice.

	We appeal to the proof of \Cref{claim:sec5 some arcs}. We did not invoke the minimality in the first part of the proof, so we again find that \( w' \) must be either the arc shown in red in \Cref{fig:J013 intro c} or the arc shown in red in \Cref{fig:J013 intro d}. 

	If \( w' \) is the arc shown in \Cref{fig:J013 intro c}, we apply \Cref{cor:twoss} to the disk shown in red in \Cref{fig:J013 intro e}. We see that the arc shown in blue must be included in \( \mathcal{A} \). Then, since the boundary of this disk is contained in the arcs \( u, w' \) we apply \Cref{lem:twos} and it must be that any arc which intersects the arc shown in blue must be included in \( \mathcal{A} \). This contradicts the assumption that \( u \) is minimally intersected.

	Similarly, if \( w' \) is the arc shown in \Cref{fig:J013 intro d}, we apply \Cref{cor:twoss} and then \Cref{lem:twos} to the disk shown in red in \Cref{fig:J013 intro f} and we find the the arc shown in blue must be in \( \mathcal{A} \) and contradicts the minimality of \( u \). 

	Since there is no arc \( w' \in \mathcal{A} \) which intersects \( w \) twice, it must be that \( w \in \mathcal{A} \). 
\end{proof}

Now, since we have assumed $u$ is intersected minimally, there must be some arc $w' \in \mathcal{A}$ which intersects $w$ but not $u$. We may apply \Cref{lem:twos} to the disk shown in red in \Cref{fig:J013 third a} and we determine that the path shown in green must be a subarc of \( w' \). There are three arcs which have this path as a subarc, shown in red in \Cref{fig:J013 third b,fig:J013 third c,fig:J013 third d}.

\begin{figure}
	\centering
	\begin{subfigure}{0.23\textwidth}
		\centering
		\begin{tikzpicture}
			\fill[color=white!70!red] (1,1) to[out=-100,in=25] (0,0) to[out=75,in=-135] (0.7,1.3);
			\draw[blue] (0,0) -- (1,1);
			\draw[dashed] (0,0) -- (2,0) -- (2,2) -- (0,2) -- (0,0);

			\draw[green] (1,0.7) to[out=180,in=180] (1,1.2);
			\draw[fill] (0,0) circle (0.05) 
			(2,2) circle (0.05)
			(0,2) circle (0.05)
			(2,0) circle (0.05)
			(1,1) circle (0.05);

			\draw (0,2) -- (1,1);
			\draw (0,0) to[out=75,in=-135] (0.7,1.3) to[out=45,in=-165] (2,2);
		\end{tikzpicture}
		\caption{}
		\label{fig:J013 third a}
	\end{subfigure}
	\begin{subfigure}{0.23\textwidth}
		\centering
		\begin{tikzpicture}
			\draw[red] (2,0) to[out=150,in=-90] (0.8,1) to[out=90,in=-150] (2,2);
			\draw[blue] (0,0) -- (1,1);
			\draw[dashed] (0,0) -- (2,0) -- (2,2) -- (0,2) -- (0,0);
			\draw[fill] (0,0) circle (0.05) 
			(2,2) circle (0.05)
			(0,2) circle (0.05)
			(2,0) circle (0.05)
			(1,1) circle (0.05);

			\draw (0,2) -- (1,1);
			\draw (0,0) to[out=75,in=-135] (0.7,1.3) to[out=45,in=-165] (2,2);
		\end{tikzpicture}
		\caption{}
		\label{fig:J013 third b}
	\end{subfigure}
	\begin{subfigure}{0.23\textwidth}
		\centering
		\begin{tikzpicture}
			\draw[red] (2,2) to[out=-120,in=-45] (0.85,0.85) to[out=135,in=-150] (2,2);

			\draw[blue] (0,0) -- (1,1);
			\draw[dashed] (0,0) -- (2,0) -- (2,2) -- (0,2) -- (0,0);
			\draw[fill] (0,0) circle (0.05) 
			(2,2) circle (0.05)
			(0,2) circle (0.05)
			(2,0) circle (0.05)
			(1,1) circle (0.05);

			\draw (0,2) -- (1,1);
			\draw (0,0) to[out=75,in=-135] (0.7,1.3) to[out=45,in=-165] (2,2);
		\end{tikzpicture}
		\caption{}
		\label{fig:J013 third c}
	\end{subfigure}
	\begin{subfigure}{0.23\textwidth}
		\centering
		\begin{tikzpicture}
			\draw[red] (2,0) to[out=120,in=45] (0.85,1.15) to[out=-135,in=150] (2,0);

			\draw[blue] (0,0) -- (1,1);
			\draw[dashed] (0,0) -- (2,0) -- (2,2) -- (0,2) -- (0,0);
			\draw[fill] (0,0) circle (0.05) 
			(2,2) circle (0.05)
			(0,2) circle (0.05)
			(2,0) circle (0.05)
			(1,1) circle (0.05);

			\draw (0,2) -- (1,1);
			\draw (0,0) to[out=75,in=-135] (0.7,1.3) to[out=45,in=-165] (2,2);
		\end{tikzpicture}
		\caption{}
		\label{fig:J013 third d}
	\end{subfigure}
	\caption{}
	\label{fig:J013 third}
\end{figure}

The arcs in \Cref{fig:J013 third c,fig:J013 third d} both bound monogons, so if either of these arcs were included in \( \mathcal{A} \), we would apply \Cref{cor:help ones}, and we would conclude that \( J \ne \emptyset \). This is a contradiction. So, the only arc that could be in \( \mathcal{A} \) which intersects \( w \) but not \( u \) is the arc shown in red in \Cref{fig:J013 third b}.

However, this implies that \( u \) is intersected by at least as many arcs of \( \mathcal{A} \) as \( w \). Since by assumption \( u \) is minimally intersected, this means that \( w \) is also minimally intersected. Note that \( w \) is a non-loop arc, giving a contradiction.

\textbf{Case 4.}
In this case, we again put in some work to derive a contradiction. 
\begin{figure}
	\centering
	\begin{subfigure}{0.18\textwidth}
		\centering
		\begin{tikzpicture}

			\draw (0,0) -- (2,0) (2,2) -- (0,2); 
			\draw[dashed] (0,0) -- (0,2) (2,2) -- (2,0);
			\draw[blue] (1,0) -- (1,2);

			\draw[fill] (0,0) circle (0.06)
			(2,0) circle (0.06)
			(0,2) circle (0.06)
			(2,2) circle (0.06)
			(1,1) circle (0.06);

			\node[blue] at (1.2,1.5) {\( v \)};
			\node at (1.5,1.8) {\( u \)};
		\end{tikzpicture}
		\caption{}
		\label{fig:J014 intro a}
	\end{subfigure}
	\begin{subfigure}{0.18\textwidth}
		\centering
		\begin{tikzpicture}

			\draw[blue] (0,0) -- (2,2)
			(2,0) -- (0,2);

			\draw (0,0) -- (2,0) (2,2) -- (0,2); 
			\draw[dashed] (0,0) -- (0,2) (2,2) -- (2,0);
			\draw (1,0) -- (1,2);

			\draw[fill] (0,0) circle (0.06)
			(2,0) circle (0.06)
			(0,2) circle (0.06)
			(2,2) circle (0.06)
			(1,1) circle (0.06);

			\draw[fill,color=white] (0.5,0.5) circle (0.15)
			(1.5,0.5) circle (0.15)
			(1.5,1.5) circle (0.15)
			(0.5,1.5) circle (0.15);

			\node[blue] at (0.5,1.5) {$w$};
			\node[blue] at (1.5,1.5) {$x$};
			\node[blue] at (0.5,0.5) {$y$};
			\node[blue] at (1.5,0.5) {$z$};
		\end{tikzpicture}
		\caption{}
		\label{fig:J014 intro a5}
	\end{subfigure}
	\begin{subfigure}{0.18\textwidth}
		\centering
		\begin{tikzpicture}
			\fill[color=white!70!red] (1,1) arc(-90:-180:1) -- (1,2);

			\draw[blue] (1,1) -- (0,2);

			\draw[red] (0.2,1.1) -- (0.3,2) (0.3,0) -- (0.2,0.2)
			(0.5,1) -- (1.3,1.3);

			\draw (0,0) -- (2,0) (2,2) -- (0,2); 
			\draw[dashed] (0,0) -- (0,2) (2,2) -- (2,0);
			\draw (1,0) -- (1,2);

			\draw[fill] (0,0) circle (0.06)
			(2,0) circle (0.06)
			(0,2) circle (0.06)
			(2,2) circle (0.06)
			(1,1) circle (0.06);
		\end{tikzpicture}
		\caption{}
		\label{fig:J014 intro b}
	\end{subfigure}
	\begin{subfigure}{0.18\textwidth}
		\centering
		\begin{tikzpicture}
			\draw[blue] (1,1) -- (0,2);

			\draw[red] (1,1) -- (2,1.3) (0,1.3) -- (2,1.7) (0,1.7) -- (0.3,2) (0.3,0) -- (0,0.3) (2,0.3) -- (1,1);
			\draw (0,0) -- (2,0) (2,2) -- (0,2); 
			\draw[dashed] (0,0) -- (0,2) (2,2) -- (2,0);
			\draw (1,0) -- (1,2);

			\draw[fill] (0,0) circle (0.06)
			(2,0) circle (0.06)
			(0,2) circle (0.06)
			(2,2) circle (0.06)
			(1,1) circle (0.06);
		\end{tikzpicture}
		\caption{}
		\label{fig:J014 intro c}
	\end{subfigure}
	\begin{subfigure}{0.23\textwidth}
		\centering
		\begin{tikzpicture}
			\draw[blue] (1,1) -- (0,2);

			\draw[red] (1,1) -- (2,1.3) (0,1.3) -- (2,1.7) (0,1.7) -- (0.3,2) (0.3,0) -- (1,1);
			\draw (0,0) -- (2,0) (2,2) -- (0,2); 
			\draw[dashed] (0,0) -- (0,2) (2,2) -- (2,0);
			\draw (1,0) -- (1,2);

			\draw[fill] (0,0) circle (0.06)
			(2,0) circle (0.06)
			(0,2) circle (0.06)
			(2,2) circle (0.06)
			(1,1) circle (0.06);
		\end{tikzpicture}
		\caption{}
		\label{fig:J014 intro f}
	\end{subfigure}
	\caption{}
	\label{fig:J014 intro}
\end{figure}
Let $w, x, y, z$ be the arcs as shown in \Cref{fig:J014 intro a5}.

\begin{claim}
	$w, x, y, z \in \mathcal{A}$.
\end{claim}
\begin{proof}
	Suppose for contradiction that \( w \notin \mathcal{A} \). Since \( \mathcal{A} \) is saturated, there must be some arc \( w' \in \mathcal{A} \) which intersects \( w \) twice. 

	We apply \Cref{lem:twos} to the disk shown in red in \Cref{fig:J014 intro b}. We determine two paths which must be subarcs of \( w' \). 

	We see that \( w' \) intersects \( u \). Therefore, from Case 1, Case 2, and Case 3 above, we may assume that \( w' \) is a loop based at the same marked point as \( v \). So, from this we see that \( w' \) must be either the arc shown in red in \Cref{fig:J014 intro c} or \Cref{fig:J014 intro f}. 

	First suppose \( w' \) is the arc shown in \Cref{fig:J014 intro c}. We apply \Cref{cor:twoss} to the disk shown in red in \Cref{fig:J014 intro d}, and we conclude that the arc \( z \), shown in blue in \Cref{fig:J014 intro d}, must be in \( \mathcal{A} \). Then, by assumption, \( u \) is minimally intersected, so there must be some arc \( z' \in \mathcal{A} \) which intersects \( z \) but not \( u \). We apply \Cref{lem:twos} and we conclude that the path shown in green in \Cref{fig:J014 intro e} must be a subarc of \( z' \). However, we see that this path cannot be extended to an arc without intersecting \( u \) or intersecting twice some arc of \( \mathcal{A} \). So, \( z \) contradicts the assumption that \( u \) is minimally intersected. 

\begin{figure}
	\centering
	\begin{subfigure}{0.23\textwidth}
		\centering
		\begin{tikzpicture}
			\fill[color=white!70!red] (1,0) -- (1,1) -- (2,0.3) -- (2,0);
			\fill[color=white!70!red] (0,0) -- (0,0.3) -- (0.3,0);

			\draw[blue] (1,1) -- (2,0);

			\draw[red] (1,1) -- (2,1.3) (0,1.3) -- (2,1.7) (0,1.7) -- (0.3,2) (0.3,0) -- (0,0.3) (2,0.3) -- (1,1);
			\draw (0,0) -- (2,0) (2,2) -- (0,2); 
			\draw[dashed] (0,0) -- (0,2) (2,2) -- (2,0);
			\draw (1,0) -- (1,2);

			\draw[fill] (0,0) circle (0.06)
			(2,0) circle (0.06)
			(0,2) circle (0.06)
			(2,2) circle (0.06)
			(1,1) circle (0.06);
		\end{tikzpicture}
		\caption{}
		\label{fig:J014 intro d}
	\end{subfigure}
	\begin{subfigure}{0.18\textwidth}
		\centering
		\begin{tikzpicture}
			\draw[green] (0.8,0.7) -- (1.3,1);

			\draw[blue] (1,1) -- (2,0);

			\draw[red] (1,1) -- (2,1.3) (0,1.3) -- (2,1.7) (0,1.7) -- (0.3,2) (0.3,0) -- (0,0.3) (2,0.3) -- (1,1);
			\draw (0,0) -- (2,0) (2,2) -- (0,2); 
			\draw[dashed] (0,0) -- (0,2) (2,2) -- (2,0);
			\draw (1,0) -- (1,2);

			\draw[fill] (0,0) circle (0.06)
			(2,0) circle (0.06)
			(0,2) circle (0.06)
			(2,2) circle (0.06)
			(1,1) circle (0.06);
		\end{tikzpicture}
		\caption{}
		\label{fig:J014 intro e}
	\end{subfigure}
	\caption{}
	\label{fig:J014 wp1}
\end{figure}

	So, \( w' \) must be the arc shown in \Cref{fig:J014 intro f}. We apply \Cref{cor:twoss} to the disk shown in red in \Cref{fig:J014 intro g} and we conclude that another arc \( a \), shown in blue in \Cref{fig:J014 intro g}, must be in \( \mathcal{A} \). 

\begin{figure}
	\centering
	\begin{subfigure}{0.23\textwidth}
		\centering
		\begin{tikzpicture}
			\draw[fill,color=white!70!red] (0,1.3) -- (0,1.7) -- (0.3,2) -- (1,2) -- (1,0) -- (0.3,0) -- (1,1) -- (2,1.3) -- (2,1.7);

			\draw[blue] (1,1) -- (2,1.5) (0,1.5) -- (0.6,2) (0.6,0) -- (1,1);

			\draw[red] (1,1) -- (2,1.3) (0,1.3) -- (2,1.7) (0,1.7) -- (0.3,2) (0.3,0) -- (1,1);
			\draw (0,0) -- (2,0) (2,2) -- (0,2); 
			\draw[dashed] (0,0) -- (0,2) (2,2) -- (2,0);
			\draw (1,0) -- (1,2);

			\draw[fill] (0,0) circle (0.06)
			(2,0) circle (0.06)
			(0,2) circle (0.06)
			(2,2) circle (0.06)
			(1,1) circle (0.06);
		\end{tikzpicture}
		\caption{}
		\label{fig:J014 intro g}
	\end{subfigure}
	\begin{subfigure}{0.23\textwidth}
		\centering
		\begin{tikzpicture}
			\draw[green] (0,0) to[out=70,in=180] (1,1.2) to[out=0,in=110] (2,0); 

			\draw[blue] (1,1) -- (2,1.5) (0,1.5) -- (0.6,2) (0.6,0) -- (1,1);

			\draw[red] (1,1) -- (2,1.3) (0,1.3) -- (2,1.7) (0,1.7) -- (0.3,2) (0.3,0) -- (1,1);
			\draw (0,0) -- (2,0) (2,2) -- (0,2); 
			\draw[dashed] (0,0) -- (0,2) (2,2) -- (2,0);
			\draw (1,0) -- (1,2);

			\draw[fill] (0,0) circle (0.06)
			(2,0) circle (0.06)
			(0,2) circle (0.06)
			(2,2) circle (0.06)
			(1,1) circle (0.06);
		\end{tikzpicture}
		\caption{}
		\label{fig:J014 intro h}
	\end{subfigure}
	\begin{subfigure}{0.23\textwidth}
		\centering
		\begin{tikzpicture}
			\draw[green] (1,1) -- (2,0.6) (0,0.6) -- (2,0);

			\draw[blue] (1,1) -- (2,1.5) (0,1.5) -- (0.6,2) (0.6,0) -- (1,1);

			\draw[red] (1,1) -- (2,1.3) (0,1.3) -- (2,1.7) (0,1.7) -- (0.3,2) (0.3,0) -- (1,1);
			\draw (0,0) -- (2,0) (2,2) -- (0,2); 
			\draw[dashed] (0,0) -- (0,2) (2,2) -- (2,0);
			\draw (1,0) -- (1,2);

			\draw[fill] (0,0) circle (0.06)
			(2,0) circle (0.06)
			(0,2) circle (0.06)
			(2,2) circle (0.06)
			(1,1) circle (0.06);
		\end{tikzpicture}
		\caption{}
		\label{fig:J014 intro i}
	\end{subfigure}
	\begin{subfigure}{0.23\textwidth}
		\centering
		\begin{tikzpicture}
			\draw[green] (1,1) -- (0,0.7) (2,0.7) -- (0,0);

			\draw[blue] (1,1) -- (2,1.5) (0,1.5) -- (0.6,2) (0.6,0) -- (1,1);

			\draw[red] (1,1) -- (2,1.3) (0,1.3) -- (2,1.7) (0,1.7) -- (0.3,2) (0.3,0) -- (1,1);
			\draw (0,0) -- (2,0) (2,2) -- (0,2); 
			\draw[dashed] (0,0) -- (0,2) (2,2) -- (2,0);
			\draw (1,0) -- (1,2);

			\draw[fill] (0,0) circle (0.06)
			(2,0) circle (0.06)
			(0,2) circle (0.06)
			(2,2) circle (0.06)
			(1,1) circle (0.06);
		\end{tikzpicture}
		\caption{}
		\label{fig:J014 intro j}
	\end{subfigure}
	\caption{}
	\label{fig:J014 wp2}
\end{figure}

	In the system we have constructed, \( u \) intersects three other arcs, and \( a \) intersects one other arc. By assumption, \( u \) is minimally intersected, so there must be at least two arcs in \( \mathcal{A} \) which intersect \( a \) but not \( u \). Using \Cref{lem:twos}, we may determine that the only such arcs are the three arcs shown in green in \Cref{fig:J014 intro h,fig:J014 intro i,fig:J014 intro j}. However, these arcs each intersect pairwise twice. So, at most one of them can be in \( \mathcal{A} \). This contradicts the assumption that \( u \) is minimally intersected. 

	So, no such arc \( w' \) can be in \( \mathcal{A} \), and therefore \( w \in \mathcal{A} \). We conclude by symmetry that \( x, y, z \), as shown in \Cref{fig:J014 intro a5}, must be in \( \mathcal{A} \). 
\end{proof}

By the assumption that $u$ is minimally intersected, there must be some arc $w' \in \mathcal{A}$ which intersects $w$ but not $u$. Since \( w \) is a non-loop arc, by assumption \( w \) is not minimally intersected, so there must be at least two distinct arcs in \( \mathcal{A} \) which intersect \( w \) but not \( u \). 

There are three possibilities for $w'$. Let $w'_1, w'_2, w'_3$ be the arcs shown in \Cref{fig:J014 wp trues}.

\begin{figure}
	\centering
	\begin{subfigure}{0.23\textwidth}
		\centering
		\begin{tikzpicture}
			\node at (1.2,1.5) {\( v \)};
			\node at (1.5,1.8) {\( u \)};

			\draw[blue] (1,1) -- (0,2);

			\draw[red] (0,0) to[out=70,in=180] (1,1.2) to[out=0,in=110] (2,0);

			\draw (0,0) -- (2,0) (2,2) -- (0,2); 
			\draw[dashed] (0,0) -- (0,2) (2,2) -- (2,0);
			\draw (1,0) -- (1,2);

			\draw[fill] (0,0) circle (0.06)
			(2,0) circle (0.06)
			(0,2) circle (0.06)
			(2,2) circle (0.06)
			(1,1) circle (0.06);

			\draw[fill,white] (1.6,1) circle (0.2);
			\node[red] at (1.6,1) {$w'_1$};
		\end{tikzpicture}
		\caption{}
		\label{fig:J014 wp trues a}
	\end{subfigure}
	\begin{subfigure}{0.23\textwidth}
		\centering
		\begin{tikzpicture}

			\draw[blue] (1,1) -- (0,2);

			\draw[red] (0,0) to[out=70,in=-160] (2,2);

			\draw (0,0) -- (2,0) (2,2) -- (0,2); 
			\draw[dashed] (0,0) -- (0,2) (2,2) -- (2,0);
			\draw (1,0) -- (1,2);

			\draw[fill] (0,0) circle (0.06)
			(2,0) circle (0.06)
			(0,2) circle (0.06)
			(2,2) circle (0.06)
			(1,1) circle (0.06);

			\draw[fill,white] (0.4,1) circle (0.2);
			\node[red] at (0.4,1) {$w'_2$};
		\end{tikzpicture}
		\caption{}
		\label{fig:J014 wp trues b}
	\end{subfigure}
	\begin{subfigure}{0.23\textwidth}
		\centering
		\begin{tikzpicture}

			\draw[blue] (1,1) -- (0,2);

			\draw[red] (1,1) -- (2,1.3) (0,1.3) -- (2,2);

			\draw (0,0) -- (2,0) (2,2) -- (0,2); 
			\draw[dashed] (0,0) -- (0,2) (2,2) -- (2,0);
			\draw (1,0) -- (1,2);

			\draw[fill] (0,0) circle (0.06)
			(2,0) circle (0.06)
			(0,2) circle (0.06)
			(2,2) circle (0.06)
			(1,1) circle (0.06);

			\draw[fill,white] (1.5,1.15) circle (0.2);
			\node[red] at (1.5,1.15) {$w'_3$};
		\end{tikzpicture}
		\caption{}
		\label{fig:J014 wp trues c}
	\end{subfigure}
	\caption{}
	\label{fig:J014 wp trues}
\end{figure}

Write \( x'_i \) for the arc obtained from \( w'_i \) by reflection across the vertical axis. Write \( y'_i, z'_i \) for the arcs obtained from \( w'_i, x'_i \), respectively, by reflection across the horizontal axis. 

Note the following: $w'_1 = x'_1$ and $y'_1 = z'_1$, and $w'_1$ intersects $y'_1$ twice. Additionally, the following pairs intersect twice:

\( (w'_2, x'_3), (w'_3, x'_2), (w'_3, x'_3), (y'_2, z'_3), (y'_3, z'_2), (y'_3, z'_3) \).

\textbf{Step 1.} For this step, we assume $w'_1 \in \mathcal{A}$. This means $y'_1 \notin \mathcal{A}$. Since there must be two distinct arcs in \( \mathcal{A} \) which intersect \( y \) but not \( u \), we have \( y_2', y_3' \in \mathcal{A} \). However, since \( y_3' \) intersects twice each \( z_2' \) and \( z_3' \), we have \( z_1', z_2', z_3' \notin \mathcal{A} \). This contradicts the assumption that \( u \) is minimally intersected.

\textbf{Step 2.} For this step, we assume \( w_1' \notin \mathcal{A} \). Since there must be two arcs in \( \mathcal{A} \) which intersect \( w \) but not \( u \), it must be that \( w_2', w_3' \in \mathcal{A} \). However, since \( w_3' \) intersects twice each \( x_2' \) and \( x_3' \), we have \( x_1', x_2', x_3' \notin \mathcal{A} \). This contradicts the assumption that \( u \) is minimally intersected.

\clearpage

\section{Acknowledgements}

I would like to thank my supervisor Piotr Przytycki for everything he taught me, and repeatedly reading and offering corrections. 

Thank you the the Centre de Recherches Math\'ethiques for funding me to travel to Paris, where I was able to work with my supervisor in person.

Thank you to Patricia Sorya, Giacomo Bascape, and anyone who stood still and let me present about some or all of my project to them, and asked questions.

Thank you to Antoine Poulin for help with French.

Thank you to Zachary Feng for buying me coffee.

\clearpage

\bibliography{bibl/proj}

\begin{bibdiv}
\begin{biblist}

\bib{aougab_2019}{article}{
      author={Aougab, Tarik},
      author={Biringer, Ian},
      author={Gaster, Jonah},
       title={Packing curves on surfaces with few intersections},
        date={2019},
        ISSN={1073-7928},
     journal={Int. Math. Res. Not. IMRN},
      number={16},
       pages={5205\ndash 5217},
         url={https://doi.org/10.1093/imrn/rnx270},
}

\bib{farb_2012}{book}{
      author={Farb, Benson},
      author={Margalit, Dan},
       title={A primer on mapping class groups},
   publisher={Princeton University Press},
        date={2012},
}

\bib{greene_2019}{article}{
      author={Greene, Joshua~Evan},
       title={On loops intersecting at most once},
        date={2019},
        ISSN={1016-443X},
     journal={Geom. Funct. Anal.},
      volume={29},
      number={6},
       pages={1828\ndash 1843},
         url={https://doi.org/10.1007/s00039-019-00517-0},
}

\bib{juvan_1996}{article}{
      author={Juvan, M.},
      author={Malni\v{c}, A.},
      author={Mohar, B.},
       title={Systems of curves on surfaces},
        date={1996},
        ISSN={0095-8956},
     journal={J. Combin. Theory Ser. B},
      volume={68},
      number={1},
       pages={7\ndash 22},
         url={https://doi.org/10.1006/jctb.1996.0053},
}

\bib{przytycki_2015}{article}{
      author={Przytycki, Piotr},
       title={Arcs intersecting at most once},
        date={2015},
     journal={Geometric and Functional Analysis},
      volume={25},
      number={2},
       pages={658–670},
}

\bib{tee_2021}{article}{
      author={Tee, Paul},
       title={Arcs intersecting at most once on the 4-punctured sphere},
        date={2021},
      eprint={arXiv:2107.06462},
         url={https://arxiv.org/abs/2107.06462},
}

\end{biblist}
\end{bibdiv}


\begin{bibdiv}
\begin{biblist}

\bib{aougab_2019}{article}{
      author={Aougab, Tarik},
      author={Biringer, Ian},
      author={Gaster, Jonah},
       title={Packing curves on surfaces with few intersections},
        date={2019},
        ISSN={1073-7928},
     journal={Int. Math. Res. Not. IMRN},
      number={16},
       pages={5205\ndash 5217},
         url={https://doi.org/10.1093/imrn/rnx270},
}

\bib{greene_2019}{article}{
      author={Greene, Joshua~Evan},
       title={On loops intersecting at most once},
        date={2019},
        ISSN={1016-443X},
     journal={Geom. Funct. Anal.},
      volume={29},
      number={6},
       pages={1828\ndash 1843},
         url={https://doi.org/10.1007/s00039-019-00517-0},
}

\bib{juvan_1996}{article}{
      author={Juvan, M.},
      author={Malni\v{c}, A.},
      author={Mohar, B.},
       title={Systems of curves on surfaces},
        date={1996},
        ISSN={0095-8956},
     journal={J. Combin. Theory Ser. B},
      volume={68},
      number={1},
       pages={7\ndash 22},
         url={https://doi.org/10.1006/jctb.1996.0053},
}

\bib{przytycki_2015}{article}{
      author={Przytycki, Piotr},
       title={Arcs intersecting at most once},
        date={2015},
     journal={Geometric and Functional Analysis},
      volume={25},
      number={2},
       pages={658–670},
}

\bib{tee_2021}{article}{
      author={Tee, Paul},
       title={Arcs intersecting at most once on the 4-punctured sphere},
        date={2021},
      eprint={arXiv:2107.06462},
         url={https://arxiv.org/abs/2107.06462},
}

\end{biblist}
\end{bibdiv}


\begin{bibdiv}
\begin{biblist}

\bib{aougab_2019}{article}{
      author={Aougab, Tarik},
      author={Biringer, Ian},
      author={Gaster, Jonah},
       title={Packing curves on surfaces with few intersections},
        date={2019},
        ISSN={1073-7928},
     journal={Int. Math. Res. Not. IMRN},
      number={16},
       pages={5205\ndash 5217},
         url={https://doi.org/10.1093/imrn/rnx270},
      review={\MR{4001026}},
}

\bib{greene_2018}{misc}{
      author={Greene, Joshua~Evan},
       title={On curves intersecting at most once},
   publisher={arXiv},
        date={2018},
         url={https://arxiv.org/abs/1807.05658},
}

\bib{juvan_1996}{article}{
      author={Juvan, M.},
      author={Malni\v{c}, A.},
      author={Mohar, B.},
       title={Systems of curves on surfaces},
        date={1996},
        ISSN={0095-8956},
     journal={J. Combin. Theory Ser. B},
      volume={68},
      number={1},
       pages={7\ndash 22},
         url={https://doi.org/10.1006/jctb.1996.0053},
      review={\MR{1405702}},
}

\bib{przytycki_2015}{article}{
      author={Przytycki, Piotr},
       title={Arcs intersecting at most once},
        date={2015},
     journal={Geometric and Functional Analysis},
      volume={25},
      number={2},
       pages={658–670},
}

\bib{tee_2021}{misc}{
      author={Tee, Paul},
       title={Arcs intersecting at most once on the 4-punctured sphere},
   publisher={arXiv},
        date={2021},
         url={https://arxiv.org/abs/2107.06462},
}

\end{biblist}
\end{bibdiv}

\end{document}